\documentclass[12pt,leqno]{amsart}

\usepackage[margin=1in,lmargin=1in,rmargin=1in]{geometry}
\usepackage{amsmath,latexsym,amssymb,amsthm,graphicx}
\usepackage{mathtools}
\usepackage{subfig}
\usepackage{xfrac}
\usepackage{fix-cm}
\usepackage{marginnote}
\usepackage{color}
\usepackage{verbatim}

\usepackage{doi}
\usepackage{stmaryrd}
\usepackage{multirow}
\usepackage{float}
\usepackage{tabularx}
\usepackage[shortlabels]{enumitem}
\usepackage{pdfsync}
\usepackage{tikz}
\usepackage{tikz-cd}
\usepackage{float}
\usepackage[utf8]{inputenc}
\usepackage[textsize=small]{todonotes}

\usepackage{hyphenat}
\usepackage{booktabs}

\usepackage[onehalfspacing]{setspace}
\allowdisplaybreaks

\usepackage{latexsym,enumitem}
\usepackage{amssymb}
\usepackage{psfrag}
\usepackage{amsthm}
\usepackage{amscd}
\usepackage{amsmath}
\usepackage{amsfonts}
\usepackage{graphics,caption}
\usepackage[all]{xy}
\usepackage{etoolbox}
\patchcmd{\quote}{\rightmargin}{\leftmargin 2em \rightmargin}{}{}
\usepackage[sort,numbers]{natbib} 

\graphicspath{{mpdraws/}{draws/}}

\usepackage[notextcomp]{kpfonts}

\usepackage[english]{babel}
\usepackage{comment}
\usepackage{color}
\usepackage{xcolor}

\DeclareMathOperator{\tr}{tr} 
\DeclareMathOperator{\diag}{diag} 

\let\phi\varphi

\newcommand{\f}[2]{\frac{#1}{#2}} 

\let\epsilon\varepsilon
\let\subset\subseteq
\newcommand{\be}{\begin{equation*}}
 \newcommand{\ee}{\end{equation*}}
\newcommand{\bpf}{\begin{dimo}}
\newcommand{\epf}{\end{dimo}}
\newcommand{\bdefi}{\begin{defin}}
\newcommand{\edefi}{\end{defin}}
\newcommand{\bthm}{\begin{thm}}
\newcommand{\ethm}{\end{thm}}
\newcommand{\blem}{\begin{lem}}
\newcommand{\elem}{\end{lem}}
\newcommand{\bcor}{\begin{cor}}
\newcommand{\ecor}{\end{cor}}

\newcommand{\bprop}{\begin{prop}}
\newcommand{\eprop}{\end{prop}}
\newcommand{\bese}{\begin{ese}}
\newcommand{\eese}{\end{ese}}
\newcommand{\brem}{\begin{rem}}
\newcommand{\erem}{\end{rem}}
\newcommand{\bpfc}{\begin{dimoclaim}}
\newcommand{\epfc}{\end{dimoclaim}}

\newcommand{\rar}{\rightarrow} 

\newcommand{\diff}{\mathop{}\!\mathrm{d}} 
	
\newcommand{\abs}[1]{\left\lvert#1\right\rvert}						
\newcommand{\norm}[1]{\left\lVert#1\right\rVert}					
\newcommand{\set}[1]{\left\{#1\right\}}					
\newcommand{\ceiling}[1]{\left\lceil#1\right\rceil}					
\newcommand{\floor}[1]{\left\lfloor#1\right\rfloor}					
\newcommand{\tonde}[1]{\left(#1\right)}							
	
\DeclareMathOperator{\interior}{int}					
\DeclareMathOperator{\diam}{diam}					
			
\DeclareMathOperator{\N}{\mathbb N}			
\DeclareMathOperator{\R}{\mathbb R}			
\DeclareMathOperator{\Z}{\mathbb Z}			

\DeclareMathOperator{\id}{id} 

\newcommand{\acts}{\ \rotatebox[origin=c]{-180}{$\circlearrowright$}\ }

\DeclareMathOperator{\GL}{GL} 
\newcommand{\G}{\mathcal G}

\newcommand{\difrac}[2]{\frac{\diff#1}{\diff#2}} 

\newenvironment{quot}
{
	\vspace{-0.2cm}
	\vspace{0.2cm}
}

\theoremstyle{definition}
\newtheorem{d1}{Definition}[section]

\newenvironment{defin}
{
	\begin{quot}
		\begin{d1}
		}
		{\end{d1}
	\end{quot}

}

\theoremstyle{definition}
\newtheorem{r1}[d1]{Remark}

\newenvironment{rem}
{
	\begin{quot}
		\begin{r1}
		}
		{\end{r1}
	\end{quot}
}

\theoremstyle{definition}
\newtheorem{e1}[d1]{Exercise}

\theoremstyle{definition}
\newtheorem{ese1}[d1]{Example}

\newenvironment{ese}
{
	\begin{quot}
		\begin{ese1}
	}
	{	
		\end{ese1}
	\end{quot}
}

\theoremstyle{definition}
\newtheorem{f1}[d1]{Formulation}

\theoremstyle{definition}
\newtheorem{f2}[d1]{Fact}

\theoremstyle{definition}

\theoremstyle{definition}

\theoremstyle{definition}
\newtheorem{t1}[d1]{Theorem}

\newenvironment{thm}
{
	\begin{quot}
		\begin{t1}}
		{\end{t1}
	\end{quot}
}

\theoremstyle{definition}
\newtheorem*{T1*}{Theorem}

\newenvironment{teor*}
{
	\begin{quot}
		\begin{T1*}}
		{\end{T1*}
	\end{quot}
}

\newenvironment{dimo}
{\begin{proof}[Proof]
	}
	{\end{proof}}

	\theoremstyle{definition}
	\newtheorem{l1}[d1]{Lemma}
	
	\newenvironment{lem}
	{
		\begin{quot}
			\begin{l1}}
			{\end{l1}
		\end{quot}
	}
	\theoremstyle{definition}
	\newtheorem{p1}[d1]{Proposition}
	
	\newenvironment{prop}
	{
		\begin{quot}
			\begin{p1}}
			{\end{p1}
		\end{quot}
	}
	
	\theoremstyle{definition}
	\newtheorem{c1}[d1]{Corollary}
	
	\newenvironment{cor}
	{
		\begin{quot}
			\begin{c1}}
			{\end{c1}
		\end{quot}
	}

		\renewenvironment{abstract}
	{\list{}{\rightmargin\leftmargin}%
		\item[\textbf{Abstract:}]\relax}
	{\endlist}
	
	\newenvironment{dimoclaim}{\emph{Proof of Claim:}\;}{\hfill$\square$}

\newtheorem{innercustomthm}{Theorem}
\newenvironment{customthm}[1]
  {\renewcommand\theinnercustomthm{#1}\innercustomthm}
  {\endinnercustomthm}

 \newtheorem*{Theorem*}{Theorem}
 \newtheorem*{Proposition*}{Proposition}
 \newtheorem*{Lemma*}{Lemma}

\graphicspath{{mpdraws/}{draws/}}

\setenumerate[1]{label=(\arabic*)}

\numberwithin{equation}{section}
\numberwithin{figure}{section}

\newtheorem*{ex*}{Exercise}

\theoremstyle{definition}

\makeatletter
\newcommand\RedeclareMathOperator{%
  \@ifstar{\def\rmo@s{m}\rmo@redeclare}{\def\rmo@s{o}\rmo@redeclare}%
}
\newcommand\rmo@redeclare[2]{%
  \begingroup \escapechar\m@ne\xdef\@gtempa{{\string#1}}\endgroup
  \expandafter\@ifundefined\@gtempa
     {\@latex@error{\noexpand#1undefined}\@ehc}%
     \relax
  \expandafter\rmo@declmathop\rmo@s{#1}{#2}}
\newcommand\rmo@declmathop[3]{%
  \DeclareRobustCommand{#2}{\qopname\newmcodes@#1{#3}}%
}
\@onlypreamble\RedeclareMathOperator
\makeatother

\renewcommand\N{\ensuremath{\mathbb{N}}}
\renewcommand\Z{\ensuremath{\mathbb{Z}}}

\renewcommand\R{\ensuremath{\mathbb{R}}}

\renewcommand\G{\ensuremath{\mathcal{G}}}

\RedeclareMathOperator{\tr}{tr}
\RedeclareMathOperator{\diam}{diam}

\DeclareMathOperator{\Vol}{Vol}
\DeclareMathOperator{\sys}{sys}
\DeclareMathOperator{\Area}{Area}
\RedeclareMathOperator{\interior}{int}
\RedeclareMathOperator{\id}{id}
\newcommand{\del}{\partial}

\renewcommand{\epsilon}{\varepsilon}
\renewcommand{\bar}{\overline}

\RedeclareMathOperator{\GL}{\mathit{GL}}

\DeclareMathOperator{\PSL}{\mathit{PSL}}

\newcommand{\GC}{\mathcal{GC}}

\renewcommand{\P}{\mathcal P}

\usetikzlibrary{arrows}
\tikzset{
    labl/.style={anchor=south, rotate=90, inner sep=.5mm}
}
\tikzstyle{every picture}=[> = to]
\tikzset{cdlabel/.style={execute at begin node=$\scriptstyle,execute at end node=$}}
\tikzset{implication/.style={double equal sign distance, -implies}}
\tikzset{biimplication/.style={double equal sign distance, implies-implies}}

\renewcommand{\hat}{\widehat}
\newcommand{\RR}{\mathbb{R}}

\newcommand{\ZZ}{\mathbb{Z}}
\newcommand{\NN}{\mathbb{N}}
\newcommand{\HH}{\mathbb{H}}

\newcommand{\II}{\mathbf{1}}

\newcommand{\TT}{\mathbb{T}}
\newcommand{\mc}[1]{\mathcal{#1}}

\newcommand{\mf}[1]{\mathfrak{#1}}
\newcommand{\closure}[1]{\overline{#1}}

\newcommand{\Sob}{\operatorname{\mathcal{S}}}
\newcommand{\kl}[1]{\left ( #1 \right)}
\newcommand{\UT}{\operatorname{UT}}
\newcommand{\PT}{\operatorname{PT}}
\newcommand{\supp}{\operatorname{supp}}
\renewcommand{\tilde}[1]{\widetilde{#1}}
\newcommand{\Vect}{\mathfrak{X}}
\newcommand{\ba}{\begin{align*}}
\newcommand{\ea}{\end{align*}}
\newcommand{\Exp}{\operatorname{Exp}}
\newcommand{\vol}{\operatorname{Vol}}
\newcommand{\lquot}{\backslash}
\newcommand{\Ksob}{K_{\text{Sob}}}
\renewcommand{\sl}{\mf{sl}}
\newcommand{\Ad}{\operatorname{Ad}}
\newcommand{\inj}{\operatorname{inj}}

\begin{document}

\title{Volume bound for the canonical lift complement of a random geodesic}

\author[T. Cremaschi]{Tommaso Cremaschi}
\address{
Department of Mathematics\\
University of Southern California\\
Los Angeles, CA 90089\\
USA}
\email{cremasch@usc.edu}

    \author[Y. Krifka]{Yannick Krifka}
\address{
Max-Planck-Institut für Mathematik\\
53111 Bonn\\
Deutschland}
\email{krifka@mpim-bonn.mpg.de}

\author[D. Martínez-Granado]{Dídac Martínez-Granado}
\address{Department of Mathematics\\
         University of California, Davis,
         Davis CA 95616\\
         USA}
\email{dmartinezgranado@ucdavis.edu}

    \author[F. Vargas Pallete]{Franco Vargas Pallete}
\address{
Department of Mathematics\\
         Yale University\\
New Haven, CT 08540\\
USA}
\email{franco.vargaspallete@yale.edu}

\maketitle

\begin{abstract}
Given a filling primitive geodesic curve in a closed hyperbolic surface one obtains a hyperbolic three-manifold as the complement of the curve’s canonical lift to the projective tangent bundle. In this paper we give the first known lower bound for the volume of these manifolds in terms of the length of generic curves. We show that estimating the volume from below can be reduced to a counting problem in the unit tangent bundle and solve it by applying an exponential multiple mixing result for the geodesic flow.
\end{abstract}

\section{Introduction}
\label{sec:intro}

\subsection{Volumes of lift complements of curves}
\label{subsec:volumes}

In the following we consider a hyperbolic surface $X \cong S_{g,n}$ of genus $g$ with $n$ punctures.

Associated to $X$ is the $3$-manifold $\PT(X)$, the projectivised tangent bundle. For any finite collection of  smooth essential closed curves $\G$ on $X$, there is a \emph{canonical} lift $\hat \G$ in $\PT(X)$ realized by the set of tangent lines to $\G$.

Drilling $\hat \G$ from $\PT(X)$ produces a $3$-manifold $M_{\hat\G} = \PT(X) \setminus \hat \G$ and when $M_{\hat\G}$ is hyperbolic, by Mostow Rigidity \cite{BP1992}, any invariant of $M_{\hat\G}$ naturally becomes a mapping class group invariant of $\G$. When $\G$ is filling, its components are primitive, and $\G$ is in minimal position, Foulon and Hasselblatt \cite{FH13:HyperbolicLift} first observed that $M_{\hat\G}$ admits a complete hyperbolic metric of finite-volume, in particular $\text{Vol}(M_{\hat\G})$ is such an invariant. One should think of $\hat\G$ as a weak version of a link diagram where ``over'' and ``under'' crossings are encoded by the tangent directions to $\G$. The most general such hyperbolicity result appears in \cite{CRM19:VolumeSelfIntersection}. The authors show that if one takes a primitive filling system $\G$ in minimal position over a surface $X$ and then drills a transverse lift $\bar\G$ in a Seifert-fibered manifold $M$ then, the resulting manifold is hyperbolic. Transverse lifts of such systems will be called \emph{topological lifts}. Canonical lifts are examples of topological lifts.

In the rest of the paper we will use $\hat\G$ to denote canonical lifts of $\G\subset S$ in $\PT(S)$ and $\bar\G$ to denote topological lifts of $\G\subset X$ in a Seifert-fibered manifold $M$.

Several upper and lower bounds for $\text{Vol}(M_{\bar\G})$ in terms of invariants of $\G$ have been studied in recent literature, see \cite{BPS17:FillingGeodesics,BPS19:Numerical, RM20:LowerBound, RM21:Periods,CRM19:VolumeSelfIntersection,CRMY20:Multicurves}. 

Going back to the special sub-class that arises by considering $M$ to be $\PT(X)$, or the unit tangent bundle, and using the bundle projection map $\pi\colon\PT(X)\rar X$ to lift a filling geodesic $\gamma$ to its canonical lift $\hat \gamma$. In the case that the surface $X$ is the modular surface $\PT(X)$ can be identified with the Trefoil complement in $\mathbb S^3$  and Ghys \cite{Gh2007} showed that all Lorenz knots and links arise as canonical lifts of geodesics on the modular surface.
Moreover, the setup of canonical lifts has been extensively studied in \cite{RM21:Periods,BPS17:FillingGeodesics,BPS19:Numerical,RM20:LowerBound} and others.

In \cite{CRM19:VolumeSelfIntersection} the authors gave an upper bound which is linear in terms of the self-intersection number of $\G$, a fact reminiscent of classical results in knot theory. In \cite{BPS17:UpperBound}, it is shown that for every hyperbolic structure $X$ on $S$, there is a constant $C_X$ such that $\text{Vol}(M_{\hat\G}) \leq C_X \ell_X(\G)$, where $\ell_X(\G)$ denotes the length of the geodesic representative. Observe that $\vol(M_{\hat\G})$ is independent of the choice of $X$ and, in fact, it is mapping class group invariant (see \cite[2.1]{RM20:LowerBound}). Using this one can easily construct a sequence of filling closed curves whose volume $\vol(M_{\hat\G})$ stays bounded while its hyperbolic length goes to infinity. See \cite{RM20:LowerBound,RM21:Periods} for many other interesting examples.

There are key differences between the volumes corresponding to canonical lifts or to topological lifts and also between simple filling systems and closed filling geodesics. In \cite[Corollary~1.6]{CRM19:VolumeSelfIntersection} the authors construct examples in which the volume of a topological lift, which is not canonical, is asymptotic to the self-intersection number $\iota(\gamma_n,\gamma_n)$ of the filling curves. Fixing a hyperbolic structure $X$ on $S$, the self-intersection number is bounded above by $\ell_X(\gamma_n)^2$ by a result of Basmajian~\cite{Basmajian2013}, which is in contrast with the general length upper bound for volumes of canonical lifts of \cite{BPS17:UpperBound}.

The best known lower bound appears in \cite{RM20:LowerBound,CRM19:VolumeSelfIntersection}, where the bound is given in terms of the number of essential homotopy classes of arcs of $\G$ after cutting $X$ open along any multi-curve $\mathfrak m$\, and taking the maximum over such $\mathfrak m$. While this lower bound is shown to be sharp for some families of {\it non-simple} closed curves on the modular surface \cite{RM20:LowerBound}, it is always at most $6(3g + n)(3g -3 + n)$ whenever $\G$ is composed entirely of {\it simple} closed curves. This is addressed in \cite{CRMY20:Multicurves}.

In \cite{CRMY20:Multicurves} the authors study the setting in which one considers a filling collection of simple closed curves in minimal position instead of a primitive filling curve in minimal position. This is interesting because the only known lower bound is completely ineffective in this case (see Remark~\ref{rmk:ineffective}). In \cite[Theorem~A]{CRMY20:Multicurves} the authors relate the volume of complement of the canonical lift of a pair of filling geodesics to pants distance in the pants graph of the surface.

\subsection{Upper bounds in terms of length}
\label{subsec:upperbound}

We now describe in more details the length upper bound and refer to some numerical evidence. The length upper bound of \cite[Theorem~1.1]{BPS17:UpperBound} is:

\begin{Theorem*}  Let $X$ be a hyperbolic surface. Then, there exists $C_X$ such that for any filling primitive geodesic $\gamma\subset X$:
\be \text{Vol}(M_{\hat\gamma})\leq C_X\ell_X(\gamma)\ee
\end{Theorem*} 
and the result also works for multi-curves. The proof goes by showing the equivalent result in the case of the modular surface $Y=\mathbb H^2/SL_2(\Z)$. Then, one shows that by taking branched coverings $Z$ of $Y$ and curves $\gamma\subset Y$ one can obtain all filling primitive systems on any $S_{g,n}$. Then, the constant $C_X$ comes from considering the optimal quasi-conformal map from $X$ to $Z$, which is in general a non-trivial problem.

\brem \label{rmk:twist} There are large families of filling primitive geodesics for which the volume of the complements are uniformly bounded but whose length go to infinity. The easiest such example can be obtained by taking the mapping class group orbit of a fixed curve $\gamma$. However, there are also more interesting examples in which the curves $\gamma_n$ are not in the same mapping class group orbit, see \cite{RM20:LowerBound,RM21:Periods}. Some of these examples can be thought of as taking a filling curve $\gamma$ intersecting another curve $\alpha$ once and concatenating $\gamma$ with powers of $\alpha$. These are called \emph{twist families} and will always give rise to bounded volumes families.
\erem 

A sequence of random geodesics is, informally, a sequence of geodesics that gets more and more equidistributed with respect to the volume of $\UT(X)$ and converges, up to scaling, to the Liouville measure of $\UT(X)$. See Subsection~\ref{subsec:geometry} for precise definitions. 

\begin{figure}[h!]
    \centering
    \includegraphics[width=\textwidth]{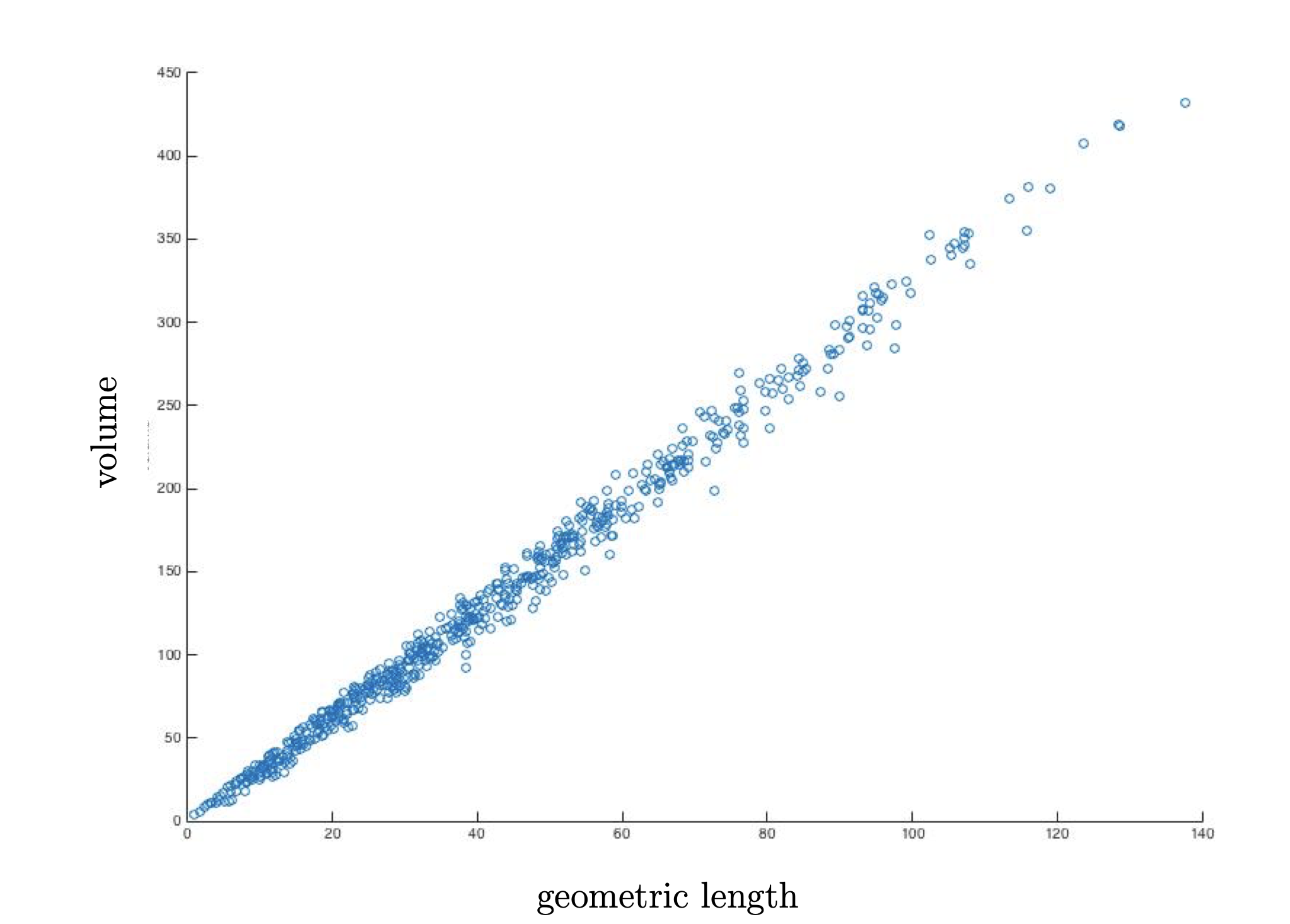}
    \caption{This figure shows the volume of the canonical lift complement as a function of hyperbolic length for a sequence of random geodesics constructed via number theoretic techniques by Duke~\cite{Du88:UniformDistrModular}. This corresponds to Figure 2 of \cite{BPS19:Numerical} }\label{lingrowth}
\end{figure}

A sequence of random geodesics in the modular surface has been considered by Duke in \cite{Du88:UniformDistrModular}. This model is constructed via number theoretic techniques. In this paper we will construct another family of random geodesics using geometry and dynamics.
It would be interesting to obtain a lower bound for the volume of Duke's random sequence.
In \cite{BPS19:Numerical} the authors compute the volumes for finitely many terms of Duke's random sequence and then give numerical evidence of the linear volume growth as a function of geodesic length (see Figure~\ref{lingrowth}).

However, by \cite[Corollary~1.2]{RM21:Periods} there exist multi-curves whose volumes are asymptotic to $\f L{W(L)}$ for $L$ the length and $W(x)$ the Lambert function. The \emph{Lambert function} is the principal branch of the inverse of $f(w)=w\log w$ which is asymptotic to $\log(x)-\log\log(x)+o(1)$.

\begin{figure}[h!]
\centering{
\resizebox{140mm}{!}{\fontsize{12pt}{12pt}\selectfont
\begingroup%
  \makeatletter%
  \providecommand\color[2][]{%
    \errmessage{(Inkscape) Color is used for the text in Inkscape, but the package 'color.sty' is not loaded}%
    \renewcommand\color[2][]{}%
  }%
  \providecommand\transparent[1]{%
    \errmessage{(Inkscape) Transparency is used (non-zero) for the text in Inkscape, but the package 'transparent.sty' is not loaded}%
    \renewcommand\transparent[1]{}%
  }%
  \providecommand\rotatebox[2]{#2}%
  \newcommand*\fsize{\dimexpr\f@size pt\relax}%
  \newcommand*\lineheight[1]{\fontsize{\fsize}{#1\fsize}\selectfont}%
  \ifx\svgwidth\undefined%
    \setlength{\unitlength}{859.67999268bp}%
    \ifx\svgscale\undefined%
      \relax%
    \else%
      \setlength{\unitlength}{\unitlength * \real{\svgscale}}%
    \fi%
  \else%
    \setlength{\unitlength}{\svgwidth}%
  \fi%
  \global\let\svgwidth\undefined%
  \global\let\svgscale\undefined%
  \makeatother%
  \begin{picture}(1,0.59854832)%
    \lineheight{1}%
    \setlength\tabcolsep{0pt}%
    \put(0,0){\includegraphics[width=\unitlength,page=1]{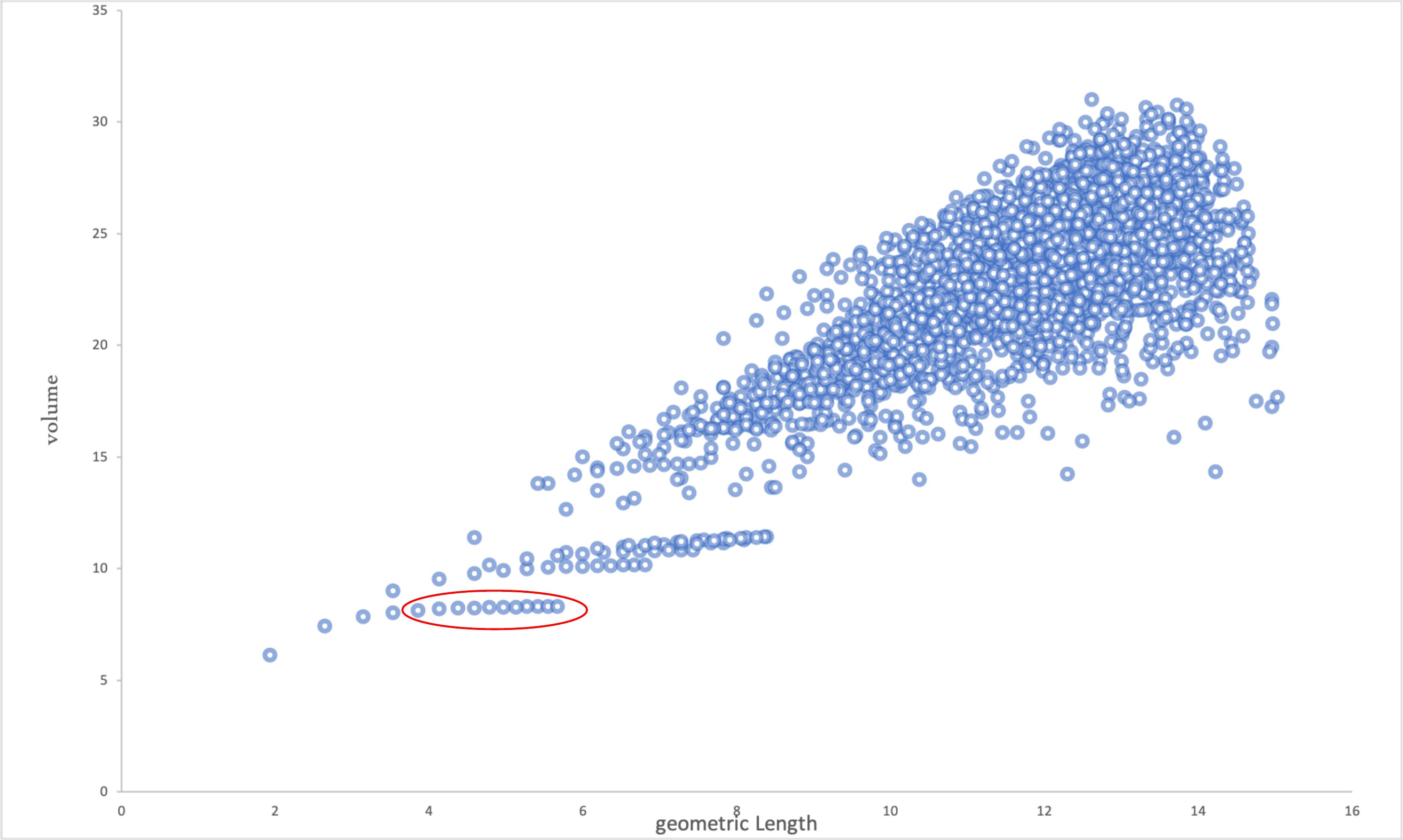}}%
    \put(0.48896736,0.11134619){\color[rgb]{1,0,0}\makebox(0,0)[lt]{\lineheight{1.25}\smash{\begin{tabular}[t]{l}$\mbox{twist families}$\end{tabular}}}}%
    \put(0,0){\includegraphics[width=\unitlength,page=2]{chart.pdf}}%
  \end{picture}%
\endgroup%
}
    \caption{This figure shows the volume of canonical lift complement as a function of hyperbolic length for geodesics of word length at most 16. }\label{volpic}
}
\end{figure}

In forthcoming work, Yarmola and Intrater compute the volumes for all geodesics of length at most $16$ in the modular surface (see Figure~\ref{volpic}). The graph  clearly still shows a linear upper bound but the situation for the lower bound is more complicated due to the presence of twist families. Such a family can be seen in the lower left corner of Figure~\ref{volpic}.

\subsection{Lower bound} 
\label{subsec:lowerbound}

\begin{figure}[h!]
    \centering{
    \resizebox{140mm}{!}{\fontsize{12pt}{12pt}\selectfont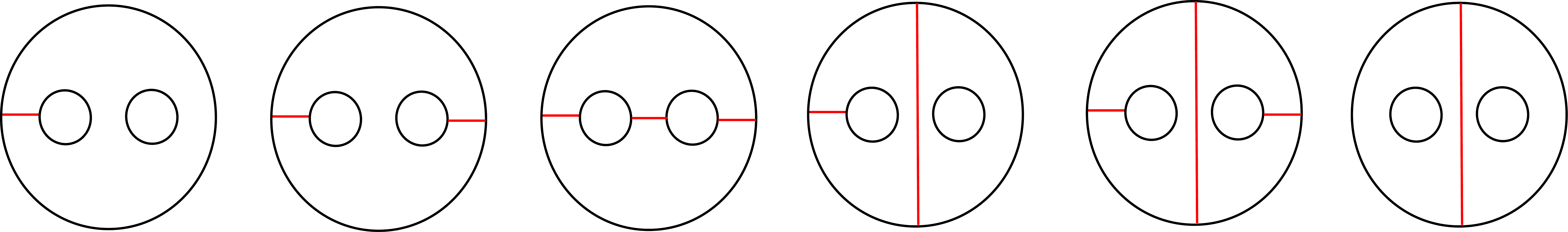}
    \caption{This figure shows the 6 $\gamma$-arcs configurations (up to homotopy) consisting of simple (undirected) $\gamma$-arcs}
    \label{fig:simple}
    }
\end{figure}

 No geometric lower bounds are currently known and we only have combinatorial ones. Using work of Agol, Storm and Thurston \cite{AgolStormThurston5}, Rodriguez-Migueles showed in \cite{RM20:LowerBound} that:

\begin{Theorem*}[Combinatorial Lower Bound] 
Let $\P$ be an essential surface decomposition of $S$ and let $\gamma$ be a filling primitive curve in minimal position with respect to $\partial \P$ and itself. Then:
\be\f{v_3}2\sum_{Q} \kl{\#{\set{\gamma\text{-arcs in }Q}}-3}\leq \text{Vol}(M_{\hat\gamma}),\ee

where $v_3$ is the volume of a regular ideal tetrahedra and we sum over all components of the pants decomposition $\P$. For a pair of pants $Q$, the \emph{$\gamma$-arcs in $Q$} are the connected components of $\gamma\cap Q$ (see Figure~\ref{fig:simple}). However, when writing ``$\#{\set{\gamma\text{-arcs in }Q}}$'' we mean the number of \emph{homotopy classes} of $\gamma$-arcs in $Q$ with endpoints gliding on the boundary.
\label{thm:lowerbound}
\end{Theorem*} 

For example if $Q$ is a pair of pants and $\gamma$ is a filling geodesic there are at most $6$ simple $\gamma$-arcs in $\gamma\cap Q$ up to homotopy.

\brem \label{rmk:ineffective}
For a simple multi-curve $\gamma$, the Combinatorial Lower Bound Theorem is ineffective since the number of $\gamma$-arcs is upper bounded by a constant independent of $\gamma$ (see Figure~\ref{fig:simple}).
\erem

In \cite{CRM19:VolumeSelfIntersection} the authors show that the above lower bound also works in the setting of Seifert-fibered spaces.

\subsection{Statement of the theorem}
\label{subsec:theorem}
From now on, we will assume we are working with a surface without punctures, $S=S_{g,0}$, unless otherwise stated.
Because of the length upper bound it would be nice to have a lower bound also depending on length, possibly in a linear fashion. The possible linearity of a lower-bound is hinted by the experimental data of Figure \ref{lingrowth}, consisting of a finite number of terms of a sequence of random geodesics that converge to the Liouville measure (see discussion on Duke's model in page 3).

However, as mentioned in Remark~\ref{rmk:twist}, one can construct sequences of geodesics $\gamma_n$ such that $\ell_X(\gamma_n)\rar\infty$ but the volume stays bounded. Similarly, there are examples of  $\gamma_n's$ where the growth is bounded above by a sub-linear function (see \cite{RM21:Periods}). This motivates to ask what the generic behavior of the volume is for long geodesics.

To this end let us recall the notion of asymptotic density for collections of geodesics. Given a hyperbolic surface $X$ we denote its set of closed geodesics by $\mc{G}$. For a collection of closed geodesics $\mc{A} \subseteq \mc{G}$ we denote by
$\mc{A}(R)$ all geodesics in $\mc{A}$ of length less than $R >0$. We say that $\mc{A}$ has \emph{asymptotic density $\rho \in [0,1]$} if
\[ \lim_{R \to \infty} \frac{ \# \mc{A}(R) }{ \# \mc{G}(R) } = \rho.\]
Intuitively, the density $\rho$ describes approximately the likelihood of drawing a geodesic in the subset $\mc{A}$ by sampling uniformly at random every closed geodesic of length at most $R$ for very large $R>0$. This definition has its origins in number theory where one is interested in the asymptotic density of infinite subsets of the natural numbers.

Our main result is the following.

\begin{customthm}{\ref{maintheorem}}
Let $X$ be a closed hyperbolic surface glued without twists from isometric pairs of pants. Further, denote by $0<\delta<1$ the Hausdorff dimension of the limit set of any one of the pairs of pants.
 
Then for every $\eta > 1 (> \delta > 0)$ and every $F(\ell) \in O(\ell^{\delta/2\eta})$ the collection of primitive, filling geodesics $\gamma$ satisfying the volume estimate
\[  F(\ell_X(\gamma)) \leq \vol(M_\gamma) \]
has asymptotic density $1$.
\end{customthm}

\brem 
Since we will only deal with canonical lifts, we simplify our notation and denote the manifold obtained by drilling the canonical lift of a closed geodesic $\gamma$ from $\PT(X)$ by $M_{\gamma}$.
\erem

\subsection{Short outline }

After reviewing some background material in Section \ref{sec:background} and giving a more detailed sketch of proof for Theorem \ref{maintheorem} in Section \ref{sec:maintheorem} we start in Section \ref{sec:random} by constructing random sequences of closed geodesics as follows.

Let us fix $\epsilon>0$ small enough and denote for every $v \in \UT(X)$ its set of $\epsilon$-return times by $R_\epsilon(v)$, i.e.\ those $t > 0$ such that $g_t(v)$ is $\epsilon$-close to $v$.
For every $t \in R_\epsilon(v)$ we can close-up the geodesic segment $\gamma_v[0,t]$ of length $t$ starting at $v$ by an $\epsilon$-short arc to obtain a broken-geodesic curve $\gamma_v(t)$. We will be interested in its geodesic representative which we shall denote by $\hat{\gamma}_v(t)$. Note that a similar model has been considered by Bonahon \cite[Page~151]{Bonahon88:GeodesicCurrent} (see also \cite[Theorem~3]{Bri04:Random}, \cite[Claim~2.3]{Sap17:Birman}).

The relevance of this construction is that by Lemma \ref{modelrel} a collection of closed geodesics $\mc{A} \subseteq \mc{G}$ has asymptotic density one if for $\mu$-almost every $v \in \UT(X)$ and all sufficiently long $\epsilon$-return times $t \geq t(v)$ the constructed closed geodesic $\hat{\gamma}_v(t)$ is in $\mc{A}$.
This observation will allow us to use dynamical properties of the geodesic flow more directly.

First, we show in Lemma \ref{lem:conv_to_Liouville} that as $t \in R_\epsilon(v)$ tends to infinity, $\hat{\gamma}_v(t)$ converges to the Liouville current, so that $\hat{\gamma}_v(t)$ will be filling for all sufficiently large $t$ (Lemma \ref{curveisfilling}). 

Next, we want to use the Combinatorial Lower Bound Theorem to find a lower bound for $\Vol(M_{\hat{\gamma}_v(t)})$. To this end we show in Proposition \ref{prop:arcs_pulled_tight} that the number (of homotopy classes) of $\gamma_v[0,t]$-arcs is roughly the number (of homotopy classes) of $\hat{\gamma}_v(t)$-arcs up to a uniform error.

In Section \ref{sec:lower_bound} we use dynamical techniques involving higher order exponential mixing of the geodesic flow to estimate the number of $\gamma_v[0,t]$-arcs (see Theorem \ref{mainestimate}).

Finally, we prove the Main Theorem \ref{maintheorem} in Section \ref{sec:proof}. Using Theorem \ref{mainestimate} we first prove Theorem \ref{auxmainthm} in Section \ref{sec:proof}, which shows that for $\mu$-almost every $v \in \UT(X)$ the constructed geodesic curves $\hat{\gamma}_v(t)$ satisfy the asserted volume estimate for sufficiently large times $t$. Theorem \ref{maintheorem} then follows from Lemma \ref{modelrel}.

\textbf{Acknowledgments:} This material is based upon work supported by the National Science Foundation under Grant No.\ DMS-1928930 while the authors participated in a program hosted by the Mathematical Sciences Research Institute in Berkeley, California, during the Fall 2020 semester. The authors would also like to thank Dunfield, Einsiedler, Prohaska, Hamenst\"adt and Blayac for helpful discussions. F. Vargas Pallete's research was supported by NSF grant DMS-2001997.
\section{Background}
\label{sec:background}

\subsection{Notation}

We will use the following notational conventions:
\begin{itemize}
    \item $A\cong B$ denotes that the two topological spaces $A,B$ are homeomorphic;
    \item $\alpha\simeq \beta$ denotes that $\alpha, \beta$ are homotopic maps or spaces, generally curves or arcs.
\end{itemize}
For two functions $f, g \colon \TT \to \RR_+$ we will use the following notations, where $\TT$ stands for $\NN$, $\ZZ$ or $\RR$:
\begin{itemize}
    \item $f(t)\sim g(t)$ means that the two functions $f, g$ are asymptotic, i.e.\ 
    \[\lim_{t \to \infty} \frac{f(t)}{g(t)} = 1;\]
    \item $f(t) \ll g(t)$ means that there is a constant $C>0$ and $T_0 \in \TT$ such that $f(t) \leq C g(t)$ for all $t \geq T_0$;
    \item $f(t) \asymp g(t)$ means that there are constants $C_1, C_2 > 0$ and $T_0 \in \TT$ such that 
    \[C_1 f(t) \leq g(t) \leq C_2 f(t)\] 
    for all $t \geq T_0$.
\end{itemize}

\subsection{Dynamics}
\label{subsec:dynamics}

We recall some concepts of dynamics that will play a key role in this paper.
Standard references are~\cite[Page~151]{HK95:Dynamical}, \cite{Wal82:ergodic}.

Let $X$ be a topological space equipped with a (Borel) probability measure $\mu$.
Suppose that $g_t \colon X \to X$ is a continuous \emph{flow}, for $t \in \mathbb{R}$, i.e.\ a family of continuous maps so that $g_0$ is the identity and
$g_{s+t}=g_s \circ g_t$ for all $s,t \in \mathbb{R}$.
We say $\mu$ is $g$-invariant if for any Borel set $A \subset X$, $\mu(A)=\mu(g_t (A))$ for all $t \in \mathbb{R}$.
There is an analogous picture for discrete time dynamical systems, where instead of a flow, we have a continuous map $T \colon X \to X$ (which need not be invertible) and $\mu$ is said to be $T$-invariant if $\mu(B)=\mu(T^{-1}(B))$ for every Borel set $B \subset X$.
One can go from a continuous to a discrete setting by taking $T=g_1$, the \emph{time-one flow map}, and we will be doing so throughout the paper.
We will focus on two properties that a continuous or discrete dynamical system can have: ergodicity and mixing.

\subsubsection{Ergodicity}

The property of ergodicity will be used to construct our geometric random model.
A flow $g_t \colon X \to X$ (resp. transformation $T \colon X \to X$) is \emph{ergodic} with respect to $\mu$ if for every subset $A \subset X$ which is $g$-invariant $g_t(A)=A$ (resp. $T$-invariant $T^{-1}(A)=A$),  then either $\mu(A)=0$ or $\mu(X \backslash A)=0$.

If $P_x$ is a property depending on a point $x \in X$, we say that $P$ holds for  \emph{$\mu$-almost every $x\in X$} if it holds for a subset $B \subset X$ so that $\mu(B)=1$.

Ergodic dynamical systems satisfy the following two properties that we will use.

\bthm[{\cite[Theorem~1.7]{Wal82:ergodic}}]
Let $X$ be a compact metric space, $T \colon X \to X$ a continuous transformation and $\mu$ a probability measure so that $\mu(U) > 0$ for all open set $U \subset X$.
      Suppose that $\mu$ is $T$-invariant and $T$ is ergodic with respect to $\mu$.
      Then for $\mu$-almost every $x \in X$ the orbit of $x$, i.e. $\{ T^n(x) : x \in X \}$, is dense.
      \label{thm:fulldense}
\ethm

\bthm[Birkhoff's Ergodic Theorem; {\cite[Theorem~1.14]{Wal82:ergodic}}] \label{birkhoff} Let $X$ be a topological space equipped with a probability measure $\mu$. For $g_t\colon X\rar X$ an ergodic flow with respect to $\mu$, and $f$ any measurable function, we have that for $\mu$-almost every $x\in X$:
\[\lim_{T\rar\infty}\frac 1 T \int_0^T f(g_t(x))dt= \int_Xf d\mu \]
\label{thm:birkhoff}
\ethm
For a discrete system, the flow is replaced by a transformation $T$, the integral by a sum, and the limit in the previous theorem takes the following form

\[
\lim_{N\to \infty} \frac{1}{N}\sum_{n=0}^{N-1} f(T^n (x)) = \int_X f d\mu
\]
for $\mu$-almost every $x \in X$.

One also says that the orbits $\{g_t(x) \, | \, t \in [0,T]\}$ resp.\ $\{x, T(x), \ldots, T^N(x) \}$ equidistribute as $T \to \infty$ resp.\ $N \to \infty$.

\subsubsection{Mixing}

We will need a stronger property than ergodicity to have enough control on our model.
We say $\mu$ is \emph{mixing} for the flow $g_t \colon X \to X$ if for all $f,g \in L^2(X,\mu)$, the \emph{correlation function}
\[
\rho(t) \coloneqq \int_X (f \circ g_t)g d\mu - \int_X f d\mu \int_X g d\mu,
\]
satisfies $\rho(t) \to 0$ as $t \to +\infty$.
For discrete systems, the correlation is written as
\[
\rho(n) \coloneqq \int_X (f \circ T^n) g d\mu - \int_X f d\mu \int_X g d\mu,
\]
and mixing means $\rho(n) \to 0$ as $n \to +\infty$.

We will be interested in extending this correlation function to more than two functions, as well as in quantifying its decay. This will lead to the notion of exponential $k$-mixing, see Appendix~\ref{app:mixing}.

Another equivalent way to phrase mixing is that, for all Borel sets $A,B \subset X$,
\[
    \lim_{n \to \infty} \mu(T^{-n}(A) \cap B) =\mu(A)\mu(B)
\]
It follows that mixing implies ergodicity.

Intuitively, mixing means that for any two Borel sets $A,B \subseteq X$ the events $x \in B$ and $T^n(x) \in A$ become more and more independent as $n \to \infty$.

In this paper, we will let $g_t$ be the geodesic flow on the unit tangent bundle $\UT(X)$, and $\mu$
will be the normalized Liouville probability measure. $T$ will be the time-one map of the geodesic flow.
The geodesic flow is mixing with respect to $\mu$ (see \cite[Theorem~1]{Bab02:Mixing}), and is by the above discussion, also ergodic.

\subsubsection{Anosov closing}

It is known that the geodesic flow on the unit tangent bundle $Y = \UT(X)$ of a closed hyperbolic surface $X$ is \emph{hyperbolic/Anosov}, i.e.\ there is a smooth $g_t$ invariant splitting $TY = E^s \oplus E^c \oplus E^u$ and constants $C \geq 1$, $\alpha > 0$ such that
\begin{itemize}
    \item $E^c_x$ is spanned by the direction of the flow $\frac{d}{dt}\big|_{t=0} g_t(x)$,
    \item $\|g_t(v)\|_{g_t(x)} \leq C e^{-\alpha t} \|v\|_x$ for all $v \in E^s_x$, i.e.\ $g_t$ is exponentially contracting along $E^s$,
    \item $\|g_{t}(w)\|_{g_{t}(x)} \geq C^{-1} e^{\alpha t} \|w \|_x$ for all $w \in E^u_x$, i.e.\ $g_t$ is exponentially expanding along $E^u$,
\end{itemize}
for all $x \in Y$, $t \geq 0$. The vectors in $E^s$ (resp.\ $E^u$) are called the \emph{stable} (resp.\ \emph{unstable}) directions.

In our concrete situation we can describe this splitting in the following way. Let $\Gamma \leq G=\PSL(2,\RR)$ be a torsion-free cocompact lattice such that $X = \Gamma \lquot \HH^2$. The action of $G$ on $\UT(\HH^2)$ is free and transitive, and thus one may identify $\UT(\HH^2) \cong G$. This descends to an identification $\Gamma \lquot \UT(\HH^2) \cong \UT(X) \cong \Gamma \lquot G$. Using the right-action of $G$ on $\Gamma \lquot G$ we can identify the Lie algebra $\mf{g} = \sl_2(\RR)$ with the tangent space $T_x (\Gamma \lquot G)$ for every $x \in \Gamma \lquot G$:
\[ \mf{g} \coloneqq \mf{sl}_2(\RR) \to T_{x} (\Gamma \lquot G), \quad X \mapsto \frac{d}{dt}\bigg|_{t=0} x \cdot \exp(t \cdot X).\]

Recall that a basis of $\sl_2(\RR)$ is given by the matrices
\[ 
    H \coloneqq \begin{pmatrix} 
            \frac{1}{2} & 0 \\  
            0 & -\frac{1}{2} 
        \end{pmatrix}, \quad
    X_+ \coloneqq \begin{pmatrix} 
            0 & 1 \\  
            0 & 0
        \end{pmatrix}, \quad
    X_- \coloneqq \begin{pmatrix} 
            0 & 0 \\  
            1 & 0
        \end{pmatrix}.
\]
This basis amounts to an Anosov splitting of $T(\Gamma \lquot G)$ for the geodesic flow $g_t$ via the identification $T_x(\Gamma \lquot G) \cong \mf{g}$:
\[ E^c_x \cong \RR \cdot H, \quad E^s_x \cong \RR \cdot X_+, \quad E^u_x \cong \RR \cdot X_-.\]
Indeed, $g_t \colon \Gamma \lquot G \to \Gamma \lquot G$ is given by right-multiplication $g_t(x) = xa_t$ with the diagonal matrix
\[a_t\coloneqq \exp(t \cdot H) = \begin{pmatrix} e^{t/2} & 0\\ 0 & e^{-t/2} \end{pmatrix},\]
and
\[ \Ad_{a_t}(X_\pm) = a_t X_\pm a_{-t} =  e^{\pm t} X_\pm. \]

Anosov flows were subject to extensive research in the past which led to a detailed understanding of their dynamical properties; we refer to \cite{Hasselblatt} for a detailed exposition. One of these properties is the \emph{Anosov Closing Lemma} which intuitively says that near every ``almost-periodic'' orbit of an Anosov flow there is a truly periodic orbit. We will only need the following special version:

\blem[Anosov Closing Lemma] \label{lem:AnosovClosing}
    Let $X$ be a closed hyperbolic surface and let $\rho > 0$. There are positive numbers $0<\epsilon = \epsilon(\rho) \leq \rho$ and $0< T = T(\rho)$ such that the following holds:
    
    For every $t \geq T$, $v \in \UT(X)$ with $d(g_t(v),v) < \epsilon$ there exist $t' \in \RR$ and $v' \in \UT(X)$ such that
    \[ |t-t'| \leq \rho, \quad d(v,v') \leq \rho, \quad g_{t'}(v') = v'.\]
    
    Moreover, the geodesics $\gamma, \gamma'$ with respective starting vectors $\dot\gamma(0)=v, \dot\gamma'(0) = v'$ stay close to each other:
    \begin{align}
        d(\gamma(s \cdot t), \gamma'(s \cdot t')) \leq 2 \rho \label{est:stayclose}
    \end{align}
    for every $s \in [0,1]$. 
    
\elem

\begin{proof}
    This follows from \cite[4.5.15 Proposition]{Ebe96:GeometryNonpos} and the convexity of the distance function in $\HH^2$.
\end{proof}

\subsubsection{Asymptotic density of collections of closed geodesics}
\label{subsec:densityone}

Given a hyperbolic surface $X$ we will denote its set of closed geodesics by $\mc{G}$. Let $\mc{A} \subseteq \mc{G}$ be a collection of closed geodesics. For every $R>0$ we will denote by
\[ \mc{A}(R) \coloneqq \{ \gamma \in \mc{A} \, | \, \ell_X(\gamma) \leq R \} \]
all closed geodesics in $\mc{A}$ of length less than $R$. We say that $\mc{A} \subseteq \mc{G}$ has \emph{asymptotic density $ \rho \in [0,1] $} if
\[ \frac{ \# \mc{A}(R) }{ \# \mc{G} (R)} \to \rho \quad (R \to + \infty).\]

\brem \label{rmk:margulis_count}
    By work of Huber~\cite{Hub60:Growth} (and Delsarte and Selberg) it is known that 
    \[ \# \mc{G}(R) \sim \frac{e^R}{R} \]
    as $R \to +\infty$; (see also Margulis~\cite{margulis} for compact negatively curved manifolds).
\erem

\subsection{Topology}
\label{subsec:topology}
In the following sections we recall some facts and definitions about the topology of surfaces and 3-manifolds. For references, see \cite{He1976,Ha2007,Ja1980}.

\bdefi A \emph{knot} in $M$ will be any embedding of $\mathbb S^1$ into a 3-manifold $M$.
\edefi 
\bdefi Given a curve $\gamma\subset S$ and the projective tangent bundle $\PT(S)$ we define the \emph{canonical lift} to be the knot $\hat\gamma\subset 
\PT(S)$ obtained by taking the lift of $\gamma$ given by the tangential line field. 
\edefi 

We say that a curve is in \emph{minimal position} if the number of self-intersections is minimal in its homotopy class.

We recall the combinatorial lower-bound Theorem:
\bthm[{\cite[Theorem 1.5]{RM20:LowerBound}}] \label{thm:lowerbound2}
Let $\P$ be an essential surface decomposition of $S$ and let $\gamma$ be a filling primitive curve in minimal position with respect to $\partial \P$ and itself. Then:
\be\f{v_3}2\sum_{Q} \kl{\#{\set{\gamma\text{-arcs in }Q}}-3}\leq \text{Vol}(M_{\hat\gamma}),\ee

where $v_3$ is the volume of a regular ideal tetrahedra and we sum over all components of the pants decomposition $\P$.

\ethm

For a pair of pants $Q$, the \emph{$\gamma$-arcs in $Q$} are the connected components of $\gamma\cap Q$ (see Figure~\ref{fig:simple}). However, when writing ``$\#{\set{\gamma\text{-arcs in }Q}}$'' we mean the number of \emph{homotopy classes} of $\gamma$-arcs in $Q$ with endpoints gliding on the boundary.

\brem \label{rem:undirected}
    It is important to note that these $\gamma$-arcs are \emph{undirected}. Moreover, even though our main results involve curve complements in the projective tangent bundle, in order to avoid talking about orientations, our proofs will involve the unit tangent bundle. This is because we want to take advantage of the geodesic flow which is naturally defined in $\UT(X)$, a double cover of $\PT(X)$.
\erem

\subsection{Geometry}
\label{subsec:geometry}

A \emph{geodesic current} is a positive finite
  Radon measure $\mu$ on
$\UT(X)$ which is invariant under the geodesic flow, in the sense
that $(g_t)_*(\mu) = \mu$ for all $t\in\mathbb{R}$, where the subscript~$*$ denotes the push-forward of measures.

A closed geodesic $\gamma$ can be seen as a geodesic current. Consider the canonical
lift~$\hat{\gamma}$ of~$\gamma$ to $\UT(X)$; this is a periodic orbit
of~$g_t$. By abuse of notation, let $\gamma$ also denote the length-normalized $\delta$-function on
this orbit. That is, for an open set~$U$ we set $\gamma(U)$ to
be the total length of $\hat{\gamma} \cap U$ with respect to the
Riemannian metric $\UT(X)$.
The geometric intersection number between closed geodesics extends continuously to a bilinear form $i(\cdot, \cdot)$ on geodesic currents \cite[Proposition~4.5]{Bonahon86:EndsHyperbolicManifolds}.
The \emph{Liouville current} $\mathcal{L}_X$ associated to the hyperbolic metric $X$ is the Liouville volume of $\UT(X)$ normalized so that $i(\mathcal{L}_X,\gamma)=\ell_X(\gamma)$ for any closed geodesic $\gamma$.

A \emph{sequence of random geodesics} is a sequence~$(\gamma_n)_{n \in \mathbb{N}}$ of closed geodesics so that, after normalization, the geodesic currents corresponding to $\gamma_n$ converge in the weak$^*$-topology to the Liouville current, i.e.

\[
\lim_{n \to \infty} \frac{4\pi^2 |\chi(S)|}{\ell(\gamma_n)} \int_{\UT(X)} f d\gamma_n \to \int_{\UT(X)} f d\mathcal{L}_{X}
\]
for any continuous and compactly supported function $f \in C_c(\UT(X))$.
\bdefi[Unit tangent bundle decomposition]
\label{def:decomposition}
Let $\P$ be a geodesic pants decomposition of $X\cong S_{g,k}$. Pick a pant $P^j\in\P$ and let $\{o_i^j\}_{i\in\N}$ be the collection of orthogeodesics in $P^j$. For $\mu$-almost every vector $v \in \UT(P^j)$ there are minimal real numbers $a,b > 0$ such that the geodesic arc $\gamma \colon [-a,b] \to P^j$ with $\dot{\gamma}(0) = v$ intersects $\del P^j$ in $\gamma(-a)$ and $ \gamma(b)$. We define $U_i^j\subset \UT(P^j)$ to be the set of all directions $v$ such that the corresponding directed arc $\gamma \colon [-a,b] \to P^j$ is freely homotopic to $o_i^j$ in $P^j$ with gliding endpoints $\gamma(-a), \gamma(b) \in \del P^j$. This amounts to the following decomposition up to measure zero (see \cite[Section~7]{B11:Dilogarithm})
\[
\UT(X)=\bigcup_j \UT(P^j) = \bigcup_{i,j} U_i^j,\]
where the second equality is up to a $\mu$-measure zero set.
\edefi

\brem\label{rem:visiting_counting}
This definition is relevant in view of the combinatorial lower-bound from Theorem \ref{thm:lowerbound2}. Indeed, for a geodesic curve $\gamma$ the number of (homotopy classes of) $\gamma$-arcs is half the number of different sets $\{U_i^j\}_{i,j}$ that $\dot{\gamma}(t), t \in [0, \ell(\gamma)], $ visits. Here, the factor $\tfrac{1}{2}$ is due to the fact that $\gamma$-arcs are undirected; see Remark \ref{rem:undirected}. 
\erem

\subsection{Counting arcs} \label{subsec:counting_arcs}

Let $P$ be a hyperbolic pair of pants with totally geodesic boundary components. Given two (possibly the same) boundary components $C_-, C_+ \subseteq \partial P$ one may consider arcs starting on $C_-$ and ending on $C_+$. In each relative homotopy class of such an arc, where we allow each endpoint to glide on the respective boundary component, there is a unique geodesic arc. It meets $C_-$ and $C_+$ perpendicularly and is called an \emph{orthogeodesic arc}.

One may now ask how many orthogeodesics running from $C_-$ to $C_+$ of length $\leq \ell$ there are. Let us denote by $N_{C_-,C_+}(\ell)$ the number of such orthogeodesic arcs. 
In \cite{parkkonen-paulin-counting} Parkkonen and Paulin consider the same counting problem in the more general setting of pinched negatively curved manifolds and properly immersed closed locally convex subsets $C_-, C_+$. Applying their result \cite[Theorem 1]{parkkonen-paulin-counting} to our situation we obtain the following corollary:

\bcor \label{cor:counting_arcs}
    Let $0 < \delta < 1$ denote the Hausdorff dimension of the limit set of $P$. Then there is a constant $C_0 >0$ such that asymptotically
    \[ N_{C_-, C_+}(\ell) \sim C_0 \cdot e^{\delta \ell} \]
    as $\ell \to \infty$. 
\ecor

We want to point out that similar counting results in varying generality were obtained before by different authors; see \cite{parkkonen-paulin-survey} and the references therein.

\subsection{Sobolev Norms} \label{sec:sobolev_norms}

Denote $G \coloneqq \PSL(2,\RR)$ and let $\Gamma \leq G$ be a cocompact lattice. 

Let $d_G \colon G \times G \to \RR_{\geq 0}$ be a left-invariant metric on $G$. 
This metric descends to a metric on the quotient $d_{\Gamma \backslash G} \colon \Gamma \backslash G \times \Gamma \backslash G \to \RR_{\geq 0}$.
We may assume that $d_G$ is induced by a left-invariant Riemannian metric $\langle \cdot, \cdot \rangle$ on $G$. 
In particular, we can equip $G$ with such a metric via the identification $G \cong \UT(\HH^2)$.
We denote by $\nu$ the corresponding (bi-)invariant Haar measure on $G$. 

The left-action of $\Gamma$ on $G$ amounts to a quotient map $\pi \colon G \to \Gamma\backslash G, g \mapsto \Gamma g$. 
Because the action $\Gamma \acts G$ is isometric the metric $d_G$ descends to a metric $d_{\Gamma \backslash G} \colon \Gamma \backslash G \times \Gamma \backslash G \to \RR_{\geq 0}$, and we obtain a Riemannian metric on the quotient, which we shall denote by $\langle \cdot , \cdot \rangle$ as well, such that $\pi \colon G \to \Gamma \backslash G$ is a Riemannian covering map. 
The corresponding volume form amounts to a right-invariant quotient measure $\mu$ on $\Gamma \backslash G$. 
After possibly rescaling we may assume that $\mu$ is a probability measure. In this way $\mu$ coincides with the normalized Liouville measure on $\UT(\Gamma \lquot \HH^2)$ via the usual identification $\Gamma \lquot G \cong \UT(\Gamma \lquot \HH^2)$.

Note that the left-action 
\begin{align*}
    G \times \Gamma \backslash G &\to G, \\ 
    (g, \Gamma h) &\mapsto g \cdot \Gamma h \coloneqq \Gamma h g^{-1},
\end{align*}
is probability measure preserving by right-invariance of the quotient probability measure $\mu$.
Thus the regular representation $\lambda \colon G \to \mc{U}(L^2(\Gamma \backslash G))$ is unitary, where we denote
\[ (\lambda_g f)(\Gamma h) = f(g^{-1} \cdot \Gamma h) = f(\Gamma h g) \]
for every $g \in G$, $f \in L^2(\Gamma \backslash G, \mu)$.

More generally, whenever there is a smooth $G$-action $G \acts M$ on a smooth manifold $M$ (e.g.\ $M= G$ or $M = \Gamma \backslash G$), there is an induced action of the universal envelopping algebra $\mc{U}(\mf{g})$ on the space of smooth functions with compact support $C_c^\infty(M)$. This action is given via differentiation of the left regular representation $\lambda \colon G \to C_c^\infty(M)$ as follows
\[ (X \cdot \phi)(x) \coloneqq \frac{d}{dt}\bigg|_{t=0} (\lambda_{\exp(t X)} \phi)(x) =  \frac{d}{dt}\bigg|_{t=0} \phi(\exp(- t X) \cdot x) \]
for all $X \in \mf{g}$, $x \in M$, $\phi \in C_c^\infty(M)$.

We interpret the Riemannian metric $\langle \cdot, \cdot \rangle$ as an inner product on the Lie algebra $\mf{g} = \mf{sl}_2(\RR)$ (of left-invariant vector fields), and pick an orthonormal basis $E_1, E_2, E_3$. 
Given a multi-index $\alpha = (i_1,\ldots,i_d) \in \{1,2,3\}^d$ of degree $\abs{\alpha} \coloneqq d$ we may define
\[ E_\alpha \cdot \phi \coloneqq E_{i_1} E_{i_2}\cdots E_{i_d} \cdot \phi\]
for every $\phi \in C_c^\infty(M)$. 
We use the convention $\alpha = \emptyset$ iff $\abs{\alpha} = 0$, and define
$E_\emptyset \cdot \phi = \phi.$

This allows us to define the \emph{(degree $d$) Sobolev norm}
\[ \Sob_d(\phi) \coloneqq \sum_{0 \leq \abs{\alpha} \leq d} \|E_\alpha \cdot \phi\|_{2} \]
for all $\phi \in C_c^\infty(\Gamma \lquot G)$. The \emph{(degree $d$) Sobolev space $H^d(\Gamma \lquot G)$} is by definition the completion of $C_c^\infty(\Gamma \lquot G)$ with respect to $\Sob_d$.

The following version of the Sobolev Embedding Theorem applies.
\bthm[{\cite[Theorem 2.10]{aubin}}] \label{thm:SET}
    If $(d-r)/3 > 1/2$ then $H^d(\Gamma \lquot G) \subseteq C^r(\Gamma \lquot G)$ and the identity operator is continuous. Here $r \geq 0$ is an integer and $C^r(\Gamma \lquot G)$ is the space of $C^r$-functions with norm $\|\phi\|_{C^r} \coloneqq \max_{0 \leq \abs{\alpha} \leq r} \| E_\alpha \cdot \phi\|_\infty$, $\phi \in C^r(\Gamma \lquot G)$.
\ethm
\bcor \label{cor:SET}
    In particular, if the degree $d=3$ and $r=1$ then there is a \emph{Sobolev constant} $\Ksob >0$ such that 
    \[ \| \phi \|_\infty \leq \norm{\phi}_{C^1} \leq \Ksob \cdot \Sob(\phi) \]
    for all $\phi \in H^3(\Gamma \lquot G)$, where we dropped the degree $d=3$ in $\Sob(\phi) = \Sob_3(\phi)$.
\ecor

Recall that $L^1(G)$ is a Banach algebra when we define multiplication by convolution:
\begin{align*}
    (f_1 * f_2)(g) 
    &\coloneqq \int_G f_1(h) \cdot f_2(h^{-1} g) \, d\nu(h) \\
    &= \int_G f_1(g h) \cdot f_2(h^{-1}) \, d\nu(h) \qquad \forall g \in G \quad \forall f_1, f_2 \in L^1(G)
\end{align*}
There is a Banach algebra action of $L^1(G)$ on $L^2(\Gamma \backslash G)$ given by convolution
\begin{align*}
    (\psi * f)(\Gamma g) 
        &\coloneqq \int_G \psi(h) \cdot (\lambda_{h^{-1}}f)(\Gamma g) \, d\nu(h) \\
        &= \int_G \psi(h) \cdot f(\Gamma g h) \, d\nu(h),
\end{align*}
for all $\Gamma g \in \Gamma \backslash G, \psi \in L^1(G), f \in L^2(\Gamma \backslash G)$. An application of Fubini shows that
\[ \| \psi * f \|_2 \leq \|\psi\|_1 \cdot \|f\|_2. \]

Moreover, we have the following lemma familiar from the situation in $\RR^n$.

\blem \label{lem:conv_properties}
    Let $\epsilon>0$, let $\psi \in C_c^\infty(G)$, and let $f \in L^2(\Gamma \backslash G)$.
    
    Then:
    \begin{enumerate}
        \item $\psi * f$ is smooth;
        \item $E_\alpha \cdot (\psi * f) = (E_\alpha \cdot \psi) * f$ for all multi-indices $\alpha$; 
        \item $\supp(\psi * f) \subseteq N_\epsilon(\supp(f))$, if $\supp(\psi) \subseteq B_\epsilon(e)$.
    \end{enumerate}
\elem

\begin{proof}
    \begin{enumerate}
        \item This will follow from (2).
        \item By induction on $\abs{\alpha}$ it is enough to show this for $\abs{\alpha}=1$. Let $E \coloneqq E_i$ be a basis vector. We compute:
        
        \begin{align*}
            (E \cdot (\psi * f))(\Gamma g) 
            &= \frac{d}{dt}\bigg|_{t=0} \int_G \psi(h) \cdot f(\Gamma g \exp(t E) h) \, d \nu(h)\\
            &= \frac{d}{dt}\bigg|_{t=0} \int_G \psi(\exp(-t E) h) \cdot f(\Gamma g h) \, d \nu(h)\\
            &= \int_G (E \cdot \psi)(h) \cdot f(\Gamma g h) \, d \nu(h)\\
            &= ((E\cdot \psi)*f)(\Gamma g)
        \end{align*}
        
        \item Recall that
        \[ (\psi * f)(\Gamma g) = \int_G \psi(h) \cdot f(\Gamma g h) \, d\nu(h). \]
        Thus, if $(\psi * f)(\Gamma g) \neq 0$ then there is $h \in \supp(\psi) \subseteq B_\epsilon(e)$ such that $\Gamma g h = h^{-1} \cdot \Gamma g \in \supp(f)$. Hence,
        \[ d_{\Gamma \backslash G} (\Gamma g h, \Gamma g) \leq d_G(gh,g) = d_G(h,e) < \epsilon. \]
        This shows that $\Gamma g \in N_\epsilon(\supp(f))$.\end{enumerate}\end{proof}
\section{Outline of the Main Theorem}
\label{sec:maintheorem}

We now outline the construction and proof of the main result. We fix a closed hyperbolic surface $X$, with injectivity radius $\rho$ and with an isometric geodesic pants decomposition $\P$. We use Lemma \ref{modelrel} to show that the set of geodesics obtained by the construction we now describe has asymptotic density one (in the sense of Subsection~\ref{subsec:densityone}). Let $\gamma_v[0,t]$ denote the ray starting at $v \in \UT(X)$, parameterised by arc length, and flowing for time $t$ so 
that $d(g_t(v),v)<\epsilon$. We close up $\gamma_v[0,t]$ by connecting its endpoints via a short arc to obtain the broken geodesic path $\gamma_v(t)$. Let 
$\hat{\gamma}_v(t)$ be the corresponding geodesic representative and let $\hat{\gamma}^p_v(t)$ be the associated primitive sub-curve 
of $\hat{\gamma}_v(t)$. We then break the proof of the main Theorem \ref{maintheorem} into the proof of the corresponding statement for the model of geodesic obtained by the geodesic flow which is Theorem \ref{auxmainthm}:

\begin{customthm}{\ref{auxmainthm}}
    Let $\eta > 1 (> \delta > 0)$ and let $F(t) \in O(t^{\delta/2\eta})$. Then for $\mu$-almost every $v\in \UT(X)$ there is $T' = T'(v)$ such that $\hat\gamma_v(t)$ is filling and
\be  F(\ell_X(\hat\gamma_v(t))) \leq \Vol(M_{\hat\gamma_v^p(t)})\ee
    for every $t \in R_\epsilon(v) \cap [T'(v),+\infty)$. Moreover, as long as $F$ is increasing we obtain that $ F(\ell_X(\hat\gamma^p_v(t)))\leq \Vol(M_{\hat\gamma_v^p(t)})$.
\end{customthm}

\begin{proof}[Sketch of proof]

\begin{enumerate}
    \item Given a pair of pants decomposition $\mathcal{P}$ of $S$, recall the decomposition of $\UT(X)$ from Definition~\ref{def:decomposition}.
    We consider the orthogeodesics $o^j_i$ contained in some pair of pants $P^j$ and define the set $U^j_i$ of all unit tangent vectors $v \in \UT(P^j)$
 whose associated geodesic is in the free homotopy class of $o^j_i$ (see Subsection~\ref{subsubsect:orthogeodsets}).
 
    Section~\ref{sec:counting} estimates the number of regions $U^j_i$ visited by a ray $g_{[0,t]}(v)$ of length $t$. We show that for every $\eta > 1$ and any (sub-linear) function 
    \[ F(t) \in O\left(t^{\delta/2\eta}\right)\] 
    this number grows for $\mu$-almost every $v \in \UT(X)$ at least like $F(t)$
    for $t$ large enough, where $0<\delta<1$ is the Hausdorff dimension of the limit set of any of the constituting (isometric) pairs of pants in $X$.
    
    \item From our construction of the curves $\hat\gamma_v(t)$, $t \in R_\epsilon(v)$, by closing up $\gamma_v[0,t]$ in section~\ref{sec:random} we get that
    \[
       t- \rho \leq \ell( \hat{\gamma}_v(t) ) \leq t + \rho.
    \]
    Moreover, we have that the primitive geodesic representative $\hat{\gamma}_v^p(t)$ satisfies
    \[
    \ell( \hat{\gamma}^p_v(t)) \leq  \ell( \hat{\gamma}_v(t) ).
    \]
    \item Section~\ref{sec:bigons} shows that:
    \[
      \# \{ \gamma_v[0,t]\mbox{-arcs} \} - 6 \leq \# \{ \hat{\gamma}_v(t)\mbox{-arcs} \} = \# \{\hat{\gamma}^p_v(t)\mbox{-arcs} \}.
    \]
    \item By (3) and Section~\ref{sec:counting}, we obtain that asymptotically
    \begin{align*}
        \f{v_3}2 F(t)-3v_3 &\leq\f{v_3}2 \# \{ \gamma_v[0,t]\mbox{-arcs} \} - 3v_3 \leq\f{v_3}2 \# \{ \hat{\gamma}_v(t)\mbox{-arcs} \}\\ 
        &=\f{v_3}2 \#\{\hat{\gamma}^p_v(t)\mbox{-arcs} \} \leq \Vol(M_{{\hat{\gamma}_v^p(t)}})+C(g),
    \end{align*}
    where $C(g)$ is some constant only depending on the topological complexity of $X$.

    \item Because $t-\rho\leq \ell_X(\hat\gamma_v)\leq t+\rho$ and $F(t) \in O\left(t^{\delta/2\eta}\right)$ there exists a $B>0$ such that:
   
   \[ \f{v_3}2F(t+\rho)-B\leq \f{v_3}2 F(t)-3v_3,
   \]
  Thus, for large $t$ we have:
   \[ \f{v_3}2 F(\ell_X({\hat{\gamma}_v(t)}))-B\leq \f{v_3}2F(t+\rho)-B\leq \f{v_3}2 F(t)-3v_3\leq \Vol(M_{\hat{\gamma}^p_v(t)})+C(g).
   \]
   \item Putting all the steps together, we obtain that
   \[
  O\tonde{   \ell_X(\hat\gamma_v(t))^{\delta/2\eta}} \leq \Vol(M_{\hat{\gamma}^p_v(t)})
   \]
   with $0<\delta<1,$ and $\eta>1$.
\end{enumerate} \end{proof}

Then the proof of the main Theorem \ref{maintheorem} follows from Theorem \ref{auxmainthm}, Lemma \ref{modelrel} and the fact that primitive geodesics have density one:

\begin{customthm}{\ref{maintheorem}}
    Let $\eta > 1 (> \delta > 0)$ and let $F(t) \in O(t^{\delta/2\eta})$. 
    
    Then the collection of filling, primitive closed geodesics satisfying
    \[  F(\ell_X(\gamma)) \leq \Vol(M_{\gamma}) \]
    has asymptotic density one.
\end{customthm}

\section{Random sequences of closed geodesics and asymptotic density}
\label{sec:random}

In this section we will discuss a natural way of constructing sequences of closed geodesics. Similar constructions have already been considered in \cite{Bonahon88:GeodesicCurrent},\cite{Bri04:Random} and \cite{Sap17:Birman}. Finally, we will see how this is related to the notion of asymptotic density.

\subsection{ Closing-up along geodesics }

We fix $0 < \rho < \tfrac{1}{2} \inj(X)$, and choose $\epsilon = \epsilon(\rho) > 0$ and $T = T(\rho) > 0$ as in the Anosov Closing Lemma \ref{lem:AnosovClosing}. Further, let $v \in \UT(X)$ and denote by $R_\epsilon(v) = \{ t > 0 \, | \, d(g_t(v),v) < \epsilon \}$ its $\epsilon$-return times as before. Therefore, for every $\epsilon$-return time $t \in R_\epsilon(v) \cap [T,+\infty)$ after $T=T(\rho) > 0$ we find a periodic vector $v' \in B_\rho(v)$ with period $t' \in [t - \rho, t + \rho]$. Moreover, the geodesics $\gamma, \gamma'$ with respective starting vectors $\dot\gamma(0)=v, \dot\gamma'(0) = v'$ stay close to each other:
\[d(\gamma(s \cdot t), \gamma'(s \cdot t')) \leq 2 \rho < \inj(X)\]
for every $s \in [0,1]$. Note that $R_\epsilon(v) \subseteq \RR$ is unbounded for $\mu$-almost every vector $v \in \UT(X)$ by ergodicity of the geodesic flow. 

\brem We will need, in Section \ref{sec:bigons}, that $\rho$ is smaller than the minimal distance $o$ between any two pants curves in $\P$. However, this is easily achieved by shrinking the cuff length $c$ in the isometric pants decomposition. Indeed, $0 < \rho < \tfrac{1}{2} \inj(X) \to 0$ and $o \to + \infty$ as $c \to 0$.
\erem

The closed geodesic $\gamma'$ can be characterized as follows. We form a broken geodesic curve $\gamma_v(t)$ by concatenating the geodesic segment $\gamma|_{[0,t]}$ with the short geodesic arc $\alpha$ connecting the endpoint $\gamma(t)$ to the starting point $\gamma(0)$. Observe that $\ell(\alpha) \leq \rho$ because $d(v,v') \leq \rho < \tfrac{1}{2} \inj(X)$. By construction $\gamma_v(t)$ is freely homotopic to $\gamma'$. Therefore, $\gamma'$ is the unique closed geodesic in the free homotopy class of $\gamma_v(t)$, i.e.\ one obtains $\gamma'$ by pulling tight $\gamma_v(t)$.

\bdefi\label{return} 
    We will use the notation $\hat{\gamma}_v(t) \coloneqq \gamma'$. Moreover, we will denote by $\hat\gamma^p_v(t)$ its primitive subcurve. With an abuse of notation we will use $M_{\hat\gamma^p_v(t)}$ to denote $\PT(X)$ minus the canonical lift of $\hat\gamma^p_v(t)$.
\edefi

In this way we obtain for $\mu$-almost every unit tangent vector $v \in \UT(X)$ an infinite collection of geodesics $\{ \hat\gamma_v(t) \}_{t \in R_\epsilon(v) \cap [T, + \infty)}$. We will study the asymptotic properties of these closed-up geodesics as $t$ tends to infinity in the following.

\brem Because periodic orbits of the geodesic flow are dense (see \cite[Theorem~3.3]{D11:Geodesic}) we have that, in general, the curve $\gamma_v(t)$ will not be primitive. This is because we can find vectors $w$ arbitrarily close to a periodic orbit $\beta$, i.e.\ they track the period orbit for a long time $t$. Thus, if $\beta$ is filling we can obtain arbitrarily long curves $\hat{\gamma}_w(t)$ that are homotopic to arbitrarily high powers of $\beta$.
\erem

Let us first show that $\hat{\gamma}_v(t)$ converges to a filling geodesic current:

\blem \label{lem:conv_to_Liouville}
For $\mu$-almost every $v \in \UT(X)$ the closed geodesics $\hat\gamma_v(t)$ satisfy
\[
\lim_{t \to \infty} \frac{\hat{\gamma}_v(t)}{\ell_X(\hat{\gamma}_v(t))} = \mathcal{L}_X
\]
as geodesic currents, where  $t \in R_\epsilon(v) \cap [T,+\infty)$ and  $\mathcal{L}_X$ denotes the Liouville current associated to the hyperbolic structure $X$ (see definitions in Subsection~\ref{subsec:geometry})
\label{lem:modeltoliouville}
\elem

\bpf
By the Birkhoff Ergodic Theorem \ref{birkhoff} it suffices to prove that for every continuous function $f \colon \UT(X) \to \mathbb{R}$,

\[\lim_{t \to \infty} \frac{1}{\ell_X(\hat\gamma_v(t))} \int_{\UT(X)} f d \hat\gamma_v(t) =  \lim_{t \to \infty} \frac{1}{t} \int_0^t f(g_s(v))ds.\]

Take a transversal $D$ to the leaf of the geodesic foliation the current $\hat\gamma_v(t)/\ell_X(\hat\gamma_v(t))$ is supported on, intersecting at points $x_1(t),\cdots, x_k(t)$.
Let $y_1(t),\cdots,y_k(t)$ be the intersection points of the leaf corresponding to $g_s(v)$ for $s \in (-\infty,\infty)$.
We claim that for $t$ long enough, $y_i(t)$ and $x_i(t)$ are close.
One way to see this is to look at the endpoints at infinity in the bands model, and note that the endpoints of the closed broken geodesic $\gamma_v(t)$ and the corresponding closed geodesic $\hat\gamma_v(t)$ are the same, and are in a neighborhood of the boundary at infinity of the endpoints of $\gamma_v(t)$ of radius $\epsilon$ that goes to $0$ as $t$ goes to infinity. 
\epf

As a result we obtain that for large enough $t \in R_\epsilon(v) \cap [T, + \infty)$ the geodesics $\hat\gamma_v(t)$ is filling:

\blem\label{curveisfilling}
For $\mu$-almost every $v \in \UT(X)$ there is $t^* \coloneqq t^*(v) \geq T(\rho) > 0$ such that the curve $\hat\gamma_v(t)$ is filling for every $t \in R_\epsilon(v) \cap [t^*,+\infty)$.
\elem
\bpf
For a geodesic current $\mu$, let $\sys \coloneqq \inf \{ i(\mu,c) : c \mbox{ closed curve } \}$ be the systole of $\mu$. 
By Lemma~\ref{lem:modeltoliouville}, $\lim_{t \to \infty} \hat\gamma_v(t)/\ell(\hat\gamma_v(t)) \to \mathcal{L}_X$. Since $\mathcal{L}_X$ is filling, by \cite[Theorem~1.3]{BIPP20:SystolesCurrents}), we have $\sys(\mathcal{L}_X)>0$. The result then follows since $\sys \colon \GC(S) \to \mathbb{R}_{\geq 0}$ is a continuous function on currents (\cite[Corollary~1.5]{BIPP20:SystolesCurrents}).
\epf

\bdefi[Filling time] \label{def:filling_time}
We call the above time $t^*(v)$ the \emph{filling time} of $v$.
\edefi

We conclude the section by observing that, by ergodicity, it follows that by closing-up along almost every geodesic the volume of the corresponding complements goes to infinity.

\blem For $\mu$-almost every $v\in \UT(X)$ we have that $\text{Vol}(M_{ \hat\gamma_v(t)}) \to \infty$ as $t \to \infty$.
\elem
\bpf 

Consider the decomposition of $\UT(X)$ into the subsets $U_i^j \subset \UT(X)$ of vectors corresponding to the homotopy classes of all orthogeodesics from Definition~\ref{def:decomposition}. Therefore, for a closed geodesic $\gamma$ the number of different homotopy classes of (undirected) $\gamma$-arcs is half the number of different sets $\{U_i^j\}_{i,j}$, that $\dot\gamma(t), t \in [0,\ell(\gamma)],$ visits.

By Birkhoff's Ergodic Theorem \ref{thm:birkhoff},
\[ \lim_{T \to \infty} \frac{1}{T} \int_0^T \chi_{U^j_i}(g_t(v)) dt = \int_{\UT(X)} \chi_{U^j_i}(w) d \mu(w) = \mu(U^j_i) > 0 \]
for $\mu$-almost every $v \in \UT(X)$ and all $i,j$. Consequently, the geodesic ray $\gamma_v[0,T]$ in $\mu$-almost every direction $v \in \UT(X)$ will eventually intersect every subset $U^j_i$ such that $\gamma_v[0,T]$ will eventually have a subarc in every homotopy class $o^j_i$. 

Later, in section \ref{sec:bigons}, we will see that the number of homotopy classes of subarcs of $\gamma_v[0,T]$ is roughly the number of homotopy classes of subarcs of $\hat{\gamma}_v(T)$; see Proposition \ref{prop:arcs_pulled_tight}. By the previous argument this number goes to infinity as $T$ tends to infinity. It follows from Theorem \ref{thm:lowerbound2} that
\[ \Vol(M_{\hat{\gamma}_v(T)}) \to + \infty \qquad (T \to + \infty).\]
\epf

\subsection{Relation to asymptotic density}

Let us see now how the previous construction of closing-up along a geodesic is related to the notion of asymptotic density:

Let $\mc{A} \subseteq \mc{G}$ be a collection of closed geodesics. Consider the set $\Omega \subseteq \UT(X)$ of starting vectors whose geodesic segments eventually close-up to a geodesic in $\mc{A}$:
\[ \Omega_{\mc{A}} \coloneqq \{ v \in \UT(X) \, | \, \exists T=T(v) > 0 \forall t \in R_\epsilon(v) \cap [T,+\infty): \hat{\gamma}_v(t) \in \mc{A} \}.\]

\blem\label{modelrel}
    A collection $\mc{A} \subseteq \mc{G}$ has asymptotic density one  if the set $\Omega_{\mc{A}}$ has full measure.
\elem

\begin{proof}
    For any subset $\mc{B} \subseteq \mc{G}$ we will use the following notation
    \begin{align*}
        \mc{B}[T,T') &\coloneqq \{ \gamma \in \mc{B} \, | \, T \leq \ell(\gamma) < T' \},\\
        \Omega_{\mc{B}}[T,T') &\coloneqq \{ v \in \UT(X) \, | \,\hat{\gamma}_v(t) \in \mc{B} \quad \forall t \in R_\epsilon(v) \cap [T,T') \}
    \end{align*}
    for all $0<T<T'$.

    For every closed geodesic $\gamma \in \mc{G}$ we consider the set $N(\gamma)$ of vectors $v \in \UT(X)$ such that there is $t_0 \in \RR$ with $d(v,\dot\gamma(t_0)) < \epsilon/2$ and $d(g_t(v),\dot\gamma(t+t_0)) < \epsilon/2$ for all $t \in [0, \ell(\gamma)]$. In particular, for every $v \in N(\gamma)$ we have that $t \coloneqq \ell(\gamma) \in R_\epsilon(v)$ and the broken geodesic curve $\gamma_v(t)$ is freely homotopic to $\gamma$, whence $\hat\gamma_v(t) = \gamma$. 
    
    If $v \in N(\gamma_1) \cap N(\gamma_2)$ and $|\ell(\gamma_1) - \ell(\gamma_2)| < \epsilon < \tfrac{1}{2} \inj(X)$ then $\gamma_1 = \gamma_2$, for any $\gamma_1, \gamma_2 \in \mc{G}$. Indeed, both broken geodesic curves $\gamma_v(\ell(\gamma_1))$ and $\gamma_v(\ell(\gamma_2))$ are in the same free homotopy class. Because there is only one closed geodesic in any free homotopy class we obtain 
    \[\gamma_1 = \hat\gamma_v(\ell(\gamma_1)) = \hat\gamma_v(\ell(\gamma_2)) = \gamma_2.\]
    
    Moreover, a direct computation using the flowbox coordinates and the Anosov properties of the flow yields that there is a universal constant $C>0$ such that 
    \[\mu(N(\gamma)) \geq C \cdot \epsilon \cdot \ell(\gamma) \cdot e^{-\ell(\gamma)}.\]
   We have that
    \[ \Omega_{\mc{A}} = \bigcup_{T>0} \Omega_{\mc{A}}[T, + \infty) \]
    and $\Omega_{\mc{A}}[T,\infty) \subseteq \Omega_{\mc{A}}[T',\infty)$ for all $0<T\leq T'$. Thus,
    \[ \mu(\Omega_{\mc{A}}[T,\infty)) \to \mu(\Omega_{\mc{A}}) = 1 \quad (T \to \infty).  \]
    
    Suppose $w \in N(\gamma)$ for $\gamma \in \mc{A}^c[k \cdot \epsilon, (k+1) \cdot \epsilon)$, $k \in \NN$. Then $\ell(\gamma) \in R_\epsilon(w) \cap [k \cdot \epsilon, (k+1) \epsilon)$ and $\gamma = \hat\gamma_w(\ell(\gamma)) \notin \mc{A}$. Thus, $w \notin \Omega_{\mc{A}}[k \cdot\epsilon, + \infty)$. Hence,
    \[ \bigsqcup_{\gamma \in \mc{A}^c[k \cdot \epsilon, (k+1) \cdot \epsilon)} N(\gamma) \subseteq \Omega_{\mc{A}}[k \cdot\epsilon, + \infty)^c, \]
    and
    \[ \sum_{\gamma \in \mc{A}^c[k \cdot \epsilon, (k+1) \cdot \epsilon)} \mu(N(\gamma)) \leq \mu(\Omega_{\mc{A}}[k \cdot\epsilon, + \infty)^c) \eqqcolon \alpha_k \to 0 \quad (k \to \infty). \]
    Note that $\alpha_k \geq \alpha_{k+1}$ for all $k \in \NN$ and $0 \leq \alpha_k \leq 1$.
    
    On the other hand,
    \begin{align*}
        \sum_{\gamma \in \mc{A}^c[k \cdot \epsilon, (k+1) \cdot \epsilon)} \mu(N(\gamma)) &\geq \sum_{\gamma \in \mc{A}^c[k \cdot \epsilon, (k+1) \cdot \epsilon)} C \cdot \epsilon \cdot \ell(\gamma) e^{-\ell(\gamma)} \\
    &\geq C \cdot (k+1) \epsilon \cdot  e^{-(k+1) \epsilon} \cdot \#\mc{A}^c[k \cdot \epsilon, (k+1) \cdot \epsilon).
    \end{align*}

    Thus,
    \begin{align}
        \#\mc{A}^c[k \cdot \epsilon, (k+1) \cdot \epsilon) \leq \frac{e^{(k+1) \epsilon}}{C \cdot (k+1) \cdot \epsilon} \cdot \alpha_k. \label{est:Ac}
    \end{align}
    
    We need to show that
    \begin{align}
        \frac{\# \mc{A}(R)}{\# \mc{G}(R)} = \frac{\# \mc{A}[0,R)}{\# \mc{G}[0,R)} = 1 - \frac{\# \mc{A}^c[0,R)}{\# \mc{G}[0,R)} \to 1 \quad (R \to \infty). \label{eq:limitgoal}
    \end{align}
    
    Let $n \in \NN$ such that $n \epsilon \leq R < (n+1) \epsilon$. By (\ref{est:Ac}), we obtain
    \begin{align*}
        \#\mc{A}^c[0,R) &\leq \sum_{k=0}^n \#\mc{A}^c[k \epsilon,(k+1) \epsilon)
        \leq\frac{1}{C} \sum_{k=0}^n \frac{\alpha_k}{(k+1)\epsilon} e^{(k+1) \epsilon} \\
        &\leq \frac{1}{C} \sum_{k=0}^{\floor{n/2}} \frac{\alpha_k}{(k+1)\epsilon} e^{(k+1) \epsilon} 
        + \frac{1}{C} \sum_{k=\floor{n/2}+1}^{n} \frac{\alpha_k}{(k+1)\epsilon} e^{(k+1) \epsilon}.
    \end{align*}
    
    We can estimate the first summand as follows
    \[ \frac{1}{C} \sum_{k=0}^{\floor{n/2}} \frac{\alpha_k}{(k+1)\epsilon} e^{(k+1) \epsilon} \leq \frac{1}{C \epsilon} \sum_{k=0}^{\floor{n/2}} e^{(k+1) \epsilon} \ll e^{\floor{n/2} \epsilon}   \ll e^{R/2} \in o \kl{\frac{e^R}{R}}.\]
    
    Regarding the second we get
    \[ \frac{1}{C} \sum_{k=\floor{n/2}+1}^{n} \frac{\alpha_k}{(k+1)\epsilon} e^{(k+1) \epsilon} \ll \frac{\alpha_{\floor{n/2}}}{(\floor{n/2}+2)\epsilon} \sum_{k=\floor{n/2}+1}^{n} e^{(k+1)\epsilon} \ll \alpha_{\floor{n/2}} \frac{e^R}{R} \in o \kl{\frac{e^R}{R}},\]
    because $\alpha_{\floor{n/2}} \to 0$ as $R \to \infty$.
    
    By a result of Huber's~\cite{Hub60:Growth} it is known that 
    \[\# \mc{G}[0,R) \sim \frac{e^R}{R}\] 
    as $R \to \infty$; see Remark \ref{rmk:margulis_count}. Therefore, the previous two estimates imply
    \[ \frac{\# \mc{A}^c[0,R)}{\# \mc{G}[0,R)} \to 0 \quad (R \to \infty)\]
    proving (\ref{eq:limitgoal}). This concludes the proof.
\end{proof}

\section{Arcs of segment vs arcs of closed-up geodesic representative}
\label{sec:bigons}

Let $\P$ be a geodesic pants decomposition of $X$ and denote by $P_1, \ldots, P_m$ the different pairs of pants. We will be interested in the free homotopy types of \emph{directed} arcs connecting two (possibly the same) boundary components of a pair of pants. Here, the free homotopy is allowed to move the arc's endpoints along the respective boundary component.

It will be useful to have a characterization of these homotopy types in terms of the universal covering $\pi \colon \HH^2 \to X = \Gamma \lquot \HH^2$. 
To this end note that the bi-infinite geodesics in the preimage $\pi^{-1}(\P)$ divide $\HH^2 = \bigcup_{i=1}^m \bigcup_{j\in\NN} \tilde{P}^j_i$ into convex domains $\{\tilde{P}^j_i\}_{j \in \NN}$ each covering a pair of pants $P_i \subseteq X$, $i=1,\ldots,m$. 
Given a properly immersed directed arc $\alpha$ in a pair of pants $P_i$ we may consider a lift $\tilde{\alpha} \subseteq \tilde{P}^j_i$ in some domain $\tilde{P}^j_i$. 
The lift $\tilde{\alpha}$ has its starting point (resp.\ end point) on the unparametrized lift $\tilde{\alpha}_-$ (resp.\ $\tilde{\alpha}_+$) of some pants curve. 
It turns out that the tuple $(\tilde{\alpha}_-, \tilde{\alpha}_+)$ determines the free homotopy type of $\alpha$ uniquely up to deck transformations:

\blem \label{lem:charact_htpy_arcs}
    The map
    \begin{align*}
        \{ \alpha \subseteq X \text{ a directed arc in some } P_i \}/\text{free homotopy} &\longrightarrow \Gamma \lquot (\pi^{-1}(\P) \times \pi^{-1}(\P)), \\
        [\alpha] & \longmapsto [\tilde{\alpha}_-, \tilde{\alpha}_+] \coloneqq \Gamma \cdot (\tilde{\alpha}_-, \tilde{\alpha}_+),
    \end{align*}
    is well-defined and injective.
\elem

Let us denote the geodesic segment of length $t>0$ starting at the unit tangent vector $v \in \UT(X)$ by $\gamma_v[0,t]$. Observe that for $\mu$-almost every $v \in \UT(X)$ the corresponding bi-infinite geodesics $\gamma_v(-\infty,\infty)$ intersects every pants curve in $\P$ transversely. In this case we can extend $\gamma_v[0,t]$ to the minimal geodesic segment $\gamma_v[-a,b]$, $a,b>0$, whose end points lie on some pants curves. This segment intersects the pants curves $\P$ at times $-a = t_0 < t_1 < \ldots < t_n = b$ and we obtain corresponding (directed) subarcs $\gamma_i \coloneqq \gamma_v[t_{i-1},t_i]$, $i=1,\ldots,n$ each of which is contained in some pair of pants. 

\bdefi
    A directed subarc $\gamma_i$ as above is called a \emph{ directed $\gamma$-arc}.
    We denote by
    \[ A_\P(\gamma_v[0,t]) \coloneqq \# \{ [\gamma_i] \, | \, i=1,\ldots,n \} \]
    the number of different free homotopy classes of directed $\gamma$-arcs. 
    
    Since a closed geodesic $\gamma$ intersects the geodesic pants decomposition $\P$ transversely we obtain, without the need of an extension, a similar decomposition into directed $\gamma$-arcs. In the same fashion, we denote by $A_\P(\gamma)$ the number of different free homotopy classes of its directed $\gamma$-arcs as well.
\edefi

\brem \label{rem:counting_directed}
    Note that we consider \emph{directed} $\gamma$-arcs here instead of undirected $\gamma$-arcs as in Theorem \ref{thm:lowerbound2}; see Remark \ref{rem:undirected}. Since every homotopy class of an undirected $\gamma$-arc corresponds to two homotopy classes of directed $\gamma$-arcs we have
    \[ \sum_Q \#{\set{\gamma\text{-arcs in }Q}} = \frac{1}{2} A_\P(\gamma),  \]
    where we sum over all components of the pants decomposition $\P$.
\erem

On a pair of pants the simple (undirected) $\gamma$-arcs are described in Figure \ref{fig:simple}.

In section \ref{sec:random} we discussed how to construct sequences of closed geodesics by closing-up along geodesics: For $v \in \UT(X)$ and an $\epsilon$-return time $t \in R_\epsilon(v) \cap [T,+\infty)$ we may close-up the geodesic segment $\gamma_v[0,t]$ and consider the unique geodesic curve $\hat{\gamma}_v(t)$ in its homotopy class. Using Lemma \ref{lem:charact_htpy_arcs} we will show that $A_\P(\hat{\gamma}_v(t))$ is equal to $A_\P(\gamma_v[0,t])$ up to a uniform error: 

\bprop \label{prop:arcs_pulled_tight}
    For $\mu$-almost every $v \in \UT(X)$ and every $t \in R_\epsilon(v) \cap [t^*(v),+\infty)$
    \[ \abs{A_\mc{P}(\gamma_v[0,t]) - A_\mc{P}(\hat\gamma_v(t)) } \leq 6 \]
    holds, where $t^*(v)$ denotes the filling time of $v$ (see Lemma \ref{curveisfilling}).
\eprop

\begin{proof}
Let $\tilde{v} \in \UT(\HH^2)$ be a preimage of $v \in \UT(X)$. By the Anosov Closing Lemma \ref{lem:AnosovClosing} we find  $\tilde{v}' \in \UT(\HH^2)$, $t' \in [t-\rho,t+\rho]$ and $g \in \Gamma$ such that $g_{t'}(\tilde{v}') = g \cdot \tilde{v}'$ and
\begin{equation} \label{eq:dist_ends_vec}
    d(\tilde{v}',\tilde{v}) = d(g_{t'}(\tilde{v}'),g \cdot \tilde{v}) < \rho.
\end{equation}

By definition the geodesic segment $\tilde{\gamma}'$ of length $t'$ starting from $\tilde{v}'$ is a lift of $\gamma' \coloneqq \hat\gamma_v(t)$. 
Likewise, the geodesic segment $\tilde{\gamma}$ of length $t$ starting from $\tilde{v}$ is a lift of $\gamma_v[0,t]$. 
Let $\tilde{\alpha} \coloneqq (\tilde\alpha_1, \ldots, \tilde\alpha_l) \subseteq \pi^{-1}(\P)$ (resp.\ $\tilde{\beta} \coloneqq (\tilde\beta_1, \ldots, \tilde\beta_m) \subseteq \pi^{-1}(\P)$) be the sequence of lifts of pants curves that $\tilde{\gamma}$ (resp.\ $\tilde{\gamma}'$) intersects in this order.

The following observation is key for our proof.

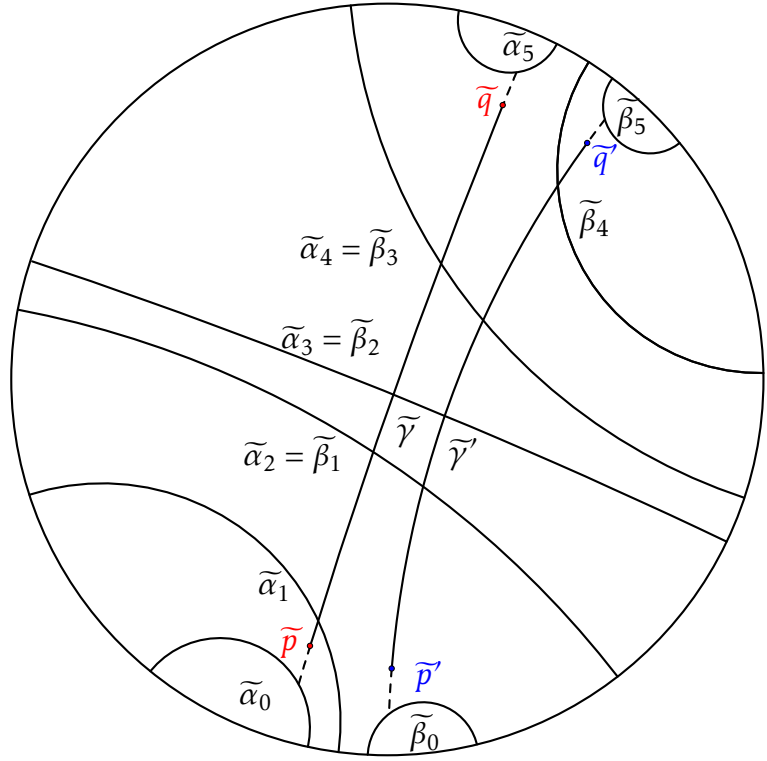
\begin{figure}
    \centering
        \definecolor{ffqqqq}{rgb}{1,0,0}
        \definecolor{qqqqff}{rgb}{0,0,1}
        \begin{tikzpicture}[line cap=round,line join=round,>=triangle 45,x=5cm,y=5cm]
        \clip(-1.14761428275442,-1.3904255811339519) rectangle (1.2945888360291156,1.2254761947465724);
        \draw [line width=0.8pt] (0,0) circle (1);
        \draw [shift={(0.9906621204897947,0.5560739237019391)},line width=0.8pt]  plot[domain=2.576722263836339:4.729441334215951,variable=\t]({1*0.5391007749898026*cos(\t r)+0*0.5391007749898026*sin(\t r)},{0*0.5391007749898026*cos(\t r)+1*0.5391007749898026*sin(\t r)});
        \draw [shift={(-0.7575213549117695,-0.9094291924783569)},line width=0.8pt]  plot[domain=-0.13006954555393513:1.8826222095320857,variable=\t]({1*0.6331666915427561*cos(\t r)+0*0.6331666915427561*sin(\t r)},{0*0.6331666915427561*cos(\t r)+1*0.6331666915427561*sin(\t r)});
        \draw [shift={(-0.44473670499791107,-0.9269854658597059)},line width=0.8pt]  plot[domain=-0.2127869031284808:2.459718959973281,variable=\t]({1*0.23894097741395545*cos(\t r)+0*0.23894097741395545*sin(\t r)},{0*0.23894097741395545*cos(\t r)+1*0.23894097741395545*sin(\t r)});
        \draw [shift={(0.09498756581881354,-1.0063519935820653)},line width=0.8pt]  plot[domain=0.2405888357682369:3.089222215122676,variable=\t]({1*0.14753634347774663*cos(\t r)+0*0.14753634347774663*sin(\t r)},{0*0.14753634347774663*cos(\t r)+1*0.14753634347774663*sin(\t r)});
        \draw [shift={(0.3266180880919779,0.9553448593067232)},line width=0.8pt]  plot[domain=2.950429216179031:5.815494306716082,variable=\t]({1*0.13915162835066278*cos(\t r)+0*0.13915162835066278*sin(\t r)},{0*0.13915162835066278*cos(\t r)+1*0.13915162835066278*sin(\t r)});
        \draw [shift={(0.698512640485078,0.7264446229040732)},line width=0.8pt]  plot[domain=2.5002147042750096:5.3929658677382335,variable=\t]({1*0.1250667784172839*cos(\t r)+0*0.1250667784172839*sin(\t r)},{0*0.1250667784172839*cos(\t r)+1*0.1250667784172839*sin(\t r)});
        \draw [shift={(0.9906621204897947,0.5560739237019391)},line width=0.8pt]  plot[domain=2.576722263836339:4.729441334215951,variable=\t]({1*0.5391007749898026*cos(\t r)+0*0.5391007749898026*sin(\t r)},{0*0.5391007749898026*cos(\t r)+1*0.5391007749898026*sin(\t r)});
        \draw [shift={(1.4338725417236902,1.1454104698362195)},line width=0.8pt]  plot[domain=3.239345676552821:4.391885413811152,variable=\t]({1*1.5388163016811278*cos(\t r)+0*1.5388163016811278*sin(\t r)},{0*1.5388163016811278*cos(\t r)+1*1.5388163016811278*sin(\t r)});
        \draw [shift={(-6.055739018896399,-15.002916638413224)},line width=0.8pt]  plot[domain=1.1253101496108355:1.2490061606119074,variable=\t]({1*16.14804888288863*cos(\t r)+0*16.14804888288863*sin(\t r)},{0*16.14804888288863*cos(\t r)+1*16.14804888288863*sin(\t r)});
        \draw [shift={(-1.472713283982658,-2.405357639119981)},line width=0.8pt]  plot[domain=0.6589644372629591:1.3838518253104457,variable=\t]({1*2.6371632086186527*cos(\t r)+0*2.6371632086186527*sin(\t r)},{0*2.6371632086186527*cos(\t r)+1*2.6371632086186527*sin(\t r)});
        \draw [shift={(16.962980382607647,-6.011803842335925)},line width=0.8pt]  plot[domain=2.7569866193251356:2.8420193747582805,variable=\t]({1*17.968986863467197*cos(\t r)+0*17.968986863467197*sin(\t r)},{0*17.968986863467197*cos(\t r)+1*17.968986863467197*sin(\t r)});
        \draw [shift={(2.7563752780688997,-0.9936638062425178)},line width=0.8pt]  plot[domain=2.5114118294543673:3.0599482235316935,variable=\t]({1*2.754082866107295*cos(\t r)+0*2.754082866107295*sin(\t r)},{0*2.754082866107295*cos(\t r)+1*2.754082866107295*sin(\t r)});
        \draw [shift={(16.962980382651295,-6.011803842353901)},line width=0.8pt, dashed]  plot[domain=2.7517569540704847:2.7569866193251196,variable=\t]({1*17.968986863514402*cos(\t r)+0*17.968986863514402*sin(\t r)},{0*17.968986863514402*cos(\t r)+1*17.968986863514402*sin(\t r)});
        \draw [shift={(16.96298038219711,-6.01180384221086)},line width=0.8pt,dashed]  plot[domain=2.842019374758188:2.8478866460236256,variable=\t]({1*17.968986863038037*cos(\t r)+0*17.968986863038037*sin(\t r)},{0*17.968986863038037*cos(\t r)+1*17.968986863038037*sin(\t r)});
        \draw [shift={(2.756375278067005,-0.9936638062424064)},line width=0.8pt,dashed]  plot[domain=3.059948223531678:3.103952076471879,variable=\t]({1*2.754082866105398*cos(\t r)+0*2.754082866105398*sin(\t r)},{0*2.754082866105398*cos(\t r)+1*2.754082866105398*sin(\t r)});
        \draw [shift={(2.7563752780701125,-0.9936638062434237)},line width=0.8pt,dashed]  plot[domain=2.4829872348607624:2.511411829454361,variable=\t]({1*2.754082866108809*cos(\t r)+0*2.754082866108809*sin(\t r)},{0*2.754082866108809*cos(\t r)+1*2.754082866108809*sin(\t r)});
        \draw[color=black] (-0.3,-0.55) node {$\widetilde{\alpha}_1$};
        \draw[color=black] (-0.35,-0.8415378247074938) node {$\widetilde{\alpha}_0$};
        \draw[color=black] (0.1,-0.94) node {$\widetilde\beta_0$};
        \draw[color=black] (0.35,0.8850268268364909) node {$\widetilde{\alpha}_5$};
        \draw[color=black] (0.65,0.6974322771717522) node {$\widetilde{\beta}_5$};
        \draw[color=black] (0.55,0.4264623721004628) node {$\widetilde{\beta}_4$};
        \draw[color=black] (-0.1,0.35) node {$\widetilde{\alpha}_4=\widetilde{\beta}_3$};
        \draw[color=black] (-0.15,0.11) node {$\widetilde{\alpha}_3=\widetilde{\beta}_2$};
        \draw[color=black] (-0.25,-0.2) node {$\widetilde{\alpha}_2 = \widetilde{\beta}_1$};
        \draw [fill=qqqqff] (0.011466414877307596,-0.7690580045409268) circle (1pt);
        \draw[color=qqqqff] (0.11,-0.7911652141493695) node {$\widetilde{p}'$};
        \draw [fill=ffqqqq] (-0.2057128518858029,-0.7089308927626892) circle (1pt);
        \draw[color=ffqqqq] (-0.26,-0.70) node {$\widetilde{p}$};
        \draw [fill=qqqqff] (0.5312939952942188,0.6292920521934268) circle (1pt);
        \draw[color=qqqqff] (0.58,0.584528150058715) node {$\widetilde{q}'$};
        \draw [fill=ffqqqq] (0.30669416814401274,0.7300518696565383) circle (1pt);
        \draw[color=ffqqqq] (0.26,0.7408569414459973) node {$\widetilde{q}$};
        \draw[color=black] (0.05,-0.14) node {$\widetilde{\gamma}$};
        \draw[color=black] (0.2,-0.2) node {$\widetilde{\gamma}'$};
    \end{tikzpicture}
    \caption{This figure shows two geodesic segments coming from two nearby trajectories of the geodesic flow, and indicates the finite sequence of lifts of pants intersected by each one of them, labelled as $\alpha_i$'s and $\beta_i$'s. The figure illustrates that the two sequences of lifts of pant curves intersected by the two geodesic arcs is the same up to possibly the first and last elements. This is the content of the proof of Lemma~\ref{lem:sequ_roughly_coincide}.}
    \label{fig:alt_proof_arc_estimate}
\end{figure}

\blem \label{lem:sequ_roughly_coincide}
    Up to deleting the first or last elements in $\tilde{\alpha}$ or $\tilde{\beta}$ the two sequences coincide; see Figure \ref{fig:alt_proof_arc_estimate}.
\elem

Before providing a proof let us first see how this implies the asserted estimate. By Lemma \ref{lem:charact_htpy_arcs} we have that
\begin{align*}
    A_\P(\gamma_v[0,t]) &= \# \{ [\tilde{\alpha}_0,\tilde{\alpha}_1], \ldots, [\tilde{\alpha}_l,\tilde{\alpha}_{l+1}] \}, \\
    A_\P(\hat\gamma_v(t)) &= \# \{ [\tilde{\beta}_0,\tilde{\beta}_1], \ldots, [\tilde{\beta}_m,\tilde{\beta}_{m+1}] \},
\end{align*}
where $\tilde{\alpha}_0, \tilde{\alpha}_{l+1}$ and $\tilde{\beta}_0, \tilde{\beta}_{m+1}$ are the lifted pants curves that the extensions of $\tilde{\gamma}$ and $\tilde{\gamma}'$ end on, respectively. By Lemma \ref{lem:sequ_roughly_coincide} one obtains shortened sequences $\tilde{\alpha}', \tilde{\beta}'$ by deleting the first or last element from $\tilde{\alpha}, \tilde{\beta}$ such that $\tilde{\alpha}' = \tilde{\beta}'$. 

For example, suppose that $\tilde{\alpha}' = (\tilde{\alpha}_2, \ldots, \tilde{\alpha_l}) = (\tilde{\beta}_1, \ldots, \tilde{\beta}_{m-1}) = \tilde{\beta}'$.
Then 
\[ \# \{ [\tilde{\alpha}_2,\tilde{\alpha}_3], \ldots, [\tilde{\alpha}_{l-1},\tilde{\alpha}_{l}] \} = \# \{ [\tilde{\beta}_1,\tilde{\beta}_2], \ldots, [\tilde{\beta}_{m-2},\tilde{\beta}_{m-1}] \}.\]
Because
\begin{align*}
    \{ [\tilde{\alpha}_0,\tilde{\alpha}_1], \ldots, [\tilde{\alpha}_l,\tilde{\alpha}_{l+1}] \} \triangle \{ [\tilde{\alpha}_2,\tilde{\alpha}_3], \ldots, [\tilde{\alpha}_{l-1},\tilde{\alpha}_{l}] \} &\subseteq \{ [\tilde\alpha_0, \tilde\alpha_1], [\tilde\alpha_1,\tilde\alpha_2], [\tilde\alpha_l, \tilde\alpha_{l+1}] \} \text{ and} \\
    \{ [\tilde{\beta}_0,\tilde{\beta}_1], \ldots, [\tilde{\beta}_m,\tilde{\beta}_{m+1}] \} \triangle \{ [\tilde{\beta}_1,\tilde{\beta}_2], \ldots, [\tilde{\beta}_{m-2},\tilde{\beta}_{m-1}] \} &\subseteq \{ [\tilde\beta_0, \tilde\beta_1], [\tilde\beta_{m-1},\tilde\beta_m], [\tilde\beta_m, \tilde\beta_{m+1}] \},
\end{align*}
it follows that 
\[ \abs{A_\mc{P}(\gamma_v[0,t]) - A_\mc{P}(\hat\gamma_v(t)) } \leq 6.\]

A similar reasoning applies for the other possible shortenings of $\tilde{\alpha}$ and $\tilde{\beta}$.
\end{proof}

\begin{proof}[Proof of Lemma \ref{lem:sequ_roughly_coincide}]
Let $\tilde{p}$ and $\tilde{q}$ (resp.\ $\tilde{p}'$ and $\tilde{q}'$) be the starting and end point of $\tilde{\gamma}$ (resp.\ $\tilde{\gamma}'$). The curves $\pi^{-1}(\P)$ divide $\HH^2$ into convex domains and denote by $\tilde{P}$, $\tilde{Q}$, $\tilde{P}'$, $\tilde{Q}'$ the domains containing $\tilde{p}, \tilde{q}, \tilde{p}', \tilde{q}'$, respectively. 

Let us first assume that $\tilde{P} = \tilde{P}'$ and $\tilde{Q}=\tilde{Q}'$. In this case both sequences coincide $\tilde{\alpha} = \tilde{\beta}$: If $\tilde{\alpha}_1 \neq \tilde{\beta}_1$ then the endpoints $\tilde{q}$ and $\tilde{q}'$ would be contained in disjoint half-spaces $H_1$ and $H'_1$ bordering on $\tilde{\alpha}_1$ and $\tilde{\beta}_1$, respectively. Hence, $d(\tilde{q},\tilde{q}') \geq d(H,H') \geq o > \rho$; a contradiction. Therefore, one may shorten both $\tilde{\gamma}$ and $\tilde{\gamma}'$ at the beginning such that the resulting arcs both start on $\tilde{\alpha}_1 = \tilde{\beta}_1$. Iterating this procedure $r$-times until there are no more intersections with $\pi^{-1}(\P)$ for $\tilde{\gamma}$ or $\tilde{\gamma}'$, yields that $\tilde{\alpha}_i=\tilde{\beta}_i$ for all $i=1,\ldots,r$. If $r$ was not equal to both $l$ and $m$, one of the two arcs would intersect another curve in $\pi^{-1}(\P)$ whence it would leave the component $\tilde{Q}=\tilde{Q}'$ and the two endpoints $q$ and $q'$ would not be in the same connected component; a contradiction. 

Now, suppose that $\tilde{P} \neq \tilde{P}'$. It follows from (\ref{eq:dist_ends_vec}) that $d(\tilde{p},\tilde{p}') < \rho$. Since $\rho$ is less than the minimal distance $o$ between any two pants curves in $\P$, $\tilde{P}$ and $\tilde{P}'$ are adjacent convex domains separated by some $\tilde{\eta} \in \pi^{-1}(\P)$. We claim that either $\tilde{\alpha}_1 = \tilde{\eta}$ or $\tilde{\beta}_1 = \tilde{\eta}$ holds: Note that $\gamma'$ is filling, so that it intersects every pants curve at least once and, in particular, $m \geq 3$. 

Therefore, if both $\tilde{\alpha}_1 = \tilde{\eta}$ and $\tilde{\beta}_1 = \tilde{\eta}$ holds, then the endpoint $\tilde{q}$ is contained in the half-space $H$ bordering on $\eta$ that does not contain $\tilde{p}$ and $\tilde{q}'$ is contained in the half-space $H'$ bordering on $\tilde{\alpha}_2$ which is disjoint from $H$. Thus, $d(\tilde{q},\tilde{q}') \geq o > \rho$; a contradiction.

Likewise, if both $\tilde{\alpha}_1 \neq \tilde{\eta}$ and $\tilde{\beta}_1 \neq \tilde{\eta}$ then the endpoints $\tilde{q}$ and $\tilde{q}'$ are contained in disjoint half-spaces $H$ and $H'$ bordering on $\tilde{\alpha}_1$ and $\tilde{\beta}_1$, respectively. Thus, $d(\tilde{o},\tilde{o}') \geq o > \rho$; a contradiction. 

Because either $\tilde{\alpha}_1 = \tilde{\eta}$ or $\tilde{\beta}_1 = \tilde{\eta}$, we may either shorten $\tilde{\gamma}$ or $\tilde{\gamma}'$ and delete $\tilde{\alpha}_1$ or $\tilde{\beta}_1$ from $\tilde{\alpha}$ or $\tilde{\beta}$, respectively, to obtain two arcs starting in the same convex domain $\tilde{P} = \tilde{P}'$. Similarly, we may shorten either $\tilde{\gamma}$ or $\tilde{\gamma}'$ at the end to obtain arcs that end in the same domain $\tilde{Q}=\tilde{Q}'$. This reduces the situation to the previous case. This proves that $\tilde{\alpha} = \tilde{\beta}$ up to deletion at the beginning or the end.
\end{proof}

\section{Lower Bound on the number of arcs} \label{sec:lower_bound}

In this section we derive our main dynamic estimate which yields a lower bound on the number or arcs appearing in a random geodesic.

\subsection{Notation}

We denote by $\mu$ the normalized Liouville (probability) measure on $\UT(X)$, i.e., $\mathcal{L}_X/\mathcal{L}_X(\UT(X))$.
Furthermore, we denote by $T \colon \UT(X) \to \UT(X)$ the time-one map of the geodesic flow (see Subsection~\ref{subsec:dynamics}).
Also, we will use the notation
\[\mu(f) = \int_{\UT(X)} f \, d\mu\]
for $f \in L^1(\UT(X),\mu)$. 
By $\Sob(\phi)$ we will denote the Sobolev norm of a smooth function $\phi \in C_c^\infty(\UT(X))$ as defined in section \ref{sec:sobolev_norms}.

\subsection{Smooth approximations of orthogeodesic sets} \label{subsubsect:orthogeodsets}

Recall that $X$ is a hyperbolic surface obtained from gluing $\abs{\chi(X)}$-many isometric pairs of pants along their cuffs without twists. 
In each pair of pants there are infinitely many geodesic arcs that connect any two boundary components. 
We consider free homotopy classes of these arcs where the endpoints are allowed to glide freely along their respective boundary component.
In each of these free homotopy classes there is a unique orthogeodesic arc that meets both boundaries perpendicularly. 
For such an orthogeodesic arc $o \subset X$ we denote its length by $\ell(o)$.
Now, we enumerate all orthogeodesic arcs in all pairs of pants $\{o_i\}_{i \in \NN}$  by increasing length, i.e.\ $\ell(o_i) \leq \ell(o_j)$ if $i \leq j$.
The length of the $i$-th orthogeodesic arc $o_i$ will be denoted by $\ell_i = \ell(o_i)$, and for every orthogeodesic arc $o_i$ we denote by $P_i$ the pair of pants that contains it.
Moreover, we consider the set $B_i$ of all unit tangent vectors $v \in \UT(P_i)$ whose associated geodesic $\gamma_v \cap P_i$ is in the free homotopy class of $o_i$. In this way the sets $\{B_i\}_{i \in \NN}$ are a certain enumeration of the sets $\{U_k^l\}_{k,l}$ in section \ref{subsec:geometry}.

Given such an orthogeodesic arc $o_i$ we may lift it to an orthogeodesic arc $\tilde{o}_i$ in the universal covering $\HH^2$, such that $\tilde{o}_i$ is a subarc of the imaginary axis in the upper half-plane model. 
Furthermore, we may assume that the lifts of its boundary components are half-circles $C_0(\ell_i)$ and $C_1(\ell_i)$ centered at $0$ of radius $e^{-\ell_i /2}$ and $e^{\ell_i/2}$, respectively. 
Then the universal covering $\pi \colon \HH^2 \to X$ induces a bijection between the set $B_i$ and the set $B(\ell_i)$ of all unit tangent vectors $v \in \UT(\HH^2)$, whose base point lies between $C_0(\ell_i)$ and $C_1(\ell_i)$ and whose induced geodesic $\gamma_v \subset \HH^2$ intersect $C_0(\ell_i)$ and $C_1(\ell_i)$.

    \begin{figure}[h!]
\centering{
\resizebox{160mm}{!}{\fontsize{12pt}{12pt}\selectfont
\begingroup%
  \makeatletter%
  \providecommand\color[2][]{%
    \errmessage{(Inkscape) Color is used for the text in Inkscape, but the package 'color.sty' is not loaded}%
    \renewcommand\color[2][]{}%
  }%
  \providecommand\transparent[1]{%
    \errmessage{(Inkscape) Transparency is used (non-zero) for the text in Inkscape, but the package 'transparent.sty' is not loaded}%
    \renewcommand\transparent[1]{}%
  }%
  \providecommand\rotatebox[2]{#2}%
  \newcommand*\fsize{\dimexpr\f@size pt\relax}%
  \newcommand*\lineheight[1]{\fontsize{\fsize}{#1\fsize}\selectfont}%
  \ifx\svgwidth\undefined%
    \setlength{\unitlength}{327.76167982bp}%
    \ifx\svgscale\undefined%
      \relax%
    \else%
      \setlength{\unitlength}{\unitlength * \real{\svgscale}}%
    \fi%
  \else%
    \setlength{\unitlength}{\svgwidth}%
  \fi%
  \global\let\svgwidth\undefined%
  \global\let\svgscale\undefined%
  \makeatother%
  \begin{picture}(1,0.4884478)%
    \lineheight{1}%
    \setlength\tabcolsep{0pt}%
    \put(0,0){\includegraphics[width=\unitlength,page=1]{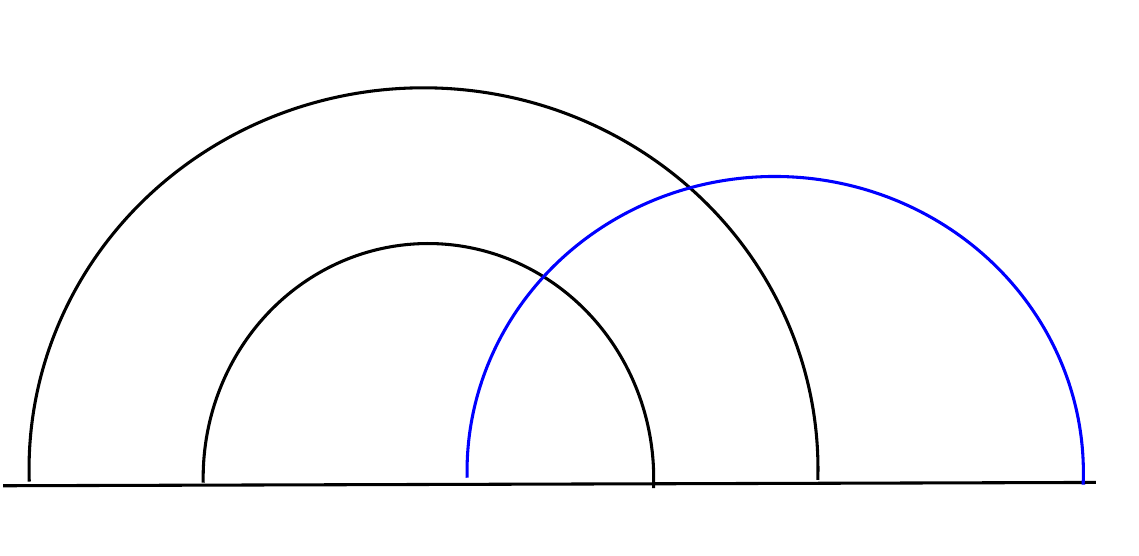}}%
    \put(0.14911014,0.01401247){\color[rgb]{0,0,0}\makebox(0,0)[lt]{\lineheight{40}\smash{\begin{tabular}[t]{l}$e^{-l_i/2}$\end{tabular}}}}%
    \put(0.69636747,0.01240593){\color[rgb]{0,0,0}\makebox(0,0)[lt]{\lineheight{40}\smash{\begin{tabular}[t]{l}$e^{l_i/2}$\end{tabular}}}}%
    \put(0,0){\includegraphics[width=\unitlength,page=2]{bi.pdf}}%
    \put(0.37520373,0.33168116){\color[rgb]{0,0,0}\makebox(0,0)[lt]{\lineheight{40}\smash{\begin{tabular}[t]{l}$\tilde{o_i}$\end{tabular}}}}%
    \put(0.51141458,0.24501754){\color[rgb]{0,0,0}\makebox(0,0)[lt]{\lineheight{40}\smash{\begin{tabular}[t]{l}$v$\end{tabular}}}}%
    \put(0.70720325,0.35124744){\color[rgb]{0,0,0}\makebox(0,0)[lt]{\lineheight{40}\smash{\begin{tabular}[t]{l}$\widetilde{\gamma_v}$\end{tabular}}}}%
    \put(-0.00401379,0.0104954){\color[rgb]{0,0,0}\makebox(0,0)[lt]{\lineheight{40}\smash{\begin{tabular}[t]{l}$C_1(\ell_i)$\end{tabular}}}}%
    \put(0.53978939,0.01464938){\color[rgb]{0,0,0}\makebox(0,0)[lt]{\lineheight{40}\smash{\begin{tabular}[t]{l}$C_0(\ell_i)$\end{tabular}}}}%
  \end{picture}%
\endgroup%
}
    \caption{This figure shows the set $B(\ell_i)$ of all unit tangent vectors $v \in \UT(\HH^2)$, whose base point lies between $C_0(\ell_i)$ and $C_1(\ell_i)$ and whose induced geodesic $\gamma_v \subset \HH^2$ intersect $C_0(\ell_i)$ and $C_1(\ell_i)$. }\label{orthopic}
}
\end{figure}

There is a formula for the (normalized) volume of $B(\ell)$:

\blem \label{lem:volume_B_ell}
    \begin{enumerate}
        \item \label{lem:volume_B_ell:formula} The measure of $B(\ell)$ is given by
        \[ \mu(B(\ell)) = \frac{1}{2 \pi \abs{\chi(X)}} \cdot \kl{ Li_2\left(\frac{1}{\cosh ^2\left(\frac{\ell}2\right)}\right)+\f 1 2\log\left(\frac{1}{\cosh ^2\left(\frac{\ell}2\right)}\right)\log\abs{1-\frac{1}{\cosh ^2\left(\frac{\ell}2\right)}} }.\]
        \item \label{lem:volume_B_ell:asymptotic} Asymptotically, this yields 
        \[ \mu(B(\ell)) \sim \frac{1}{\pi \abs{\chi(X)} } \cdot \ell e^{- \ell}\]
        as $\ell \to \infty$. 
    \end{enumerate}    
\elem
\bpf
\begin{enumerate}
    \item This follows from Bridgeman's computations in \cite[Section~9]{B11:Dilogarithm}.
    
    \item Observe that 
    \begin{align*}
        \mu(B(\ell)) 
        &= \frac{1}{2 \pi \abs{\chi(X)}} \cdot \kl{ Li_2\left(\frac{1}{\cosh ^2\left(\frac{\ell}2\right)}\right)+\f 1 2\log\left(\frac{1}{\cosh ^2\left(\frac{\ell}2\right)}\right)\log\abs{1-\frac{1}{\cosh ^2\left(\frac{\ell}2\right)}} } \\
        &= \frac{1}{2 \pi \abs{\chi(X)}} \cdot \kl{ Li_2\left(\frac{1}{\cosh ^2\left(\frac{\ell}2\right)}\right)
        -\log\left(\cosh \left(\frac{\ell}2\right)\right)\log \kl{ \frac{\sinh^2\kl{\frac{\ell}2}}{\cosh^2 \kl{\frac{\ell}2}} } } \\
        &= \frac{1}{2 \pi \abs{\chi(X)}} \cdot \kl{ Li_2\left({\cosh^{-2}\left(\frac{\ell}2\right)}\right)
        +2\log\left(\cosh \left(\frac{\ell}2\right)\right)\log \kl{ \frac{\cosh\kl{\frac{\ell}2}}{\sinh\kl{\frac{\ell}2}} } }
    \end{align*}
    
    Recall that
    \[
        \log(\cosh(x)) \sim x
    \]
    as $x \to \infty$. Moreover,
    \[
        \log \kl{ \frac{\cosh\kl{x}}{\sinh\kl{x}} } = \log \kl{ 1 + e^{-2x} } - \log \kl{ 1 - e^{-2x} } \sim 2 e^{-2x}
    \]
    as $x \to \infty$.
    Thus, with $x = \frac{\ell}{2}$, we have that
    \begin{align} \label{eq:asympt1}
        2\log\left(\cosh \left(\frac{\ell}2\right)\right)\log \kl{ \frac{\cosh\kl{\frac{\ell}2}}{\sinh\kl{\frac{\ell}2}} } \sim 2 \frac{\ell}{2} \cdot  2 e^{-\ell} = 2 \ell e^{-\ell}
    \end{align}
    as $\ell \to \infty$.
    
    On the other hand 
    \[ Li_2(z) = \sum_{k=1}^\infty \frac{z^k}{k^2} = z + O(z^2), \]
    and
    \[ \cosh^{-2} \kl{\frac{\ell}{2}} \sim 4 e^{-\ell}. \]
    Thus,
    \begin{align} \label{eq:asympt2}
    Li_2\left({\cosh^{-2}\left(\frac{\ell}2\right)}\right) \sim 4 e^{-\ell}.
    \end{align}
    
    Using (\ref{eq:asympt1}) and (\ref{eq:asympt2}) we get that
    \[ \mu(B(\ell)) \sim \frac{1}{2 \pi \abs{\chi(X)}} \cdot \kl{ 2 \ell e^{-\ell} + 4 e^{-\ell} } \sim \frac{1}{\pi \abs{\chi(X)}} \cdot \ell e^{-\ell}\]
    as $\ell \to \infty$.
\end{enumerate}
\epf

\brem
    It follows from Lemma \ref{lem:volume_B_ell} \ref{lem:volume_B_ell:formula} that 
   $\ell_i\leq \ell_{i+1}$ implies $\mu(B(\ell_i))\geq \mu(B(\ell_{i+1}))$. In particular, the sequence $\{ \mu(B(\ell_i))\}_{i \in \NN}$ is decreasing.
\erem

Using Corollary \ref{cor:counting_arcs} we obtain the following asymptotic for the length of the $i$-th orthogeodesic arc.

\blem \label{lem:asymptoticLength}
    Let $0 < \delta <1$ denote the Hausdorff dimension of the limit set for (any-) one of the constituent pairs of pants of $X$. Then
    \[ \ell_i \sim \frac{1}{\delta} \cdot \log(i).\]
    
    Consequently, there are sequences $(\theta_i)_{i \in \NN}, (\theta'_i)_{i \in \NN}$ converging to $1$ as $i \to \infty$, such that
    \be
        \mu(B(\ell_i)) = \frac{1}{\pi \abs{\chi(X)} \delta} \log(i) i^{-\theta_i/\delta} \theta'_i.  
    \ee
\elem
\bpf
    Let us fix two boundary curves $C_-, C_+$ of one of the constituent pairs of pants $P$ and denote by $N_{C_-,C_+}(\ell)$ the number of orthogeodesic arcs in $P$ of length $\leq \ell$ connecting the two given boundary curves.
    By Corollary \ref{cor:counting_arcs} there is a constant $C_0>0$ such that
    \begin{align} \label{eq:count_one_pants}
        N_{C_-,C_+}(\ell) \sim C_0 \cdot e^{\delta \ell}
    \end{align}
    as $\ell \to \infty$, where $\delta > 0$ is the Hausdorff dimension of the limit set of $P$. 
    
    Let us now denote by $N(\ell)$ the \emph{total} number of \emph{all} orthogeodesic arcs in any constituent pair of pants of length $\leq \ell$. It follows from (\ref{eq:count_one_pants}) that there is a constant $C>0$ such that
    \[N(\ell) \sim C \cdot e^{\delta \ell}.\]
    In particular, there are constants $0< C_1 < C_2$ such that
    \[ C_1 \cdot e^{\delta \ell} \leq N(\ell) \leq C_2 \cdot e^{\delta \ell}\]
    for all $\ell >0$.
    
    Because we have enumerated the orthogeodesics $\{ o_i \}_{i \in \NN}$ by increasing length $\{\ell_i\}_{i \in \NN}$, there is a constant $D>0$ such that
    \[ N(\ell_i) - D \leq i \leq N(\ell_i) + D\]
    for all $i \in \NN$. Thus,
    \[ i \leq N(\ell_i) + D \leq C_2 e^{\delta \ell_i} + D, \]
    equivalently
    \[ \frac{1}{\delta} \log(i - D) - \frac{1}{\delta} \log(C_2) \leq \ell_i, \]
    and similarly
    \[ \ell_i \leq \frac{1}{\delta} \log(i + D) - \frac{1}{\delta} \log(C_1).\]
    Therefore,
    \[ \ell_i \sim \frac{1}{\delta} \log(i) \]
    as $i \to \infty$. 
\epf
Moreover, we can approximate each of the sets $B(\ell)$ by a smooth bump function $\phi_\ell \in C^\infty_c(\UT(\HH^2))$. In the rest of this section we will explicitly construct for every $\ell>0$ a smooth function $\phi_\ell \in C^\infty_c(\UT(\HH^2))$, that satisfies $0 \leq \phi_\ell \leq 1$ and $\supp \phi_\ell \subseteq B(\ell)$.

Let $\epsilon > 0$. 
For any subset $A \subseteq \UT(X)$ we denote by $N_{\epsilon}(A)$ the $\epsilon$-neighborhood of $A$, i.e.\
\[ N_{\epsilon}(A) \coloneqq \{ x \in \UT(X) \, | \, d(x,A) < \epsilon \}, \]
where $d(x,A) \coloneqq \inf \{ d(x,y) \, | \, y \in A \}$ denotes the distance from $x$ to the set $A$. We can define the \emph{$\epsilon$-interior of $A$} as 
\[ \interior_\epsilon(A) := A \setminus N_\epsilon(A^c),\]
i.e.\ all points in $A$ that have distance more than $\epsilon$ from its boundary.

Let $\{\psi_{\epsilon}\}_{0<\epsilon<\epsilon_0} \subseteq C_c^\infty(G)$ be the family of smooth approximate convolutional identities as constructed in Appendix \ref{sec:approx_id}. In particular, they satisfy
\begin{enumerate}
    \item $\psi_\epsilon \geq 0$, and
    \item $\int_G \psi_\epsilon \, d\nu = 1$, and
    \item $\supp(\psi_\epsilon) \subseteq B_{\epsilon}(e)$,
\end{enumerate}
for all $\epsilon > 0$. In addition, for every $d>0$ there is a constant $K = K(d) > 0$ such that 
\[\norm{ E_\alpha \cdot \psi_\epsilon}_1 \leq K \cdot \kl{ \frac{1}{\epsilon}}^d \]
for all $\abs{\alpha} \leq d$, where $E_\alpha$ denotes the multi-index derivative by left-invariant vector fields (see section \ref{sec:sobolev_norms} and Lemma \ref{lem:approx_id_L1}).

For every $\ell > 0$ and $0< \epsilon = \epsilon(\ell) < \epsilon_0$, that may depend on $\ell$, we set
\[ \phi_\ell \coloneqq \psi_{\epsilon/2} * \II_{\interior_{\epsilon/2}(B(\ell))} \in C_c^\infty(\UT(X)).\]
These smooth approximations satisfy the following lemma.

\blem\label{lem:smoothapprox}
    Let $\epsilon = \epsilon(\ell) > 0$. Then there are constants $C_1, C_2, C_3 >0$, such that
    \begin{enumerate}
        \item $0 \leq \phi_\ell \leq 1;$
        \item $\supp \phi_\ell \subseteq B(\ell);$
        \item $\| \II_{B(\ell)} - \phi_\ell \|_1 \leq C_1 \epsilon(\ell);$ \label{lem:smoothapprox:L1estimate}
        \item $ \Sob(\phi_\ell) \leq C_2 \cdot \mu(B(\ell)) \cdot \kl{ \frac{1}{\epsilon(\ell)}}^d$, where $d$ denotes the degree of the Sobolev norm;
        \item $\Sob(\bar{\phi}_\ell) \leq C_3 \cdot \kl{ \frac{1}{\epsilon(\ell)}}^d$, where we denote $\bar{\phi}_\ell \coloneqq 1 - \phi_\ell$.
    \end{enumerate}
\elem

\begin{proof}
    \begin{enumerate}
        \item Follows immediately from the construction.
        
        \item By definition $\supp (\II_{\interior_{\epsilon/2}(B(\ell))}) = \closure{\interior_{\epsilon/2}(B(\ell))}$. Thus every point in $\supp( \II_{\interior_{\epsilon/2}(B(\ell))})$ has distance at least $\epsilon/2$ from $\partial B(\ell)$. Because $\supp (\psi_{\epsilon/2}) \subseteq B_{\epsilon/2}(\id)$ and $d_{\UT(X)}(g v,w) \leq d_G(g,\id)$ for all $v,w \in \UT(X)$, $g \in G$, it follows by the definition of convolution that $\supp (\phi_\ell) \subseteq B(\ell)$.
        
        \item By construction $\phi_\ell$ differs from $\II_{B(\ell)}$ only on the \emph{inner $\epsilon$-tube neighborhood of $B(\ell)$}
        \[ T_\epsilon \coloneqq B(\ell) \setminus \interior_{\epsilon}(B(\ell)) \subseteq B(\ell).\]
        Because both functions take values in $[0,1]$, one obtains the trivial bound
        \[ \| \II_{B(\ell)} - \phi_\ell \|_1 \leq \mu(T_\epsilon).\]
        
        It can be shown that both the principal curvature of the smooth boundary components of $\partial B(\ell)$ and its area $\Area(\partial B(\ell))$ are uniformly bounded. Thus, there exists $C_1>0$ independent of $\ell$ and $0<\epsilon<\epsilon_0$ such that 
        \[ \mu(T_\epsilon) \leq C_1 \epsilon(\ell)\]
        by the Tube Lemma \ref{lem:tube_lemma}.
   
        It is well-known that geodesic trajectories in $\UT(\HH^2)$ are given unit vector fields of constant slope along geodesic lines in $\HH^2$. Note that the unit vectors on $B(\ell)$ have basepoint in a given ideal quadrilateral, such that forward and backwards trajectories intersect two designated opposite sides of the ideal quadrilateral, which we call \textit{opposing sides}.
        
       \begin{figure}[h!]
\centering{
\resizebox{160mm}{!}{\fontsize{12pt}{12pt}\selectfont
\begingroup%
  \makeatletter%
  \providecommand\color[2][]{%
    \errmessage{(Inkscape) Color is used for the text in Inkscape, but the package 'color.sty' is not loaded}%
    \renewcommand\color[2][]{}%
  }%
  \providecommand\transparent[1]{%
    \errmessage{(Inkscape) Transparency is used (non-zero) for the text in Inkscape, but the package 'transparent.sty' is not loaded}%
    \renewcommand\transparent[1]{}%
  }%
  \providecommand\rotatebox[2]{#2}%
  \newcommand*\fsize{\dimexpr\f@size pt\relax}%
  \newcommand*\lineheight[1]{\fontsize{\fsize}{#1\fsize}\selectfont}%
  \ifx\svgwidth\undefined%
    \setlength{\unitlength}{327.76167982bp}%
    \ifx\svgscale\undefined%
      \relax%
    \else%
      \setlength{\unitlength}{\unitlength * \real{\svgscale}}%
    \fi%
  \else%
    \setlength{\unitlength}{\svgwidth}%
  \fi%
  \global\let\svgwidth\undefined%
  \global\let\svgscale\undefined%
  \makeatother%
  \begin{picture}(1,0.4884478)%
    \lineheight{1}%
    \setlength\tabcolsep{0pt}%
    \put(0,0){\includegraphics[width=\unitlength,page=1]{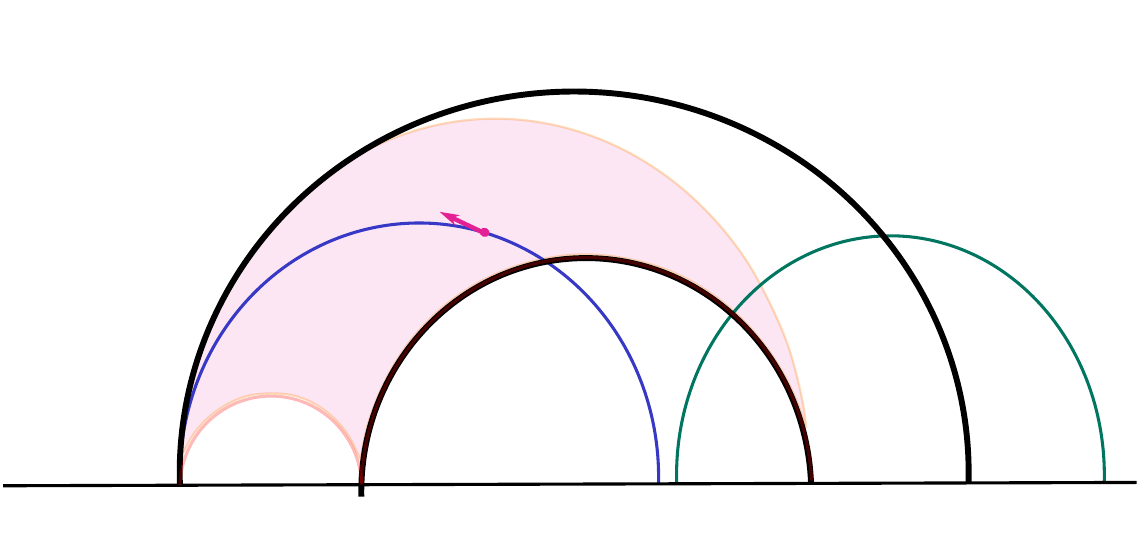}}%
    \put(0.46649295,0.27896873){\color[rgb]{0,0,0}\makebox(0,0)[lt]{\lineheight{40}\smash{\begin{tabular}[t]{l}$v$\end{tabular}}}}%
    \put(0,0){\includegraphics[width=\unitlength,page=2]{faces.pdf}}%
    \put(0.67041196,0.19948411){\color[rgb]{0,0,0}\makebox(0,0)[lt]{\lineheight{40}\smash{\begin{tabular}[t]{l}$w$\end{tabular}}}}%
  \end{picture}%
\endgroup%
}
    \caption{This figure shows a vector $v$ is in a side face of $B(\ell)$. The basepoint of $v$ is located within an ideal triangle determined by 3 vertices of $B(\ell)$. The vector $w$ is in an opposite face of $B(\ell)$, with basepoint one of the two geodesics delimiting $B(\ell)$ and geodesic crossing the other delimiting geodesic. }\label{fig:faces}
}
\end{figure}
        
        The boundary of $\partial B(\ell)$ can be described as the unions of 6 faces: 4 which we call \textit{side faces} and 2 that we call \textit{opposing faces}, which we will now describe.
        Side faces are composed by unit vectors whose basepoints are in an ideal triangle formed by 3 out of the 4 ideal vertices associated to $B(\ell)$ (see vector $v$ in Figure~\ref{fig:faces} with basepoint in the shaded ideal triangle therein), and their geodesic trajectories end at the ideal vertex that does not pair to form an opposing side in the triangle. Hence, side faces are totally geodesic along the geodesic lines ending at the ideal vertex. In the orthogonal direction we have horocycles tangent at the ideal vertex and the unit vectors are normal to those horocycles. Since $\UT(\HH^2)$ is a Riemannian submersion, these horocyclic trajectories have curvature equal to $1$. Because of the symmetries of side face, these horocyclic trajectories are lines of curvature, which tells us the the principal curvatures of a side face are $0$ (the totally geodesic direction) and $1$ (the horocyclic direction). 
        
        On the other hand, opposing faces are composed of vectors whose basepoints are on one of the opposing sides of the ideal quadrilateral, such that their geodesic trajectory intersects the other opposing side (see vector $w$ in Figure~\ref{fig:faces}). Since unit vectors based at a geodesic are a totally geodesic submanifold (by the geodesic description of $\UT(\HH^2)$), it follows that opposing faces are totally geodesic in $\UT(\HH^2)$, because these faces are open sets of totally geodesic submanifolds.

        \item Recall that
        \[ \Sob(\phi_\ell) = \sum_{\abs{\alpha} = 0}^d \norm{ E_\alpha \cdot \phi_\ell }_{2}  .\]
        
        By Lemma \ref{lem:conv_properties} we get
        \begin{align*}
            \norm{ E_\alpha \cdot \phi_\ell}_2 
            &= \norm{E_\alpha \cdot (\psi_{\epsilon/2} * \II_{\interior_{\epsilon/2}(B(\ell))})}_2 \\
            &= \norm{(E_\alpha \cdot \psi_{\epsilon/2}) * \II_{\interior_{\epsilon/2}(B(\ell))}}_2 \\
            &\leq \norm{E_\alpha \cdot \psi_{\epsilon/2}}_{1} \cdot \norm{ \II_{\interior_{\epsilon/2}(B(\ell))} }_2 \\
            &\leq K \cdot \kl{ \frac{2}{\epsilon(\ell)} }^d \cdot \mu(B(\ell)).
        \end{align*}
        Hence,
        \[ \Sob(\phi_\ell) =  \sum_{\abs{\alpha} = 0}^d \norm{ E_\alpha \cdot \phi_\ell }_{2}  \leq C_2 \cdot \mu(B(\ell)) \cdot \kl{ \frac{1}{\epsilon(\ell)} }^d \]
        for $C_2= d K 2^d > 0$.
        
        \item Observe that
        \[ \bar{\phi}_\ell = 1 - \bar\phi_\ell = 1 - \psi_{\epsilon/2} * \II_{\interior_{\epsilon/2}(B(\ell))} = \psi_{\epsilon/2} * (1 - \II_{\interior_{\epsilon/2}(B(\ell))}). \]
        
        As before we get
        \begin{align*}
            \norm{ E_\alpha \cdot \bar\phi_\ell}_2 
            &\leq \norm{E_\alpha \cdot \psi_{\epsilon/2}}_{1} \cdot \norm{1-  \II_{\interior_{\epsilon/2}(B(\ell))} }_2 \\
            &\leq K \cdot \kl{ \frac{2}{\epsilon(\ell)} }^d.
        \end{align*}
        Hence,
        \[ \Sob(\bar\phi_\ell) =  \sum_{\abs{\alpha} = 0}^d \norm{ E_\alpha \cdot \bar\phi_\ell }_{2}   \leq C_3  \cdot \kl{ \frac{1}{\epsilon(\ell)} }^d \]
        for some $C_3 > 0$.
    \end{enumerate}
    
\end{proof}

\subsection{Counting the number of visited sets} \label{sec:counting}
The orthogeodesic sets $\{B_i\}_{i \in \NN}$ defined in section \ref{subsubsect:orthogeodsets} partition the unit tangent bundle $\UT(X)$ (up to a null set). For a given unit tangent vector $v \in \UT(X)$ let $C_N(v)$ denote the number of \emph{different} sets in this partition that the orbit $v, Tv, \ldots, T^{N-1} v$ has visited, i.e.
\[ C_N(v) \coloneqq \# \{ i \in \NN \, | \, \exists 0 \leq n \leq N-1: T^n v \in B_i \}. \]

It follows immediately from the definition that
\[ A_\P(\gamma_v[0,N]) \geq C_N(v);\]
see also Remark \ref{rem:visiting_counting}.
Therefore, we would like to find a certain growth rate $F(N)$ for the number of visited sets such that for $\mu$-almost every $v \in \UT(X)$ we have
\[ C_N(v) \geq F(N)\]
for $N$ large enough.
This is established by the following theorem.

\bthm\label{mainestimate}
    Let $\eta > 1$ and let $F \colon \NN \to (0,+\infty)$ be a function satisfying 
    \[F(N) \in O(N^{\delta/2\eta}).\]
    
    Then, for $\mu$-almost every $v \in \UT(X)$ there is an $N_0=N_0(v) \in \NN$ such that
    \[C_N(v) \geq F(N)\]
    for all $N \geq N_0$.
\ethm

\begin{proof}
    We need to prove that the set
    \[ \Omega \coloneqq \set{ v \in \UT(X) \, \bigg| \, \exists N_0 \in \NN \, \forall N \geq N_0:C_N(v) \geq F(N)} \]
    has full measure.
    
    Let $N \in \NN$ and denote by $\Omega_N \subseteq \UT(X)$ the set of all vectors $v \in \UT(X)$ such that 
    \[C_N(v) \geq F(N).\]
    With this notation we have
    \[ \Omega = \liminf_{N \to \infty} \, \Omega_N = \bigcup_{N_0 \in \NN} \bigcap_{N \geq N_0} \Omega_N.\]
    In order to show that $\Omega$ has full measure we need to see that its complement
    \[ \Omega^c = \limsup_{N \to \infty} \, \Omega^c_N = \bigcap_{N_0 \in \NN} \bigcup_{N\geq N_0} \Omega^c_N\]
    has measure zero. By the Borel--Cantelli Lemma this will follow if the sequence $(\mu(\Omega^c_N))_{N \in \NN}$ is summable, i.e.\ 
    \begin{align} \label{eq:summable}
        \sum_{N=1}^\infty \mu(\Omega^c_N) < + \infty.
    \end{align}
    
    Before we proceed with verifying (\ref{eq:summable}) we make some general observations first: if the orbit $v, Tv, \ldots, T^{N-1} v$ has visited less than $m$ sets, i.e.\ $C_N(v) < m$, then $v$ has not visited all of the first $m$ sets $B_1, \ldots, B_m$ in particular. Thus,
    \[ \mu(C_N < m) \leq \mu \kl{ \bigcup_{i=1}^m \bigcap_{k=0}^{N-1} T^{-k}(B^c_i)} \leq  \sum_{i=1}^m \mu \kl{ \bigcap_{k=0}^{N-1} T^{-k}(B^c_i)}. \]

Notice that
\[\mu \kl{ \bigcap_{k=0}^{N-1} T^{-k}(B^c_i)} 
        = \int_{\UT(X)} \prod_{k=0}^{N-1} \kl{\II_{B^c_i} \circ T^k} \, d\mu.\]
We want to apply our exponential mixing theorem \ref{thm:expmixing} to the right-hand-side.
For this we will approximate the sets $B_i$ by smooth functions $\phi_i \in C_c^\infty(\UT(X))$ as described in section \ref{subsubsect:orthogeodsets}.
For each $i \in \NN$ we let $\epsilon_i>0$ be a positive number that we will specify later.
Given $\epsilon(\ell_i) = \epsilon_i$ we obtain smooth functions $\phi_{\ell_i} \in C^\infty_c(\UT(\HH^2))$ that approximate the sets $B(\ell_i) \subseteq \UT(\HH^2)$. 
Via the isometry $B_i \cong B(\ell_i)$ these pull-back to smooth approximations of the set $B_i$ that satisfy the assertions of Lemma \ref{lem:smoothapprox}.
Finally, we set $\bar\phi_i \coloneqq 1 - \phi_i \geq \II_{B^c_i}$.

Thus, we get
\begin{align*}
\int_{\UT(X)} \prod_{k=0}^{N-1} \kl{\II_{B^c_i} \circ T^k} \, d\mu 
    & \leq \int_{\UT(X)} \prod_{k=0}^{N-1} \kl{\bar\phi_i \circ T^k} \, d\mu
\end{align*}

In order to obtain an error term that tends to zero from Theorem \ref{thm:expmixing}, we will take the product only over an equally spaced subset of $\{0, \ldots, N-1\}$. Let $n \in \NN$, which we will also specify later, and set 
\[l \coloneqq \floor{ \frac{N-1}{n-1}  }.\] 
By Theorem \ref{thm:expmixing} we obtain
\begin{align*} 
\int_{\UT(X)} \prod_{k=0}^{N-1} \kl{\bar\phi_i \circ T^k} \, d\mu 
    &\leq \int_{\UT(X)} \prod_{k=0}^{n-1} \kl{\bar\phi_i \circ T^{k \cdot l}} \, d\mu \\
    &\leq \mu(\bar\phi_i)^n + L \cdot (n-1) \cdot D^n \cdot \Sob(\bar\phi_i)^n \cdot e^{-\beta \cdot l}
\end{align*}
for all $l\geq K_0$, where $\Sob$ denotes the degree $d=3$ Sobolev norm and the constants $\alpha, \beta>0$ are as in Theorem \ref{thm:expmixing}.

Therefore,
\begin{align*}
    \mu(C_N < m) 
        &\leq \sum_{i=1}^m \mu \kl{ \bigcap_{k=0}^{N-1} T^{-k}(B^c_i)} \\
        &\leq \sum_{i=1}^m \kl{ \mu(\bar\phi_i)^n + L \cdot (n-1)\cdot D^n \cdot \Sob(\bar\phi_i)^n \cdot e^{-\beta \cdot l} }\\
        &\leq \underbrace{\sum_{i=1}^m \mu(\bar\phi_i)^n}_{ \eqqcolon A_1(N) } + \underbrace{L \cdot (n-1) \cdot D^n \cdot  e^{-\beta \cdot l} \cdot  \sum_{i=1}^m \Sob(\bar\phi_i)^n}_{\eqqcolon A_2(N)}. \\
\end{align*}

Thus, we are left with proving that both sequences $\{A_1(N)\}_{N\in\NN}$ and $\{A_2(N)\}_{N \in \NN}$ are summable for $m = m(N) \coloneqq \ceiling{F(N)} = O( N^{\delta/2 \eta})$ and suitable choices of $\epsilon_i$ and $n=n(N)$. To this end, we set
\begin{align}
    \epsilon_i &\coloneqq \frac{\mu(B_i)}{ C_1 (2 + \log(i)) }, &\forall i \in \NN, \label{eq:defepsiloni}
\end{align}
and choose $n(N) \in \NN$ such that 
\begin{align}
    n(N) \sim N^{1/2\eta'}
\end{align}
for some fixed $1< \eta' < \eta$.

Moreover, note that we may assume without loss of generality that 
\begin{align}
m(N) \asymp N^{\delta/2\eta}, \label{eq:asymptotic_m}
\end{align}
i.e.\ there are positive constants $A,B >0$ and $N_0 \in \NN$ such that $A \cdot N^{\delta/2\eta} \leq m(N) \leq B \cdot N^{\delta/2\eta}$ for all $N \geq N_0$. Indeed, for every function $m_1(N) \in O(N^{\delta/2 \eta})$ there is a function $m_2(N) \asymp N^{\delta/2 \eta}$ such that $m_1(N) \leq m_2(N)$ for sufficiently large $N$. In this case, 
\[\mu(C_N < m_1(N)) \leq \mu(C_N < m_2(N)),\] 
whence it suffices to prove that $(\mu(C_N < m(N)))_{N\in\NN}$ is summable whenever $m(N) \asymp N^{\delta/2\eta}$. In particular, we may assume that $m(N) \to + \infty$ as $N \to +\infty$.

\textbf{Regarding $A_1(N)$.} Observe that
\begin{align}
    \mu(\bar{\phi}_i) \leq 1 - \kl{\mu(B_i) - C_1 \epsilon_i}
    \leq 1 - \tfrac{1}{2} \mu(B_i) \label{eq:barphi_estimate}
\end{align}
for every $i \in \NN$ by (\ref{eq:defepsiloni}).
Because the sets $\{ B_i \}_{i \in \NN}$ are enumerated in such a way that the sequence $(\mu(B_i))_{i\in\NN}$ is monotonically decreasing, we obtain that
\[ A_1(N) \leq \sum_{i=1}^m \mu(\bar{\phi}_i)^n \leq \sum_{i=1}^m \kl{1 - \tfrac{1}{2} \mu(B_i)}^n \leq m \kl{1 - \tfrac{1}{2} \mu(B_m)}^n. \]

By Lemma \ref{lem:asymptoticLength} there are sequences $(\theta_i)_{i \in \NN}, (\theta'_i)_{i \in \NN}$ converging to $1$ as $i \to \infty$, such that
\begin{align}
    \mu(B_i) = \mu(B(\ell_i)) = \frac{1}{\pi \abs{\chi(X)} \delta} \log(i) i^{-\theta_i/\delta} \theta'_i 
    \label{eq:measure_Bm}
\end{align}
for all $i \in \NN$.

Regarding the logarithm we get
\begin{align*}
    \log(A_1(N)) &\leq \log(m) + n \log \kl{1 - \tfrac{1}{2} \mu(B_m) } \\
    & \stackrel{(\diamond)}{\leq} \log(m) - \frac{n}{2} \mu(B_m) \\
    & = \log(m) - \frac{\theta'_m}{2 \pi \abs{\chi(X)} \delta} n \log(m) m^{-\theta_m/\delta} \\
    &= - \log(m) \frac{n}{m^{\theta_m/\delta}} \kl{ \frac{\theta'_m}{2\pi \abs{\chi(X)} \delta } - \frac{m^{\theta_m/\delta}}{n}}{\qquad (\star)}
\end{align*}
where we used the fact that $\log(x+1) \leq x$ for all $x \in \RR$ in $(\diamond)$ and plugged-in (\ref{eq:measure_Bm}). 

From $m=m(N) \asymp N^{\delta/2\eta}$ it follows that
\begin{align}
    \log(m(N)) \sim \frac{\delta}{2\eta} \log(N). \label{eq:log(m(N))}
\end{align}
Furthermore, $n = n(N) \sim N^{1/2\eta'}$, $\eta'<\eta$, $m = m(N) \to + \infty$ as $N \to + \infty$ and $\theta_m \to 1$ as $m \to + \infty$, such that $-\f1{2\eta'}+\f{\theta_m}{2\eta}<0$ for large $m$ and hence
\begin{align}
    \frac{m^{\theta_m/\delta}}{n} \ll \frac{N^{\theta_m/2\eta}}{N^{1/2\eta'}} = N^{-\tfrac{1}{2\eta'} + \tfrac{\theta_m}{2\eta} } \to 0
    \label{as:A1:1}
\end{align}
as $N \to + \infty$. 

Using both (\ref{eq:log(m(N))}) and (\ref{as:A1:1}) in $(\star)$ we obtain that there is a positive constant $C'>0$ such that 
\[ \log(A_1(N)) \leq - C' \frac{n}{m^{\theta_m/\delta}} \log(N)\]
for sufficiently large $N$. Again by (\ref{as:A1:1}) we have that, as $N \to +\infty$, $\frac{n}{m^{\theta_m/\delta}} \to + \infty$ whence
\[ \log(A_1(N)) \leq - 2 \log(N) \]
for sufficiently large $N$. Equivalently,
\[ A_1(N) \leq \frac{1}{N^2}\]
for sufficiently large $N$, which implies that the sequence $(A_1(N))_{N \in \NN}$ is summable.

\textbf{Regarding $A_2(N)$.} Let us consider the logarithm of $A_2(N)$:
\begin{align*}
    \log(A_2(N)) &= \log(L) + \log(n-1) + \log(D) n - \beta l + \log \kl{\sum_{i=1}^m \Sob(\bar{\phi}_i)^n}\\
    &= - \beta l \cdot \kl{1 - \frac{\log(D)}{\beta} \cdot \frac{n}{l} - \frac{\log(n-1)}{\beta l} - \frac{\log(L)}{\beta l} } + \log \kl{\sum_{i=1}^m \Sob(\bar{\phi}_i)^n}
\end{align*}
Note that
\[ l = \floor{\frac{N-1}{n-1}} \sim \frac{N}{n} \sim N^{1- 1/2\eta'} \to + \infty,\]
whence
\[ \frac{n}{l} \sim \frac{N^{1/2\eta'}}{N^{1- 1/2\eta'}} = N^{1/\eta' - 1} \to 0 \]
and
\[ \frac{\log(n-1)}{\beta l} \ll \frac{n}{l} \to 0, \quad \frac{\log(L)}{\beta l} \to 0 \]
as $N \to + \infty$. Therefore,
\begin{align}
    - \beta l \cdot \kl{1 - \frac{\log(D)}{\beta} \cdot \frac{n}{l} - \frac{\log(n-1)}{\beta l} - \frac{\log(L)}{\beta l} } \sim - \beta \cdot N^{1- 1/2\eta'}. \label{as:A2:1}
\end{align}

By Lemma \ref{lem:smoothapprox} there is $C_3 > 0$ such that
\begin{align*}
    \Sob(\bar\phi_i) &= \Sob(\bar{\phi}_{\ell_i}) \leq C_3 \cdot \kl{ \frac{1}{\epsilon_i}}^d
\end{align*}
for all $i \in \NN$. 

Recall that, we set
\[\epsilon_i = \frac{\mu(B_i)}{ C_1 (2 + \log(i)) } \]
and that the sets $B_i$ are enumerated in such a way that their measures $\mu(B_i)$ are decreasing in $i$. Therefore, the sequence $\epsilon_i$ is decreasing in $i$, too. Thus, we obtain
\[\sum_{i=1}^m \Sob(\bar{\phi}_i)^n \leq \sum_{i=1}^m C_3^n \epsilon_i^{-d n} \leq m C_3^n \epsilon_m^{-d n}.\]

From (\ref{eq:measure_Bm}) it follows that $\epsilon_m \asymp m^{-\theta_m/\delta}.$ Thus, there is a positive constant $C_4>0$ such that
\[ m C_3^n \epsilon_m^{-d n} \leq m C_4^n m^{\theta_m n d/\delta}. \]

Regarding the logarithm this yields
\begin{align}
    \log \kl{\sum_{i=1}^m \Sob(\bar{\phi}_i)^n} 
        \leq \log(m) + \log(C_4) n + \theta_m \frac{d}{\delta} n \log(m)
        \ll n \log(m) \ll  N^{1/2\eta'} \cdot \log(N). \label{as:A2:2}
\end{align}

Combining the estimates (\ref{as:A2:1}) and (\ref{as:A2:2}) there is a positive constant $C_5 > 0$ such that
\begin{align*}
    \log(A_2(N)) 
        &\leq C_5 \cdot \kl{ - N^{1- 1/2\eta'} + N^{1/2\eta'} \cdot \log(N)} \\
        &= - C_5 \cdot N^{1- 1/2\eta'} \cdot \underbrace{\kl{ 1 - \frac{\log(N)}{N^{1-1/\eta'}} }}_{\to 1}
\end{align*}

Hence, there is a positive constant $C_6>0$ and $\gamma \coloneqq 1- \tfrac{1}{2\eta'} > 0$, such that
\[ A_2(N) \leq \exp \kl{-C_6 \cdot N^\gamma} \]
for sufficiently large $N$. Since the sequence $\kl{\exp \kl{-C_6 \cdot N^\gamma}}_{N \in \NN}$ is summable it follows that $(A_2(N))_{N \in \NN}$ is summable, too.

All in all, this shows that the sequence $(\mu(\Omega^c_N))_{N \in \NN}$ is summable, and we conclude by the Borel--Cantelli Lemma.
\end{proof}
\section{Proof of the Main Theorem} \label{sec:proof}

We now have all the ingredients needed to prove our main volume estimate for volumes of canonical lifts of generic filling geodesics.

We choose a hyperbolic structure $X$ on our surface $S_g$ so that all pair of pants $P\in\P$ are isometric and are glued without twists and all the conditions of Section \ref{sec:random} are satisfied. Moreover, we let $0<\delta<1$ denote the Hausdorff dimension of the limit set of $P$.

\bthm\label{auxmainthm}
    Let $\eta > 1 (> \delta > 0)$ and let $F(t) \in O(t^{\delta/2\eta})$. Then for $\mu$-almost every $v\in \UT(X)$ there is $T' = T'(v)$ such that $\hat\gamma_v(t)$ is filling and
\be  F(\ell_X(\hat\gamma_v(t))) \leq \Vol(M_{\hat\gamma_v^p(t)})\ee
    for every $t \in R_\epsilon(v) \cap [T'(v),+\infty)$. Moreover, as long as $F$ is increasing we obtain that $ F(\ell_X(\hat\gamma^p_v(t)))\leq \Vol(M_{\hat\gamma_v^p(t)})$.
\ethm

\begin{proof}
    Let $v \in \UT(X)$ and $t \in R_\epsilon(v) \cap [t^*(v), + \infty)$. We consider the filling closed geodesic $\hat\gamma_v = \hat\gamma_v(t)$ and its primitive subcurve $\hat\gamma^p_v = \hat\gamma^p_v(t)$. Then, $M_{\hat\gamma_v^p}$ is hyperbolic and we now want to estimate its volume. We know that by \cite{RM20:LowerBound} we have;
\be \f{v_3}2 \sum_{P\in\P} \kl{\#\set{\hat\gamma_v^p\text{-arcs in }P} -3}\leq \Vol(M_{\hat\gamma_v^p}).\ee

Since the arcs are not counted with multiplicity we have that
\be \frac{1}{2} A_\P(\hat\gamma_v^p)= \sum_{P\in\P} \#\set{\hat\gamma_v^p\text{-arcs in }P}=  \sum_{P\in\P} \# \set{\hat\gamma_v\text{-arcs in }P}=\frac{1}{2} A_\P(\hat\gamma_v),\ee
where the factor $\tfrac{1}{2}$ is due to the fact the number of arcs in the projective bundle is at least half the number of arcs in the unit tangent bundle; see Remark \ref{rem:counting_directed}. 

By Proposition \ref{prop:arcs_pulled_tight} we may compare the number $A_{\P}(\hat{\gamma}_v)$ of homotopy classes of arcs for $\hat\gamma_v$  to the number $A_\P(\gamma_v[0,t])$ of homotopy classes of arcs in the geodesic ray $\gamma_v[0,t]$ of length $t$ starting from $v$:
\[ A_\P(\gamma_v[0,t]) - 6 \leq A_\P(\hat\gamma_v). \]

With the notation of section \ref{sec:counting} we may further estimate
\[C_{\floor{t}}(v) \leq  A_\P(\gamma_v[0,t]).\]

Thus, there are constants $A, B > 0$ (that are functions of $\chi(X)$ and $v_3$) such that
\[ A \cdot C_{\floor{t}}(v) - B \leq \Vol(M_{\hat\gamma^p_v}). \]

It follows from the Anosov Closing Lemma \ref{lem:AnosovClosing} that the length/period of $\hat\gamma_v(t)$ is close to $t$
\be \abs{\ell_X(\hat\gamma_v(t)) - t} \leq \rho, \ee
whence
\[ \floor{t} = \ell_X(\hat\gamma_v(t)) + \Delta \]
for some $\Delta \in [ -\rho -1,\rho + 1]$.

Applying Theorem \ref{mainestimate} we obtain that for any function $G(N) \in O(N^{\delta/2\eta})$ there is a subset of full measure $\Omega \subseteq \UT(X)$ satisfying that for every $v \in \Omega$ there is $N_0 = N_0(v) \in \NN$ such that 
\[ C_N(v) \geq G(N) \]
for every $N \geq N_0(v)$.

Putting everything together we obtain
\begin{equation} \label{eq:last_estimate}
    \Vol(M_{\hat\gamma_v^p}) \geq A \cdot G(\floor{t}) - B = A \cdot G(\ell_X(\hat\gamma_v(t)) + \Delta) - B 
\end{equation}
for every $v \in \Omega$, $t \in R_\epsilon(v) \cap [\max\{N_0(v),t^*(v) \}, + \infty)$ and some $\Delta \in [-\rho-1,\rho+1]$.

Note that for any $F(t) \in O(t^{\delta/2\eta})$ the function
\[ G(N) \coloneqq \tfrac{1}{A} \sup_{\Delta} F(N + \Delta) + \tfrac{B}{A}, \]
where the supremum is taken over $\Delta \in [-\rho-1,\rho+1]$,
is in $O(N^{\delta/2\eta})$, too. With this choice of $G$ it follows from (\ref{eq:last_estimate}) that
\[ \Vol(M_{\hat\gamma_v^p}) \geq  A \cdot G(\ell_X(\hat\gamma_v(t)) + \Delta) - B \geq F(\ell_X(\hat\gamma_v(t)))\]
for every $v \in \Omega$, $t \in R_\epsilon(v) \cap [T'(v), + \infty)$, where we set $T'(v) \coloneqq \max\{N_0(v),t^*(v) \}$. 

The case for $F$ increasing follows from the fact that $\ell_X(\hat\gamma_v(t))\geq \ell_X(\hat\gamma_v^p(t))$. This concludes the proof.\end{proof}

In combination with Lemma \ref{modelrel} this gives us:

\bthm\label{maintheorem}
    Let $\eta > 1 (> \delta > 0)$ and let $F(t) \in O(t^{\delta/2\eta})$. 
    
    Then the collection of filling, primitive closed geodesics satisfying
    \[  F(\ell_X(\gamma)) \leq \Vol(M_{\gamma}) \]
    has asymptotic density one.
\ethm

\begin{proof}
    A direct application of Lemma \ref{modelrel} yields that the collection
    \[ \mc{A}' \coloneqq \{ \gamma \in \mc{G} \,| \, \gamma \text{ is filling and } F(\ell_X(\gamma)) \leq \Vol(M_{\gamma^p}) \}\]
    has asymptotic density one since 
    \[ \Omega_{\mc{A}'} = \{ v \in \UT(X) \, | \, \exists T >0 \forall t \in R_\epsilon(v) \cap [T, \infty): \hat\gamma_v(t) \in \mc{A}' \} \]
    has full measure.
    
    Finally, observe that the collection of primitive geodesics $\mc{G}' \subseteq \mc{G}$ has asymptotic density one. Indeed, it is known that
    \[ \#\mc{G}(R) = \#\{ \gamma \in \mc{G} \, | \, \ell_X(\gamma) \leq R \} \sim \frac{e^R}{R}\]
    as $R \to + \infty$; see Remark \ref{rmk:margulis_count}.
    Any non-primitive geodesic traverses some other geodesic at least twice whence
    \[ \#\mc{G}(R) - \#\mc{G}'(R) \leq \# \mc{G}(R/2) \sim \frac{e^{R/2}}{R/2},\]
    and
    \[ \frac{\#\mc{G}(R) - \#\mc{G}'(R)}{\#\mc{G}(R)} \sim 2 \frac{e^{R/2}}{e^R} \to 0\]
    as $R \to +\infty$.

    Therefore, the intersection $\mc{A}' \cap \mc{G}'$ has asymptotic density one, too.
\end{proof}

\appendix

\section{Explicit construction of approximate convolutional identities} \label{sec:approx_id}

We call $\{ \psi_\epsilon \}_{0<\epsilon<\epsilon_0} \subseteq L^1(G)$ a \emph{family of approximate convolutional identities} if
\begin{enumerate}
    \item $\psi_\epsilon \geq 0$, and
    \item $\int_G \psi_\epsilon \, d\nu = 1$, and
    \item $\supp(\psi_\epsilon) \subseteq B_{\epsilon}(e)$,
\end{enumerate}
for all $0 < \epsilon < \epsilon_0$.

In the rest of this section we will construct an explicit family of approximate convolutional identities $\{ \psi_\epsilon \}_{0<\epsilon<\epsilon_0} \subseteq C_c^\infty(G)$, and estimate the $L^1$-norms of their derivatives (see Lemma \ref{lem:approx_id_L1}).

There is $\epsilon_0 > 0$ such that $\exp \colon B_{\epsilon_0}(0) \subseteq \mf{g} \to G$ is a diffeomorphism on its image. Let us denote by $\log \colon \exp(B_{\epsilon_0}(0)) \subseteq G \to B_{\epsilon_0}(0) \subseteq \mf{g}$ its inverse. We may choose a smooth bump function $\tilde{\psi}_{\epsilon_0} \in C_c^\infty(\mf{g})$, such that $\tilde{\psi}_{\epsilon_0} \geq 0$ and $\emptyset \neq \supp(\tilde{\psi}_{\epsilon_0}) \subseteq B_{\epsilon_0}(0)$. We set
\[ \bar{\psi}_{\epsilon_0}(g) 
\coloneqq 
    \begin{cases}
        \tilde{\psi}_{\epsilon_0}(\log(g)), & \text{if } g \in \exp(B_{\epsilon_0}(0)) \\
        0, &\text{else;}
    \end{cases}
\]
for every $g \in G$.
Further, we set
\[ M_{\epsilon_0} \coloneqq \int_G \bar{\psi}_{\epsilon_0} \, d\nu > 0, \]
and define
\[ \psi_{\epsilon_0} \coloneqq \frac{\bar{\psi}_{\epsilon_0}}{M_{\epsilon_0}} .\]
Then $\psi_{\epsilon_0}$ satisfies (1) and (2). 

Regarding (3) observe that for every $X \in B_{\epsilon_0}(0)$ the curve $c(t) \coloneqq \exp(t X), 0 \leq t \leq 1$ connects $e$ and $\exp(X) \in \exp(B_{\epsilon_0}(0))$. Therefore,
\[ d_G(\exp(X), e) \leq \ell(c) = \int_0^1 \| \dot{c}(t) \| \, dt  = \|X\| < \epsilon_0,  \]
and $\exp(B_{\epsilon_0}(0)) \subseteq B_{\epsilon_0}(e)$. Thus $\psi_{\epsilon_0}$ satisfies (3), too.

We will continue to define the other members of the family $\{\psi_{\epsilon}\}_{0<\epsilon<\epsilon_0}$ by rescaling $\psi_{\epsilon_0}$. Let $0<\epsilon<\epsilon_0$. We define
\[ \tilde{\psi}_{\epsilon}(X) \coloneqq \tilde{\psi}_{\epsilon_0} \kl{ \frac{\epsilon_0}{\epsilon} X} \]
for all $X \in \mf{g}$. By definition, $\supp(\tilde{\psi}_{\epsilon}) \subseteq B_{\epsilon}(0)$. 

As before, we set
\[ \bar{\psi}_{\epsilon}(g) 
\coloneqq 
    \begin{cases}
        \tilde{\psi}_{\epsilon}(\log(g)), & \text{if } g \in \exp(B_{\epsilon}(0)) \\
        0, &\text{else;}
    \end{cases}
\]
for every $g \in G$,
\[ M_{\epsilon} \coloneqq \int_G \bar{\psi}_{\epsilon} \, d\nu > 0, \]
and define
\[ \psi_{\epsilon} \coloneqq \frac{\bar{\psi}_{\epsilon}}{M_{\epsilon}}. \]
Then $\psi_{\epsilon}$ satisfies (1), (2) and (3) as before.

The following $L^1$-estimates is used in section \ref{subsubsect:orthogeodsets}.

\blem \label{lem:approx_id_L1}
    For every $d > 0$ there is a constant $K = K(d)>0$ such that 
    \[ \norm{ E_\alpha \cdot \psi_{\epsilon}}_1 \leq K \cdot \kl{\frac{1}{\epsilon}}^{d} \]
    holds for every $0 < \epsilon < \epsilon_0$ and every multi-index $\alpha$ with $\abs{\alpha} \leq d$.
\elem
\begin{proof}
    Recall that
    \[ \norm{E_\alpha \cdot \psi_\epsilon}_1= M_{\epsilon}^{-1} \norm{E_\alpha \cdot \bar{\psi}_{\epsilon}} _1.\]
    
    Let us first estimate $M_\epsilon = \int_G \bar{\psi}_\epsilon \, d\nu$, and see that $M_\epsilon \asymp \epsilon^3$.
    
    We compute
    \begin{align*}
        M_\epsilon &= \int_G \bar{\psi}_\epsilon \, d\nu
                    = \int_{\exp(B_{\epsilon}(0))} \tilde{\psi}_\epsilon(\log(g)) \, d\nu(g) \\
                &= \int_{B_{\epsilon}(0)} \tilde{\psi}_\epsilon(X) \, d (\log_* \nu)(X) \\
                &= \int_{B_{\epsilon}(0)} \tilde{\psi}_\epsilon(X) \cdot \frac{d(\log_* \nu)}{d\lambda}(X) \, d\lambda(X),
    \end{align*}
    where we denote by $\lambda$ the Lebesgue measure on $\mf{g}$ and $\frac{d(\log_* \nu)}{d\lambda}$ denotes the Radon-Nikodym derivative of $\log_*\nu$ with respect to $\lambda$. This is well-defined since the measure $\log_*(\nu)$ is clearly of Lebesgue class. Moreover, because it is the push-forward of the Haar measure $\nu$ via the diffeomorphism $\log \colon \exp(B_{\epsilon_0}(0)) \to B_{\epsilon_0}(0)$, the Radon-Nikodym derivative is smooth. Consequently there are uniform upper and lower bounds
    \[ 0< \inf_X \frac{d(\log_* \nu)}{d\lambda}(X), \qquad \sup_{X} \frac{d(\log_* \nu)}{d\lambda}(X) < + \infty,\]
    where the infimum and the supremum are taken over the compact set $\supp(\psi_{\epsilon_0}) \subseteq B_{\epsilon_0}(0)$.
    
    Hence,
    \begin{align*}
        M_\epsilon &\asymp \int_{B_\epsilon(0)} \tilde{\psi}_{\epsilon}(X) \, d\lambda(X)
            = \kl{ \frac{\epsilon}{\epsilon_0}}^3 \cdot \int_{B_\epsilon(0)} \tilde{\psi}_{\epsilon_0} \kl{ \frac{\epsilon_0}{\epsilon} X } \cdot \kl{ \frac{\epsilon_0}{\epsilon}}^3 \, d\lambda(X)\\
            &= \kl{ \frac{\epsilon}{\epsilon_0}}^3 \cdot \int_{B_{\epsilon_0}(0)} \tilde{\psi}_{\epsilon_0}  \, d\lambda \asymp \epsilon^3,
    \end{align*}
    as asserted.
    
    We turn to $\norm{E_\alpha \cdot \bar{\psi}_{\epsilon}}_1$. For a single $E_i$ we have
    \[ (E_i \cdot \bar{\psi}_{\epsilon})(g) = \frac{d}{dt}\bigg|_{t=0} \bar{\psi}_\epsilon(\exp(- t E_i) g) \]
    for every $g \in G$. Thus, if we define by $\tilde{E}_i$ the \emph{right-invariant} vector field given by $\tilde{E}_i(e) = - E_i$, then we can rewrite
    \[ E_\alpha \cdot \bar{\psi}_\epsilon = \tilde{E}_\alpha \cdot \bar{\psi}_{\epsilon} \]
    where the right hand side is understood as the usual derivation with respect to vector fields.
    
    On $\supp(\bar{\psi}_\epsilon) \subseteq \exp(B_{\epsilon}(0))$ we obtain
    \[ \tilde{E}_\alpha \cdot \bar{\psi}_{\epsilon} = \tilde{E}_\alpha \cdot (\tilde{\psi}_\epsilon \circ \log) = \kl{\log_*(\tilde{E})_\alpha \cdot \tilde{\psi}_\epsilon} \circ \log. \]
    Here we denote by $\log_*(\tilde{E})_i \in \Vect(B_{\epsilon}(0))$ the push-forward of the vector field $\tilde{E}_i$ along the map $\log \colon \exp(B_{\epsilon}(0)) \to B_{\epsilon}(0) \subseteq \mf{g}$.
    
    We may express $\log_*(\tilde{E})_i$ with respect to the standard vector fields $\partial_1, \ldots, \partial_3$ of $\mf{g} \cong \RR^3$.
    In this way, we obtain polynomials $P_{\alpha} = \sum_{0\leq \abs{\beta} \leq \abs{\alpha}} c^{\beta}_\alpha \cdot \partial_\beta \in \mc{U}(\Vect(B_{\epsilon}(0)))$ such that
    \[\log_*(\tilde{E})_\alpha = P_{\alpha}, \]
    where $c^{\beta}_\alpha \colon B_{\epsilon}(0) \to \RR $ are smooth coefficient functions. If we apply these to $\tilde{\psi}_\epsilon$, we obtain
    \begin{align*}
        \kl{\log_*(\tilde{E})_\alpha \cdot \tilde{\psi}_\epsilon}(X)
            &= \kl{P_\alpha \cdot \tilde{\psi}_\epsilon}(X)
            = \sum_{0\leq \abs{\beta} \leq \abs{\alpha}} c^{\beta}_\alpha(X) \cdot \partial_\beta \tilde{\psi}_\epsilon(X)\\
            &=\sum_{0\leq \abs{\beta} \leq \abs{\alpha}} c^{\beta}_\alpha(X) \cdot \partial_\beta \kl{ \tilde{\psi}_{\epsilon_0}\kl{ \frac{\epsilon_0}{\epsilon} X } }\\
            &=\sum_{0\leq \abs{\beta} \leq \abs{\alpha}} \kl{ \frac{\epsilon_0}{\epsilon}}^{\abs{\beta}} c^{\beta}_\alpha(X) \cdot \kl{\partial_\beta  \tilde{\psi}_{\epsilon_0}}\kl{ \frac{\epsilon_0}{\epsilon} X }
    \end{align*}
    for every $X \in B_{\epsilon}(0)$.
    
    With these preliminary computations we obtain
    \begin{align*}
        \norm{E_\alpha \cdot \bar{\psi}_{\epsilon}}_1 
            &= \int_G \abs{E_\alpha \cdot \bar{\psi}_{\epsilon}} \, d\nu \\
            &= \int_G \abs{\kl{\log_*(\tilde{E})_\alpha \cdot \tilde{\psi}_\epsilon} \circ \log} \, d\nu \\
            &= \int_{B_{\epsilon}(0)}\abs{ \kl{\log_*(\tilde{E})_\alpha \cdot \tilde{\psi}_\epsilon}(X)} \cdot \abs{\frac{d\log_*(\nu)}{d\lambda}(X)}  \, d\lambda(X) \\
            &\ll \int_{B_{\epsilon}(0)}\abs{ \kl{\log_*(\tilde{E})_\alpha \cdot \tilde{\psi}_\epsilon}(X)} \, d\lambda(X)\\
            &\ll \sum_{0\leq \abs{\beta} \leq \abs{\alpha}} \abs{ \frac{\epsilon_0}{\epsilon}}^{\abs{\beta}} \int_{B_\epsilon(0)}  \abs{c^{\beta}_\alpha(X)} \cdot \abs{\kl{\partial_\beta  \tilde{\psi}_{\epsilon_0}}\kl{ \frac{\epsilon_0}{\epsilon} X }} \, d \lambda(X) \\
            &= \sum_{0\leq \abs{\beta} \leq \abs{\alpha}} \abs{ \frac{\epsilon_0}{\epsilon}}^{\abs{\beta}-3} 
            \underbrace{\int_{B_{\epsilon_0}(0)}  \abs{c^{\beta}_\alpha \kl{ \frac{\epsilon}{\epsilon_0} X }} \cdot \abs{\kl{\partial_\beta  \tilde{\psi}_{\epsilon_0}}(X)} \, d \lambda(X)}_{\text{uniformly bounded}} \\
            &\ll \kl{ \frac{\epsilon_0}{\epsilon} }^{\abs{\alpha}-3}.
    \end{align*}
    
    In conjunction with $M_\epsilon \asymp \epsilon^3$ we get that
    \[ \norm{E_\alpha \cdot \psi_{\epsilon}}_1 = M_\epsilon^{-1} \norm{E_\alpha \cdot \bar{\psi}_{\epsilon}}_1 \ll \kl{ \frac{\epsilon_0}{\epsilon} }^{\abs{\alpha}} 
    \leq \kl{ \frac{\epsilon_0}{\epsilon} }^{d} \ll \kl{ \frac{1}{\epsilon} }^{d}  .\]\end{proof}
\section{Exponential higher-order mixing}
\label{app:mixing}
Let $Y\coloneqq \UT(X) \cong\Gamma\backslash G$ denote the unit tangent bundle of a closed hyperbolic surface $X=\Gamma\backslash \mathbb H^2$ with $G=\PSL_2\R$ and $\Gamma$ a torsion-free cocompact subgroup. We set 
\[a_t\coloneqq \begin{pmatrix} e^{t/2} & 0\\ 0 & e^{-t/2} \end{pmatrix} \quad \text{and} \quad u_s\coloneqq \begin{pmatrix} 1 & s\\ 0 & 1 \end{pmatrix}\]
for all $t, s \in \RR$. 

The geodesic flow on $Y$ through $x = \Gamma g \in \Gamma \lquot G$ is given by 
\[x a_t \coloneqq \Gamma g a_t \quad \forall t \in \RR,\] 
and the horocycle flow is given by 
\[x u_s \coloneqq \Gamma g u_s \quad \forall s \in \RR.\] 

Moreover, we denote by $T \colon Y \to Y$ the time-one map of the geodesic flow, i.e.\ $ T(x) \coloneqq xa_1 $ for all $x \in Y$. Moreover, we will abbreviate integration with respect to the unique invariant probability measure (the normalized Liouville measure) on $Y$ by $dx$ or $dy$.

The following formula will be useful:
\be a_{t}u_sa_{-t}=\begin{pmatrix} e^{t/2} & 0\\ 0 & e^{-t/2} \end{pmatrix}\begin{pmatrix} 1 & s\\ 0 & 1 \end{pmatrix}\begin{pmatrix} e^{-t/2} & 0\\ 0 & e^{t/2} \end{pmatrix}=\begin{pmatrix} 1 & e^ts\\ 0 & 1 \end{pmatrix} = u_{e^ts}\ee
for all $t,s \in \RR$.

Moreover, we will use the following result due to Burger \cite{Bu1990}:

\bthm[{\cite[Theorem 2 (C)]{Bu1990}}] \label{thm:Burger}
    Let $\lambda_1<0$ be the first non-zero eigenvalue of the Laplacian of $X$ and let $0<\alpha\leq \f 12$ satisfy $\alpha(\alpha-1)\geq\lambda_1$. Then, we have for all $f\in C_c^\infty(Y)$ and $T\geq 1$
    \be \sup_{x\in Y}\abs{\f 1{2T}\int_{-T}^T f(xu_t)dt-\int_Y f(y)dy }\leq c \f{T^{-\alpha}-T^{\alpha-1}}{1-2\alpha} \Sob(f),  \ee
    where $\Sob(f)=\Sob_3(f)$ denotes the degree $d=3$ Sobolev norm; see section \ref{sec:sobolev_norms}.
\ethm
\bcor\label{cor:Burger}
    In particular, we obtain that there exists $C>0$, $\alpha>0$ such that for all $f\in  C_c^\infty(Y)$, for all $T\geq 1/2$, and for all $x\in Y$
\be \abs{\f 1T \int_0^Tf(xu_t)dt-\int_Y f(y)dy}\leq  C T^{-\alpha} \Sob(f).\ee
\ecor

We will now deduce the following exponential $k$-mixing result. Note that similar results have been deduced in more generality before; see \cite{expmixEinsiedler}. However, for our application it is important that the given constants  depend neither on the functions nor on the \emph{number of functions}. 

\begin{samepage}
\bthm[Exponential $k$-mixing] \label{thm:expmixing}
    Let $\alpha$ be as in Theorem \ref{thm:Burger}. There are constants $D, L \geq 1$ and $K_0 > 0$, such that 
    \[ \abs{ \int_{Y} \prod_{j=1}^n \kl{f_j \circ T^{-j \cdot k} } \, dx - \prod_{j=1}^n \kl{ \int_{Y} f_j \, dx }} \leq L \cdot (n-1) \cdot D^n  \cdot e^{- \frac{\alpha}{1+\alpha} k} \cdot \prod_{j=1}^n \Sob(f_j)\]
    for all $n\in\NN$, for all $f_1, \ldots, f_n \in C_c^\infty(\UT(X))$ and all integers $k \geq K_0$. Here $\Sob(f)=\Sob_3(f)$ denotes the degree $d=3$ Sobolev norm as before.
\ethm
\end{samepage}

Clearly, the above result for the geodesic flow translates into the following equivalent exponential $k$-mixing result for the time-reversed geodesic flow.
\bcor \label{cor:expmixing}
    With the same constants as above we have that
   \[ \abs{ \int_{Y} \prod_{j=1}^n \kl{f_j \circ T^{j \cdot k} } \, dx - \prod_{j=1}^n \kl{ \int_{Y} f_j \, dx }} \leq L \cdot (n-1) \cdot D^n  \cdot e^{- \frac{\alpha}{1+\alpha} k} \cdot \prod_{j=1}^n \Sob(f_j)\]
    for all $n\in\NN$, for all $f_1, \ldots, f_n \in C_c^\infty(\UT(X))$ and all integers $k \geq K_0$.
\ecor

\begin{proof}[Proof of Corollary \ref{cor:expmixing}]
    Indeed, the measure $\mu$ is invariant under $T^{- n \cdot k} \colon Y \to Y$. Thus,
    \begin{align*}
        &\quad \abs{ \int_{Y} \prod_{j=1}^n \kl{f_j \circ T^{j \cdot l} } \, d x - \prod_{j=1}^n \kl{ \int_{Y} f_j \, dx }}\\
        &=\abs{ \int_{Y} \prod_{j=1}^n \kl{f_j \circ T^{-(n-j) \cdot k} } \, d x - \prod_{j=1}^n \kl{ \int_{Y} f_j \, dx}}\\
        &\leq L \cdot (n-1) \cdot D^n  \cdot e^{- \frac{\alpha}{1+\alpha} k} \cdot \prod_{j=1}^n \Sob(f_j)
    \end{align*}
    after reversing the indices $f_j \leftrightarrow f_{n-j}$.
\end{proof}

We thank Einsiedler for suggesting the strategy for the following proof of Theorem \ref{thm:expmixing}.

\begin{proof}[Proof of Theorem \ref{thm:expmixing}.] We will prove the theorem by induction on the number of functions $n$. The induction base $n=2$ is the classical result that the geodesic flow is exponentially mixing. However, we include a proof here, since it illustrates the idea for the induction step.

Let $f_1,f_2\in C_c^\infty(Y)$, $k\in\N$, and $T \geq 1/2$. We will choose $T$ appropriately later on. In a first step we split the error into two terms $\Delta, \Delta'$ that we will then bound separately.

\begin{align*}
&\quad \abs{\int_Y f_1(x)f_2(T^{-k}(x))dx-\int_Y f_1(x)dx\int_Y f_2(x)dx}\\
&=\abs{\int_Yf_1(x)f_2(xa_{-k})dx-\int_Y f_1(x)dx\int_Y f_2(x)dx}\\
&=\abs{\f 1 T\int_0^T\int_Y f_1(x u_s)f_2(x u_sa_{-k})dxds-\int_Yf_1(x)dx\int_Yf_2(x)dx}\\
&\leq \Delta+\Delta',
\end{align*}
where
\begin{align*}
    \Delta &\coloneqq \abs{\f 1 T\int_0^T\int_Y f_1(xu_s)f_2(xu_sa_{-k})dxds-\f 1T \int_0^T\int_Yf_1(x)f_2(xu_sa_{-k})dx},\\
    \Delta' &\coloneqq \abs{\int_Y f_1(x)\left(\f 1 T\int_0^T f_2(xu_sa_{-k})ds\right)dx-\int_Yf_1(x)\left(\int_Yf_2(y)dy\right)dx}.
\end{align*}

Let us first bound $\Delta$. We denote by
\[ U \coloneqq \begin{pmatrix} 0 & 1 \\ 0 & 0 \end{pmatrix} = \frac{d}{dt}\bigg\vert_{t=0} \exp(u_t) \in \mf{sl}_2(\RR). \]
By the Sobolev Embedding Theorem \ref{thm:SET} there is a constant $D_0>0$ such that
\begin{align}
    \abs{f(xu_t) - f(x)} 
        = \abs{(\lambda_{\exp(t U)}f)(x) - f(x)} 
        \leq \int_0^t \abs{(U \cdot f)(x u_s)} \, ds 
        \leq t \cdot \| U \cdot f \|_\infty  
        \leq D_0 \cdot \Sob(f) \cdot t
        \label{eq:LinftyEst}
\end{align}
for all $f \in C_c^\infty(\Gamma \lquot G)$, $t \in \RR$, $x \in \Gamma \lquot G$.

Thus, we may estimate
\begin{align*}
    \Delta 
    &\leq \f 1 T\int_0^T\int_Y\underbrace{\abs{f_1(xu_s)-f_1(x)}}_{\leq D_0\Sob(f_1)T}\underbrace{\abs{f_2(x u_s a_{-k})}}_{\leq \norm{f_2}_\infty \leq \Ksob \Sob(f_2)}dxds\leq D_0 \Ksob \Sob(f_1) \Sob(f_2)T.
\end{align*}

Using Corollary \ref{cor:Burger} we may estimate $\Delta'$ as follows.
\begin{align*}
\Delta' &= \abs{\int_Y f_1(x)\left(\f 1 T\int_0^T f_2(x a_{-k} a_k u_s a_{-k})ds\right)dx-\int_Yf_1(x)\left(\int_Yf_2(y)dy\right)dx}\\
 &\leq \int_Y \underbrace{\abs{f_1(x)}}_{\leq \norm{f_1}_\infty}\bigg |\underbrace{\f 1 T\int_0^Tf_2(xa_ku_{e^ks})ds}_{=\f 1{e^kT}\int_0^{e^kT} f_2(xa_ku_t)dt}-\int_Yf_2(y)dy \bigg| dx\\
 &\stackrel{(\diamond)}{\leq} \Ksob \Sob(f_1) \cdot C\cdot \Sob(f_2)(e^kT)^{-\alpha} \leq C \Ksob \Sob(f_1)\Sob(f_2) T^{-\alpha} e^{-\alpha k},
\end{align*}
where we assumed in $(\diamond)$ that $e^k T \geq 1/2$.

Altogether
\[\Delta + \Delta' \leq D^2 \Sob(f_1) \Sob(f_2) \underbrace{(T + T^{-\alpha} e^{-\alpha k})}_{\eqqcolon F_2(T)}\]
for $D \coloneqq \max(1,C,\Ksob,D_0)$.

The function $F_2$ attains its global minimum at its unique critical point
\[ T_2 \coloneqq \alpha^{1/(1+\alpha)} e^{k/(1+\alpha)},\]
with
\[ F_2(T_2) = \underbrace{\kl{\alpha^{1/(1+\alpha)} + \alpha^{-\alpha/(1+\alpha)}}}_{\leq L} e^{ -\frac{\alpha}{1+\alpha} k},\]
where we set $L \coloneqq D_1^{\frac{\alpha}{1+\alpha}} \kl{\alpha^{1/(1+\alpha)} + \alpha^{-\alpha/(1+\alpha)}}$ with $D_1 \coloneqq 1/(1-e^{-1}) > 1$. (The relevance of these choices will become apparent in the induction step.)

Observe that
\[ e^k T_2 = \alpha^{1/(1+\alpha)} e^{k/(1+\alpha)} \geq \frac{1}{2} \quad \iff \quad k \geq - \log(\alpha) - (1+\alpha) \log(2). \]
In particular, $e^k T_2 \geq 1/2$ holds if
\[k \geq K_0 \coloneqq  \log(D_1)-\log(\alpha)  - (1+\alpha) \log(2).\]

All in all,
\[ \Delta + \Delta' \leq L D^2 \Sob(f_1) \Sob(f_2) e^{ -\frac{\alpha}{1+\alpha} k} \]
for all $k \geq K_0$.
This concludes the proof of the induction base.

We now proceed with the induction step, going from $n-1$ to $n$ functions. Let $f_1,\dotsc,f_n\in C_c^\infty(Y)$, $T\geq 1/2$, and $k\in\N$. As for the induction base we will determine $T$ later on. This time we will split the error in three terms $\Delta_1, \Delta_2, \Delta_3$ that we will then bound separately.

\begin{align*}
    &\quad\abs{ \int_Y f_1(x) \cdot \ldots \cdot f_n(x a_{-k(n-1)}) d x - \prod_{j=1}^n \int_Y f_j(x) dx } \\
    &= \abs{ \frac{1}{T} \int_0^T \int_Y f_1(xu_s) \cdot \ldots \cdot f_n(xu_s a_{-k(n-1)}) dx ds - \prod_{j=1}^n \int_Y f_j(x) dx } \\
    &\leq \Delta_1 + \Delta_2 + \Delta_3,
\end{align*}
where we set
\begin{align*}
    \Delta_1 
        &\coloneqq \bigg \vert \frac{1}{T} \int_0^T \int_Y f_1(xu_s) \cdot \ldots \cdot f_n(xu_s a_{-k(n-1)}) dx ds \\ 
        & \qquad - \frac{1}{T} \int_0^T \int_Y f_1(x) \cdot \ldots \cdot f_{n-1}(x a_{-k(n-2)}) \cdot f_{n}(xu_s a_{-k(n-1)}) dx ds \bigg \vert, \\
    \Delta_2
        &\coloneqq \bigg \vert \frac{1}{T} \int_0^T \int_Y f_1(x) \cdot \ldots \cdot f_{n-1}(x a_{-k(n-2)}) \cdot f_{n}(xu_s a_{-k(n-1)}) dx ds \\
        &\qquad - \int_Y f_1(x) \cdot \ldots \cdot f_{n-1}(x a_{-k(n-2)}) \cdot \kl{ \int_Y f_n(y) dy} dx \bigg \vert,\\
    \Delta_3 
        & \coloneqq \bigg \vert \int_Y f_1(x) \cdot \ldots \cdot f_{n-1}(x a_{-k(n-2)}) \cdot \kl{ \int_Y f_n(y) dy} dx \\
        &\qquad - \int_Y f_1(x) dx \cdot \ldots \cdot \int_Y f_n(x) dx \bigg \vert.
\end{align*}

Let us start with $\Delta_3$. By the induction hypothesis we have that
\begin{align*}
    \Delta_3 
        & \leq \abs{ \int_Y f_1(x) \cdot \ldots \cdot f_{n-1}(x a_{-k(n-2)}) dx - \int_Y f_1(x) dx \cdot \ldots \cdot \int_Y f_{n-1}(x) dx } \cdot \int_Y \abs{f_n(x)} dx \\
        &\leq \norm{f_n}_\infty L (n-2) D^{n-2} e^{ - \frac{\alpha}{1 + \alpha } k } \prod_{j=1}^{n-1} \Sob(f_j) \\
        &\leq L (n-2) D^{n-1}  e^{ - \frac{\alpha}{1 + \alpha } k } \prod_{j=1}^{n} \Sob(f_j).
\end{align*}

Next, we consider $\Delta_2$. By Corollary \ref{cor:Burger} we obtain
\begin{align*}
    \Delta_2
        &\leq \int_Y \abs{f_1(x)} \cdot \ldots \cdot \abs{f_{n-1}(x a_{-k(n-2)})} \cdot 
        \abs{ \frac{1}{T} \int_0^T f_n(x a_{-k(n-1)} u_{e^{(n-1)k}s}) ds - \int_Y f_n(y) dy } dx \\
        &\leq \Ksob^{n-1} \prod_{j=1}^{n-1} \Sob(f_j) \cdot \int_Y \abs{ \frac{1}{e^{(n-1)k}T} \int_0^{e^{(n-1)k}T} f_n(x a_{-k(n-1)} u_{s}) ds - \int_Y f_n(y) dy } dx \\
        &\leq C \Ksob^{n-1} e^{ -\alpha k (n-1) } T^{-\alpha} \prod_{j=1}^n \Sob(f_j) \\
        &\leq D^n e^{ -\alpha k (n-1) } T^{-\alpha} \prod_{j=1}^n \Sob(f_j),
\end{align*}
where we assumed that $e^{(n-1)k} T \geq 1/2$.

Finally, we turn to $\Delta_1$. Using a telescope sum argument one can show that
\begin{align}
    \abs{ \prod_{j=1}^{n-1} x_j - \prod_{j=1}^{n-1} y_j} \leq \sum_{i=1}^{n-1} \abs{x_i - y_i} \prod_{j\neq i} \max \kl{\vert x_j \vert ,\vert y_j \vert} 
    \label{eq:productEst}
\end{align}
for all $x_1, \ldots, x_{n-1}, y_1, \ldots, y_{n-1} \in \RR$. We will use (\ref{eq:productEst}) with
\[
    x_i \coloneqq f_i(xu_s a_{-k(i-1)}) = f_i(x a_{-k(i-1)} u_{e^{k(i-1)}s}), \qquad y_i \coloneqq f_i(x a_{-k(i-1)}),
\]
for all $i=1,\ldots, n-1$. Notice that in this case $\max \kl{\vert x_j \vert,\vert y_j \vert} \leq \norm{f_j}_\infty$ for all $j=1,\ldots, n-1$. Moreover,
\[ \abs{x_i - y_i} = \abs{ f_i(x a_{-k(i-1)} u_{e^{k(i-1)}s}) -  f_i(x a_{-k(i-1)}) } \leq D_0 \Sob(f_i) e^{k(i-1)} T \]
for all $0 \leq s \leq T$, $i=1, \ldots, n-1$ by (\ref{eq:LinftyEst}).

With these observations we may estimate $\Delta_1$ as follows.
\begin{align*}
    \Delta_1 
        &\leq \frac{1}{T} \int_0^T \int_Y  
         \abs{ f_1(xu_s) \cdot \ldots \cdot f_{n-1}(x a_{-k(n-2)} u_{ e^{k(n-2)}s }) 
             - f_1(x) \cdot \ldots \cdot f_{n-1}(x a_{-k(n-2)}) } \\
            & \qquad \qquad \cdot \abs{f_n(x u_s a_{ -k(n-1)  })} dx ds \\
        &\leq \norm{f_n}_\infty \frac{1}{T} \int_0^T \int_Y  \sum_{i=1}^{n-1} D_0 \Sob(f_i) e^{k(i-1)} T \prod_{j\neq i} \norm{f_j}_\infty dx ds \\
        &\leq D^n \frac{e^{k(n-1)}-1}{e^k-1} T \prod_{j=1}^n \Sob(f_j)  \\
        &\leq D_1 D^n  e^{k(n-2)} T \prod_{j=1}^n \Sob(f_j),
\end{align*}
where we used that $D_1 = 1/(1-e^{-1}) \geq 1/(1-e^{-k})$.

Combining these estimates we get
\begin{align*}
    &\quad \Delta_1 + \Delta_2 + \Delta_3 \\ 
    &\leq  D_1 D^n  e^{k(n-2)} T \prod_{j=1}^n \Sob(f_j) + D^n e^{ -\alpha k (n-1) } T^{-\alpha} \prod_{j=1}^n \Sob(f_j) + L (n-2) \underbrace{D^{n-1}}_{\leq D^n}  e^{ - \frac{\alpha}{1 + \alpha } k } \prod_{j=1}^{n} \Sob(f_j)\\
    &= D^n \prod_{j=1}^n \Sob(f_j) \cdot \underbrace{\kl{ D_1 T e^{k(n-2)} + e^{ -\alpha (n-1) k } T^{-\alpha} + L (n-2) e^{ - \frac{\alpha}{1+\alpha} k } }}_{ \eqqcolon F_n(T) }.
\end{align*}
The function $F_n$ attains its global minimum at its unique critical point
\[ T_n \coloneqq \alpha^{\frac{1}{1+\alpha}} D_1^{-\frac{1}{1+\alpha}} \exp \kl{ -k (n-1) + \frac{k}{1+\alpha} } \]
with
\[ F_n(T_n) = \kl{ \underbrace{ D_1^{\frac{\alpha}{1+\alpha}} \kl{\alpha^{\frac{1}{1+\alpha}} + \alpha^{-\frac{\alpha}{1+\alpha}}}}_{\leq L} + L (n-2)} e^{ - \frac{\alpha}{1+\alpha} k } \leq L (n-1) e^{ - \frac{\alpha}{1+\alpha} k }.\]

Also, notice that
\[ T_n e^{(n-1)k} = \alpha^{\frac{1}{1+\alpha}} D_1^{-\frac{1}{1+\alpha}} \exp\kl{ \frac{k}{1+\alpha} } \geq \frac{1}{2} \quad \iff \quad k \geq  \log(D_1) - \log(\alpha) - (1+\alpha) \log(2) = K_0. \] 
All in all,
\[ \Delta_1 + \Delta_2 + \Delta_3 \leq L (n-1) D^n e^{ - \frac{\alpha}{1+\alpha} k } \prod_{j=1}^n \Sob(f_j) \]
for all $k \geq  K_0$. This concludes the proof.

\end{proof}
\section{Tube Lemma}

Let $M^{n+1}$ be an oriented Riemannian manifold with bounded sectional curvature $\abs{K} \leq b$. Let $S^n \subseteq M$ be a codimension $1$ orientable hypersurface with normal field $N \colon S \to TS^{\perp} \subseteq TM$. We denote by $\omega$ the volume form on M, and by $\bar{\omega}$ the induced volume form on $S$. 

At every point $x \in S$ we can find an orthonormal basis $\mc{B} = \{v_1, \ldots, v_n\} \subseteq T_x S$ of its tangent space. We consider the symmetric matrix
\[ \bar{N}(x) = \kl{ \langle\nabla_{v_i} N, v_j \rangle }_{1 \leq i,j \leq n}. \]
Its eigenvalues $\kappa_1(x), \ldots, \kappa_n(x)$ are called the principal curvatures at $x$ of $S$ in $M$. Note that in this case $\| \bar{N}(x) \|_2 = \max_{i=1,\ldots,n} \abs{\kappa_i(x)}$, and this quantity is independent of the basis $\mc{B}$.

We can parametrize the $\epsilon$-tubular neighborhood about $S$ via the map
    \[ \varphi \colon (-\epsilon, \epsilon) \times S \to M, (t,x) \mapsto \Exp_x(t \cdot N(x)). \]
There is a smooth function $\rho \in C^\infty((-\epsilon,\epsilon)\times S)$ such that
\[ \phi^* \omega = \rho \cdot dt \wedge \bar{\omega}. \]

\blem[Tube Lemma] \label{lem:tube_lemma}
    Let $\hat{b} \coloneqq \max(1,b)$.
    There is a constant $C>0$ such that 
    \[ \abs{\rho(t,x)} \leq C \cdot (1+ \norm{\bar{N}(x)}_2)^n \cdot e^{\hat{b} n t} \]
    for all $(t,x) \in (-\epsilon,\epsilon) \times S$. 
    
    In particular, if there is a uniform bound $\kappa > 0$ on the principal curvatures, i.e.\ 
    \be \norm{\bar{N}(x)}_2 = \max_{i=1,\ldots,n} \abs{\kappa_i(x)}\leq \kappa\ee for every $x \in S$, then 
    \[\| \rho \|_\infty \leq C \cdot (1+ \kappa)^n \cdot e^{\hat{b} n \epsilon}, \]
    such that
    \[ \vol_M(\phi(B)) \leq C \cdot (1+ \kappa)^n \cdot e^{\hat{b} n \epsilon} \cdot \vol_{(-\epsilon,\epsilon)\times S}(B) \]
    for every Borel set $B \subseteq (-\epsilon,\epsilon)\times S$ with respect to the volume form $dt\wedge \bar{\omega}$.
\elem
\begin{proof}
 Let $x\in S$ be a point and denote by $c(t)\coloneqq \Exp_x(t\cdot N(x)) $ the normal geodesic at $x$. Let $v_0\coloneqq \dot c(0),v_1,\dotsc,v_n \in T_xS\subset T_xM$ be a positively oriented orthonormal basis and let $E_i(t)$ be the parallel vector field along $c(t)$ such that: $E_i(0)=v_i$. In particular, $E_0(t)=\dot c(t)$. Further, let $(t,x)\in(-\epsilon,\epsilon)\times S$. Then, $\partial_t\in T_t(-\epsilon,\epsilon)\subset T_{(t,x)}((-\epsilon,\epsilon)\times S)=T_t(-\epsilon,\epsilon) \oplus T_x S $, $v_i\in T_x S\subset T_{(t,x)}((-\epsilon,\epsilon)\times S)$, and $(\partial_t,v_1,\dotsc,v_n)$ is a positively oriented orthonormal basis for $T_{(t,x)}((-\epsilon,\epsilon)\times S)$. Finally, choose $\gamma_i:\R\rar S$ such that $\dot\gamma_i(0)=v_i$. 
 
 Then 
 \begin{align*}
     \abs{\rho(t,x)}&=\rho(t,x)(dt\wedge \bar\omega)(\partial_t,v_1,\dotsc,v_n)=\\
     &=(\phi^*\omega)(\partial_t,v_1,\dotsc,v_n)=\\
     &=\omega(d\phi(\partial_t),d\phi(v_1),\dotsc,d\phi(v_n)).
 \end{align*}
 Note that
 \be d\phi(\partial_t)=\difrac{}{t} \Exp_x(t\cdot N(x))=\dot c(t)=E_0(t), \qquad d\phi(v_i)=\difrac{}{s}\bigg\vert_{s=0} \Exp_{\gamma_i(s)}(t\cdot N(\gamma_i(s)))\eqqcolon J_i(t).\ee
As a geodesic variation $J_i(t)$ is a Jacobi field, such that it satisfies the Jacobi equation
\be
J''_i=-R(\dot c, J_i)\dot c=-R(E_0,J_i)E_0
\ee
with the initial conditions
\begin{align*}
    J_i(0)&=\difrac{}{s}\bigg\vert_{s=0} \Exp_{\gamma_i(s)}(0\cdot N(\gamma_i(s)))=\difrac{}{s}\bigg\vert_{s=0} \gamma_i(s) =v_i=E_i(0), \quad \text{and}\\
    J_i'(0)&=\f{D}{\partial t}\difrac{}{s}\bigg\vert_{t,s=0} \Exp_{\gamma_i(s)}(t\cdot N(\gamma_i(s)))\\
        &=\f{D}{\partial s}\difrac{}{t}\bigg\vert_{t,s=0} \Exp_{\gamma_i(s)}(t\cdot N(\gamma_i(s)))=\f D {\partial s}\bigg\vert_{s=0} N(\gamma_i(s))=\nabla_{v_i}N.
\end{align*}

In parallel coordinates $J_i(t)=\sum_j\alpha_i^j(t)\cdot E_j(t)$ for some function $\alpha_i^j(t)$. Therefore, $J''_i=\sum_j\ddot \alpha_i^j E_j$. Then,
\begin{align*}
    -\langle R(E_0,J_i)E_0,E_j\rangle&=-\langle R(E_0,\sum_k \alpha_i^k E_k)E_0,E_j\rangle=\\
    &=-\sum_k \alpha_i^k\langle R(E_0,E_k)E_0,E_j\rangle
\end{align*}
and we let $R_{kj}\coloneqq\langle R(E_0,E_k)E_0,E_j\rangle $. By curvature symmetries we get $R_{kj}=R_{jk}$. Using matrix notation $A_{ij}\coloneqq \alpha_i^j$ we have:
\begin{itemize}
    \item $A''=-AR$;
    \item $A(0)=\id$;
    \item $A'(0)=\left( \langle \nabla_{v_i} N,v_j\rangle\right)_{ij}=\bar N(x)$.
\end{itemize}
We reduce this to a first order system as follows:
\be (A'\;\; A'')=(A\;\; A')\cdot \begin{pmatrix} 0 & -R\\ \id & 0\end{pmatrix}.\ee

  It follows from ODE theory that
  \begin{align}
      \norm{(A(t) \;\; A'(t))}_2\leq \norm{(A(0) \; \; A'(0))}_2\cdot\exp \left(t\norm{\begin{pmatrix} 0 & -R(t)\\ \id & 0\end{pmatrix}}_2\right).
      \label{eq:ODE}
  \end{align}
  Notice that 
  \[\norm{\begin{pmatrix} 0 & -R\\ \id & 0\end{pmatrix}}_2=\max(\norm\id_2,\norm{R}_2)=\max(1,\norm{R}_2).\]
   
   Since $R^T=R$  we can diagonalise $R_{jk}=\langle R(E_0,E_j)E_0,E_k\rangle$ by $D=Q R Q^T$ with $Q$ orthogonal and $D=\diag(\lambda_1,\dotsc,\lambda_n)$. Set $F_i\coloneqq \sum_j Q_i^j E_j$. Then,
   \begin{align*}
       \langle R(E_0,F_i)E_0,F_j\rangle &=\sum_{k,\ell} Q_i^k \langle R(E_0,E_k)E_0,E_\ell\rangle Q^\ell_j=\\
       &=\left( QRQ^T\right)_{ij}=\delta_{ij}\lambda_i,
   \end{align*}
   which implies that
   \be \lambda_i=\langle R(E_0,F_i)E_0,F_i\rangle=K(\langle E_0,F_i\rangle)\ee
   for $K$ the sectional curvature.

   Therefore, we get that $\norm{R}_2=\norm D_2=\max_{i=1,\dotsc,n}\abs{K(\langle E_0,F_i\rangle)}\leq b$ by our curvature bound $\abs{K} \leq b$, and consequently
   \[ \max(1,\norm{R}_2) \leq \max(1,b) = \hat{b}. \]
   From (\ref{eq:ODE}) we obtain
   \be \norm A_2\leq \norm{(A(t) \;\; A'(t))}_2\leq \norm{(I \;\; \bar N(x))}_2 e^{\hat{b}t}\leq (1+\norm{\bar N(x)}_2)e^{\hat{b}t}.\ee
   
   Hence, 
   \begin{align*}
       \abs{\rho(t,x)} &=\omega(d\phi(\partial_t),d\phi(v_1),\dotsc,d\phi(v_n))\\
       &=\omega(E_0,J_1,\dotsc,J_n)\\
       &=\det\begin{pmatrix}1& \langle J_1,E_0 \rangle & \dotsc &\langle J_n, E_0\rangle\\0& \langle J_1,E_1 \rangle & \dotsc &\langle J_n, E_1\rangle\\
       \vdots & \vdots  & \vdots &\vdots \\
       0& \langle J_1,E_n \rangle & \dotsc &\langle J_n, E_n\rangle\\
       \end{pmatrix}\\
       &=\det(\langle J_i,E_j\rangle_{1\leq i,j\leq n})\\
       &=\underbrace{\det(A_{ij})}_{\text{polynomial of degree} \dim S}\\
       &\leq C\norm A_2^{\dim S}\\
       &\leq C (1+\norm{\bar N(x)}_2)^n\cdot e^{n \hat{b} t},\\
   \end{align*}
   for some constant $C>0$.
\end{proof}

\bibliographystyle{hamsalpha}
\bibliography{main}

\end{document}